\documentclass[12pt]{iopart}
 \expandafter\let\csname equation*\endcsname\relax
  \expandafter\let\csname endequation*\endcsname\relax
\usepackage{amsmath}
\usepackage{amsthm}
\usepackage{url}
\usepackage{float}
\usepackage{graphicx}
\usepackage{epsfig,subfigure}
\usepackage{graphics}
\usepackage{footnote}

\newtheorem{lemma}{Lemma}
\newtheorem{remark}{Remark}
\newcommand{\argmin}[1]{\textnormal{arg} \min_{#1}}
\newcommand{\argmax}[1]{\textnormal{arg} \max_{#1}}
\newcommand{\alphaopt}{\alpha_{\mathrm{opt}}}
\newcommand{\bfeta}{\mbox{\boldmath{$\eta$}}}
\newcommand{\bfEps}{\mbox{\boldmath{$\mathcal{E}$}}}
\newcommand{\bfeps}{\mbox{\boldmath{$\epsilon$}}}
\newcommand{\g}{\gamma}
\newcommand{\zetaopt}{\zeta_{\mathrm{opt}}}
\newcommand{\tildeA}{\tilde{A}}
\newcommand{\bfb}{\mathbf{b}}
\newcommand{\bfbtrue}{\mathbf{b}_{\mathrm{ex}}}
\newcommand{\bfbt}{\tilde{\mathbf{b}}}
\newcommand{\Cb}{C_{\bfeta}}
\newcommand{\Ceta}{C_{\bfeta}}
\newcommand{\bfe}{{\mathbf e}}
\newcommand{\bff}{\mathbf{f}}
\newcommand{\G}{\mathcal{G}}
\newcommand{\bfpfull}{\mathbf{p}_{\mathrm{full}}}
\newcommand{\bfr}{\mathbf{r}}
\newcommand{\bfrfull}{\mathbf{r}_{\mathrm{full}}}
\newcommand{\bfrproj}{\mathbf{r}_{\mathrm{proj}}}
\newcommand{\tmin}{t_{\mathrm{min}}}
\newcommand{\tmax}{t_{\mathrm{max}}}
\newcommand{\topt}{t_{\mathrm{opt}}}
\newcommand{\toptmin}{t_{\mathrm{opt}-\mathrm{min}}}
\newcommand{\toptp}{t_{\mathrm{opt}-\rho}}
\newcommand{\toptG}{t_{\mathrm{opt}-\G}}
\newcommand{\bfu}{\mathbf{u}}
\newcommand{\bfv}{\mathbf{v}}
\newcommand{\bfx}{\mathbf{x}}
\newcommand{\bfxtrue}{\mathbf{x}_{\mathrm{ex}}}
\newcommand{\bfxapr}{\mathbf{x}_{\mathrm{apr}}}
\newcommand{\bfxalpha}{\bfx(\alpha)}
\newcommand{\bfy}{\mathbf{y}}
\newcommand{\bfxyalpha}{\bfy}
\newcommand{\bfyalpha}{\bfy(\alpha)}
\newcommand{\bfw}{\mathbf{w}}
\newcommand{\Wb}{W_{\bfeta}}
\newcommand{\bfz}{\mathbf{z}}
\newcommand{\bfzalpha}{\bfz(\alpha)}

\begin{document}
\title[Iterated Lanczos  hybrid regularization] {Hybrid and iteratively reweighted regularization  by unbiased predictive risk and weighted GCV for projected systems \\ \today}
\author{Rosemary A Renaut$^1$,  Saeed Vatankhah$^2$, and Vahid E Ardestani $^2$}
\date{\today}
\address{$^1$ School of Mathematical and Statistical Sciences, Arizona State University, Tempe, USA. \\$^2$Institute of Geophysics, University of  Tehran,Tehran, Iran. }
\eads{\mailto{renaut@asu.edu}, \mailto{svatan@ut.ac.ir},  \mailto{ebrahimz@ut.ac.ir} }  


\begin{abstract}
Tikhonov regularization for projected solutions of large-scale ill-posed problems is considered. The Golub-Kahan iterative bidiagonalization is used to project the problem onto a subspace and regularization then applied to find a subspace approximation to the full problem. Determination of the regularization parameter for the projected problem by  unbiased predictive risk estimation,  generalized cross validation and discrepancy principle techniques is investigated.  It is shown that the regularized parameter obtained  by the unbiased predictive risk estimator can provide a good estimate for that to be used for a full problem which is moderately to severely ill-posed. A similar analysis provides the weight parameter for the weighted generalized cross validation such that the approach is also useful in these cases, and  also explains why the generalized cross validation without weighting is not always useful. All results  are independent of whether systems are over or underdetermined. Numerical simulations for standard one dimensional test problems and two dimensional data, for both image restoration and tomographic image reconstruction, support the analysis and validate the techniques.  The size of the projected problem is found using an extension of a noise revealing function for the projected problem Hn\u etynkov\'a,    Ple\u singer, and Strako\u s, [\textit{BIT Numerical Mathematics} {\bf 49} (2009), 4 pp. 669-696.].  Furthermore, an iteratively reweighted regularization approach for edge preserving regularization is extended for projected systems, providing stabilization of the solutions of the projected systems and reducing dependence on  the determination of the size of the projected subspace. 
\end{abstract}

\maketitle


\section{Introduction}
The solution of the numerically ill-posed linear system of equations 
\begin{align}\label{linear equation}
\bfb = A\bfxtrue +\bfeta, \quad \bfb \in\mathcal{R}^{m}, \quad \bfx \in \mathcal{R}^{n}, \quad A \in \mathcal{R}^{m \times n},
\end{align}
for  matrix $A$ of large dimension with $m\ge n$ or $m< n$ is considered.  Matrix $A$  is  ill-conditioned;  the singular values of $A$ decay exponentially to zero, or to the limits of the numerical precision.  Noise in the data is represented by  $\bfeta \in \mathcal{R}^{m}$, i.e. $\bfb=\bfbtrue+\bfeta$ for exact but unknown data $\bfbtrue$ that satisfies $\bfbtrue=A\bfxtrue$ for unknown exact model parameters $\bfxtrue$. Components  $\eta_i$ of $\bfeta$ are assumed to be independently sampled from a Gaussian distribution with  mean $0$ and covariance  $s_i^2$.  Given $A$ and $\bfb$ an  estimate for $\bfx$ that predicts $\bfbtrue$ is desired. 

Discrete ill-posed problems of the form \eqref{linear equation} may be obtained by discretizing linear ill-posed problems such as Fredholm integral equations of the first kind and arise in many research areas including image deblurring, geophysics, etc.  Due to the presence of the noise in the data and the  ill-conditioning of $A$  regularization is needed in order to obtain an estimate for $\bfx$ approximating $\bfxtrue$. Standard Tikhonov regularization provides 
\begin{align}\label{tikhonov1}
\bfx(\alpha)=\argmin{\bfx}{\lbrace\|\Wb(A\bfx-\bfb)\|_2^2 + \alpha^2 \|D(\bfx-\bfxapr)\|_2^2\rbrace},
\end{align}
for weighted data fidelity  term $ \|\Wb(A\bfx-\bfb)\|_2^2$  and regularization term $\|D(\bfx-\bfxapr)\|_2^2$. $D$ is a regularization matrix, assumed here to be invertible, and $\bfxapr$ allows specification of a given reference vector of  prior information for $\bfx$.  The unknown regularization parameter $\alpha$   trades-off between the data fidelity and regularization terms. The noise in the measurements ${\bfb}$ is whitened  when $\Wb=\Cb^{-1/2}$ for the covariance matrix $\Cb=\mathrm{diag}(s_1^2,\dots,s_{m}^2)$.   Introducing  $\bfbt=\Wb \bfb$,  $\tildeA=\Wb A$,  shifting by the prior information through $ \bfxyalpha=\bfx-\bfxapr$, and  assuming that the null spaces of $\tildeA$ and $D$ do not intersect, yields
\begin{align}\label{tikhonov2}
\bfyalpha&=\argmin{\bfy}\lbrace\|\tildeA \bfxyalpha-\tilde{\bfr}\|_2^2 + \alpha^2 \|D\bfxyalpha\|_2^2\rbrace, \quad \tilde{\bfr}=(\bfbt-\tildeA\bfxapr)\\ \nonumber
&=(\tildeA^T\tildeA+\alpha^2D^TD)^{-1}\tildeA^T \tilde{\bfr}.
\end{align}
Analytically when $D$ is invertible, which is not always the case,  we may write
\begin{align*}(\tildeA^T\tildeA+\alpha^2 D^TD ) = D^T((D^T)^{-1}\tildeA^T\tildeA D^{-1}+\alpha^2 I_n )D.\end{align*}
Thus when it is feasible to calculate $D^{-1}$, or to solve systems of equations defined by invertible $D$, it is convenient to introduce the  right preconditioned matrix $\tilde{\tildeA}=\tildeA D^{-1}$ and regularized inverse $\tilde{\tildeA}^\dagger(\alpha)=(  \tilde{\tildeA}^T\tilde{\tildeA} +\alpha^2 I_n)^{-1}   \tilde{\tildeA}^T$,\footnote{Note that we use in general the notation $A^\dagger(\alpha)$ for the pseudo inverse of the augmented matrix $[A;\alpha I]$.} which provides
\begin{align}\label{tikhonov3}
\bfzalpha:&=\argmin{\bfz}\lbrace\|\tilde{\tildeA} \bfz - \tilde{\bfr} \|_2^2 + \alpha^2 \|\bfz\|_2^2\rbrace, \quad \bfzalpha=D\bfyalpha, \quad \mathrm{and} \\
\bfx(\alpha)&=\bfxapr+\bfyalpha = \bfxapr+ D^{-1}\tilde{\tildeA}^\dagger(\alpha) \tilde{\bfr}.   \label{xupdate}
\end{align}
Although equivalent analytically, numerical techniques to solve \eqref{tikhonov2} and \eqref{tikhonov3} differ. For small scale problems, for example, we may solve \eqref{tikhonov2}  using the    generalized singular value decomposition (GSVD), e.g. \cite{PaigeSau1},  for the matrix pair $[\tildeA, D]$, but would use the singular value decomposition (SVD) of $\tilde{\tildeA}$ for  \eqref{tikhonov3}, e.g. \cite{GoLo:96}, as given in Appendix~\ref{SVD}, dependent on the feasibility of calculating $D^{-1}$.  Still, the use of the SVD or GSVD  is not viable computationally for large scale problems unless the underlying operators possess a specific structure. For example, if the underlying system matrix, and associated regularization matrix are expressible via   Kronecker decompositions, e.g. \cite{HNO}, then the GSVD decomposition can be found via the GSVD for each dimension separately. Here we consider the  general situation  and use of  iterative Krylov methods to estimate   $\bfx(\alpha)$.

\subsection{Numerical solution by the Golub-Kahan bidiagonalization}
In principle iterative methods such as conjugate gradients or other Krylov methods, can be employed to solve \eqref{tikhonov1}.  Results presented in \cite{HaJe}  demonstrate, however, that MINRES and GMRES should not be used as regularizing Krylov iterations due to the early transfer of noise to the Krylov basis. Here we use the   well-known Golub-Kahan bidiagonalization (GKB), implemented in the LSQR algorithm, which has been well-studied in the context of projected solutions of the  least squares problem \cite{PaigeSau2,PaigeSau3}.   Recently, there has also been some interest in the LSMR modification of LSQR, \cite{FoSa}, but due to our goal to investigate the regularization parameter $\alpha$ we focus on LSQR for which the noise regularizing properties of the iteration are better understood, \cite{HaJe, HPS}. Effectively the GKB  projects the solution of the inverse problem to a smaller subspace, say of size $t$.

Applying $t$ steps of the GKB on  matrix $A$ 
with initial vector $ \bfb $, of norm $\beta_1=\|\bfb\|_2$, and defining $\bfe^{(t+1)}_1$ to be the unit vector of length $t+1$ with a $1$ in the first entry, lower bidiagonal matrix $B_t \in\mathcal{R}^{(t+1) \times t} $ and column orthonormal matrices $H_{t+1} \in\mathcal{R}^{m \times (t+1)} $ , $G_t \in\mathcal{R}^{n \times t} $ are generated such that, see \cite{Regtools,kilmer:2001}, 
\begin{align}\label{bidiag}
A G_t=H_{t+1}B_t, \quad \beta_1 H_{t+1}\bfe^{(t+1)}_1={\bfb}.
\end{align}
For $\bfx_t=G_t \bfw_t$,  full, $\bfrfull(\bfx_t)$, and projected,  $\bfrproj(\bfw_t)$, residuals are related via
\begin{align}\label{fullresidual}
\bfrfull(\bfx_t)  =A \bfx_t -\bfb &= A G_t \bfw_t - \beta_1 H_{t+1} \bfe^{(t+1)}_1  = \\
H_{t+1} B_t \bfw_t - \beta_1 H_{t+1} \bfe^{(t+1)}_1&=H_{t+1}(B_t \bfw_t - \beta_1  \bfe^{(t+1)}_1) =H_{t+1}\bfrproj(\bfw_t), \nonumber 
\end{align}
for which,  by the column orthonormality of $H_{t+1}$,
\begin{align}\label{normres}
 \|\bfrfull(\bfx_t)\|_2^2 = \|\bfrproj(\bfw_t)\|_2^2.
\end{align}
Theoretically, therefore, an estimate for  $\bfx$ with respect to a reduced subspace may be found by finding $\bfw_t$ and then projecting back to the full problem. Matrix $B_t$ in most cases, however, inherits the ill-conditioning of the matrix $A$, \cite{PaigeSau2},  and regularization of the projected problem is needed.

By the column orthonormality of $G_t$, we have $\|\bfx_t\|_2^2 = \|G_t \bfw_t\|_2^2 = \|\bfw_t\|_2^2$.   Thus, explicitly introducing regularization parameter $\zeta$, distinct from $\alpha$ in order to emphasize regularization on the projected problem, yields the  projected Tikhonov problem
 \begin{align}\label{projectedtikhonov} 
\bfw_t(\zeta)=\argmin{\bfw \in \mathcal{R}^t}\lbrace\|B_t\bfw-\beta_1 \bfe^{(t+1)}_{1}\|_2^2 + \zeta^2 \|\bfw\|_2^2\rbrace,
\end{align}
with solution 
\begin{align}\label{regsoln}
\bfw_t(\zeta)&= \beta_1(B_t^TB_t +\zeta^2 I_t)^{-1}B_t^T \bfe^{(t+1)}_1 = \beta_1 B^\dagger_t(\zeta) \bfe^{(t+1)}_1  \\
&=(G_t^TA^T A G_t +\zeta^2 I_t)^{-1}G_t^TA^T \bfb = (AG_t)^\dagger(\zeta) \bfb.  \nonumber 
\end{align}
In practice,  one uses \eqref{regsoln} to find $\bfw_t(\zeta)$ via the SVD for $B_t$, under the assumption that $t << m^*=\min(m,n)$, noting that an explicit solution for $\bfw_t$ is immediately available, see e.g. Appendix~\ref{SVD}.

As already observed in \cite[p. 302]{Hanke}, the regularized LSQR algorithm now poses the problem of both detecting the appropriate number of steps $t$ as well as of finding the optimal parameter $\zetaopt$.  One method of regularization is simply to avoid the introduction of the regularizer in \eqref{projectedtikhonov} and find an optimal $t$ at which to stop the iteration. Although it is known that the LSQR iteration is a regularizing iteration, it also exhibits a semi-convergence behavior so that eventually regularization is also needed. This regularization may be achieved by either picking $\alpha$ in advance,  namely regularize and project, or by the hybrid approach of regularizing the projected problem,  e.g.  \cite{ChNaOl:08,kilmer:2001,RHM10}. The problem of first determining the appropriate size $t$ for the projected space is discussed in  e.g. \cite{HPS,kilmer:2001} and more recently for large scale geophysical inversion in  \cite{PaHaHaFe}.   Although the solutions obtained from the regularize then project, and project then regularize, for a given $t$ and $\alpha=\zeta$ are equivalent, \cite[Theorem 3.1]{kilmer:2001}, \cite[p 301]{Hanke}, this does not immediately mean that 
  $\zetaopt$ for the subspace problem  provides $\alphaopt$ for the full problem, \cite{kilmer:2001}.
\begin{remark}
Determining to which degree certain regularization techniques provide a good estimate for $\alphaopt$  from the subspace problem estimate, and the conditions under which this will hold,  is the topic of this work and is the reason we denote regularization parameter on the subspace by $\zeta$ distinct from $\alpha$. 
\end{remark}

\subsection{Regularization parameter estimation}
For the full problem the question of determining an optimal parameter $\alphaopt$ is well-studied, see e.g. \cite{Hansen:98,Vogel:2002}, for a discussion of methods including the Morozow discrepancy principle (MDP), the L-curve (LC), generalized cross validation (GCV)  and unbiased predictive risk estimation (UPRE). 
The use of the MDP, LC and GCV is also widely discussed for the projected problem, particularly starting with the work of Kilmer et al, \cite{kilmer:2001} and continued in \cite{ChNaOl:08}. Further, extensions for windowed regularization, and hence multi-parameter regularization,  \cite{ChEaOl:11} are also applied for the projected problem \cite{ChKiOl:15}. 
Our attention is initially on the use of the UPRE.  Effectively, the UPRE provides the correct estimate for $\alphaopt$ in the context of a filtered truncated SVD (FTSVD) solution of \eqref{tikhonov1} with $t$ terms, provided that the LSQR factorization effectively captures the dominant right singular subspace of size $t$ for matrix $A$. This  observation does not immediately extend to the GCV. Applying a similar analysis as for the UPRE, however, provides a  choice of the weighting parameter in the weighted GCV  (WGCV) introduced in \cite{ChNaOl:08}.  

We stress that the approach assumes throughout, both numerically and theoretically, that the projected system is calculated with full reorthogonalization, a point not made explicit in many discussions, although it is apparent than many references implicitly make this assumption. 
\subsection{Overview}
The paper is organized as follows. The regularization parameter estimation techniques of interest are presented in  \S\ref{sec:parameter estimation}. The discussion  in \S\ref{sec:parameter estimation} is validated with one dimensional simulations in \S\ref{sec:simulationoned}. Image restoration problems presented in \S\ref{sec:simulationtwod} illustrate the relevance for the two dimensional case. In \S\ref{sec:irr} we extend the hybrid  approach for use with an iteratively reweighted regularizer (IRR), which sharpens edges within the solution, \cite{PoZh:99,vatan:2014,VAR:2014b,VRA:2014,WoRo:07,Zhd:2002}, hence demonstrating that edge preserving regularization can be applied in the context of regularized LSQR solutions of the least squares problem on a projected subspace. 
 Finally in \S\ref{sec:walnut} we also illustrate the algorithms in the context of sparse tomographic reconstruction of a walnut data set, \cite{walnut1}, demonstrating the more general use of the approach beyond deblurring of noisy data. Our conclusions are presented in \S\ref{conclusions}. It is of particular interest that our analysis applies for both over and under determined systems of equations and is thus potentially of future use for other algorithms also in which alternative regularizers are imposed and also  require repeated Tikhonov solves at each step. Further, this work extends our analysis of the UPRE in the context of underdetermined but small scale problems in \cite{VAR:2014b,VRA:2014}, and demonstrates that IRR can be applied for projected algorithms.

\section{Regularization parameter estimation}\label{sec:parameter estimation}
In order to use any specific regularization parameter estimation method for the projected problem it is necessary  to understand the derivation on the full problem. We thus provide a brief overview of the derivations as needed. 

\subsection{Unbiased Predictive Risk Estimator}
The predictive error, $\bfpfull(\bfx(\alpha))$, for the solution $\bfxalpha$,     
is defined by 
\begin{align}\label{risk}
\bfpfull(\bfx(\alpha)) &= A \bfxalpha - \bfbtrue 
 =AA^\dagger(\alpha)\bfb  - \bfbtrue = (A(\alpha) -I_m)\bfbtrue +A(\alpha) \bfeta , \end{align}
where $A(\alpha)=AA^\dagger(\alpha)$ is the influence matrix.  The residual may also be written in terms of the $A(\alpha)$ as 
\begin{align}\label{measurable residual}
\bfrfull(\bfx(\alpha))=(A(\alpha)-I_m)\bfb  =(A(\alpha)-I_m)\bfbtrue+(A(\alpha)-I_m)\bfeta.
\end{align}
In both equations the  first term is deterministic, whereas the second is stochastic due to noise vector  ${\bfeta}$. To proceed we need the Trace Lemma e.g. \cite[Lemma 7.2]{Vogel:2002}.
\begin{lemma}\label{lemma21}
For deterministic vector $\bff$, random vector $\bfeta$ with diagonal covariance matrix $\Ceta$, matrix $F$, and expectation operator $E$
\begin{align}\nonumber
E(\|\bff +F \bfeta\|_2^2) = \|\bff\|_2^2 + \mathrm{tr}(\Ceta F^TF),
\end{align}
using $\mathrm{tr}(A)$ to denote the trace of  matrix $A$. 
\end{lemma}

Applying Lemma~\ref{lemma21} to  both \eqref{risk}  and \eqref{measurable residual} with the assumption that $\Ceta=I_m$, due to whitening of noise $\bfeta$, and using the symmetry of the influence matrix,  we obtain 
\begin{align}\label{risk1}
 E(\|\bfpfull(\bfx(\alpha))\|_2^2) &=   \|(A(\alpha) -I_m)\bfbtrue\|_2^2  +\mathrm{tr}(A^T(\alpha)A(\alpha)) \\
  E(\|\bfrfull(\bfx(\alpha))\|_2^2) &=   \|(A(\alpha) -I_m)\bfbtrue\|_2^2  +\mathrm{tr}((A(\alpha)-I_m)^T(A(\alpha)-I_m)).\label{risk2}
\end{align}
Here $E(\|\bfpfull(\bfx(\alpha))\|_2^2)/m$ is the expected value of the risk of using the solution $\bfxalpha$ to predict $\bfbtrue$.  
The first term on the right hand side  in each case  cannot be obtained, but we may use $E(\|\bfrfull(\bfx(\alpha))\|_2^2) \approx \|\bfrfull(\bfx(\alpha))\|_2^2$ in \eqref{risk2}. Thus using linearity of the trace and  eliminating the first term in the right hand side of   \eqref{risk1} gives the UPRE estimator to find $\alphaopt$
\begin{align}\label{optalpha2}
\alphaopt=\argmin{\alpha}\{U(\alpha)= \|(A(\alpha)-I_m)\bfb \|_2^2 +2\,\mathrm{tr}(A(\alpha)) -m\}.
\end{align}
Typically, $\alphaopt$   is found by evaluating \eqref{optalpha2} for a range of $\alpha$, for example by the SVD see e.g. Appendix~\ref{appB}, with the minimum found within that range of parameter values, as suggested in \cite{Regtools} for the GCV.  See also e.g. \cite[Appendix, (A.6)]{VRA:2014} for the formulae for calculating the function in terms of the SVD of matrix $A$.  

\subsubsection{Extending the UPRE for the projected problem}
We observe that we may immediately write the predictive error and the residual in terms of the solution of the projected problem explicitly depending on the  regularization parameter $\zeta$. Specifically,   defining the influence matrix $(AG_t)(\zeta)=  A G_t(AG_t)^\dagger(\zeta)$ for the projected solution we have 
\begin{align}\label{fprisk}
\bfpfull(\bfx_t(\zeta)) &=  A G_t  \bfw_t(\zeta)- \bfbtrue=  (AG_t)(\zeta)\bfb - \bfbtrue \\ \label{eq:projres}
\bfrfull(\bfx_t(\zeta))&= A G_t  \bfw_t(\zeta)- \bfb  
 =  \left((AG_t)(\zeta)-I_m\right) \bfb.
\end{align}
By comparing \eqref{fprisk} with \eqref{risk} and \eqref{eq:projres} with \eqref{measurable residual} we obtain 
\begin{align*}
U_{\mathrm{full}}(\zeta )= \|\left((AG_t)(\zeta )-I_m\right)\bfb \|_2^2 + 2 \,\mathrm{tr}\left((AG_t)(\zeta )\right) -m.\end{align*}
Now by  \eqref{normres} it is immediate that the first term can be obtained without finding $\bfx_t(\zeta)$. For the second term we observe
\begin{align*}
(AG_t)(\zeta ) &= AG_t ((AG_t)^TAG_t +\zeta^2 I_t)^{-1}(AG_t)^T \\ 
&= H_{t+1}B_t ((H_{t+1}B_t)^T(H_{t+1}B_t) +\zeta^2I_t)^{-1} (H_{t+1}B_t)^T \\
&=H_{t+1} \left(B_t(B_t^TB_t+\zeta^2 I_t)^{-1}B^T_t\right)H^T_{t+1}=H_{t+1}B_t(\zeta)H_{t+1}^T,
\end{align*}
which yields
\begin{align*}\nonumber
\Tr\left((AG_t)(\zeta )\right)=\Tr(H_{t+1} B_t(\zeta)H^T_{t+1}) =\Tr(B_t(\zeta)),
\end{align*}
where the last equality follows from the cycle property of the trace operator for consistently sized matrices. Hence 
\begin{align}\label{projfull}
U_{\mathrm{full}}(\zeta ) = \|\beta_1(B_t(\zeta) -I_{t+1})\bfe_1^{t+1}\|_2^2 + 2 \, \Tr(B_t(\zeta)) -m
\end{align}
can be evaluated without  reprojecting the solution for every $\zeta$ back to the full problem.
\begin{remark}
Although the UPRE function can be found from the projected solution alone, it is not clear whether 
\eqref{projfull} has any relevance with respect to the projected solution, i.e. does this appropriately regularize the projected solution,  otherwise it may not be appropriate to find $\zeta$ to minimize this function on the subspace. 
 \end{remark}
 
The projected solution solves the problem with system matrix $B_t$ and right hand side vector $\beta_1 \bfe_1^{t+1}=H_{t+1}^T\bfb$ which also  consists of a deterministic and stochastic part, $ H_{t+1}^T\bfbtrue+H_{t+1}^T\bfeta$, where for white noise vector $\bfeta$ and column orthogonal $H_{t+1}$, $H_{t+1}^T\bfeta$ is a random vector of length $t+1$ with covariance matrix $I_{t+1}$. Thus the UPRE for the projected problem is  
\begin{align}\label{projupre}
U_{\mathrm{proj}}(\zeta) = \|\beta_1(B_t(\zeta)-I_{t+1})\bfe_1^{(t+1)}\|_2^2 + 2 \, \Tr(B_t(\zeta))-(t+1).
\end{align}
Comparing \eqref{projfull} with \eqref{projupre},  it is immediate that minimizing \eqref{projfull} to minimize the risk for the projected solution, also minimizes the risk for the full solution with respect to the given subspace. 

It remains to determine whether there is any case in which finding $\zetaopt$ also minimizes the predictive risk (\ref{optalpha2}) for the full problem. Specifically it is not immediate that  $U_{\mathrm{full}}(\alphaopt)  \approx U_{\mathrm{full}}(\zetaopt)$ because $\alphaopt$ is needed with respect to solutions in $\mathrm{Range}(V)$, not just restricted to $\mathrm{Range}(G_t\tilde{V}_t )$. Here matrices $V$ and $\tilde{V}_t$, the column orthogonal matrices arising in the SVDs of $A$ and $B_t$, respectively, span the respective right singular subspaces.  Although in exact arithmetic the large singular values of $B_t$ provide a good approximation of  the large singular values of $A$,  \cite[Section 9.3.3]{GoLo:96}, the number of small singular values in the spectrum of $B_t$ limits how well the full problem will be regularized by regularizing the projected problem.  Adopting now the statement of \textit{full regularization} of the LSQR as given in \cite{HJ:16}, namely that the LSQR iterate with $t$ steps effectively captures the $t$-dimensional dominant right spectral space of $A$, suppose that $t$ is such that the singular values of $B_t$ approximate the $t$ largest singular values of $A$ with the natural order so that necessarily $\gamma_t > \sigma_{t^*+1}$ for $t\le t^*$.   Equivalently this requires that  $t^*$ is close to $t$ and that the spectrum of $B_t$ contains no singular value approximating a very small spectral value of $A$. It is shown in \cite[Theorem 2.3]{HJ:16} that this requirement is more likely satisfied for severely and moderately ill-posed problems, than for mildly ill-posed problems. Further, in such cases the LSQR solution on the space of size $t$ approximates the TSVD  solution of the full problem, namely the solution of the full problem with filter factors $\phi_i(\alpha)=0$ for $i>t$. 
Then, 
\begin{align} \nonumber
\Tr(A(\alpha)) 
 &= \left(t - \alpha^2 \sum_{i=1}^t (\sigma_i^2 +\alpha^2)^{-1}\right) +\left((m^*-t) - \alpha^2 \sum_{i=t+1}^{m^*} (\sigma_i^2 +\alpha^2)^{-1} \right) \\
& = \left(t - \alpha^2 \sum_{i=1}^t (\sigma_i^2 +\alpha^2)^{-1}\right) \approx  t -\alpha^2 \sum_{i=1}^t (\g_i^2 + \alpha^2)^{-1} =\Tr(B_t(\alpha)). \label{tracerelation}
\end{align}
Thus,  in the situation in which the LSQR iterate provides full regularization, determining $\zetaopt$ to minimize \eqref{projupre} will yield $\alphaopt$ which is optimal for the  filtered truncated SVD (FTSVD) solution of the full problem.
This observation also follows Theorem  3.2 \cite{kilmer:2001} which connects the use of the TSVD of $B_t$ for the solution with the solution obtained using the TSVD of $A$. To summarize:
\begin{remark}[UPRE]
If $t$ is such that $\phi_i(\alpha) \approx 0$ for $i>t$, and such that the LSQR iterate provides \textit{full regularization}  so that $\Tr(A(\alpha)) \approx  \Tr(B_t(\alpha))$ and $\mathrm{Range}(G_t \tilde{V}_t)$  approximates $\mathrm{Range}(V_t)$,  then $\zetaopt \approx \alphaopt$ when obtained using the UPRE.  Further, the estimate is found without projecting the solution back to the full space, namely by minimizing \eqref{projupre}.
\end{remark} 

When the LSQR does not provide \textit{full regularization}, denoted as \textit{partial regularization} in \cite{HJ:16}, the above result will not hold, and  $B_t$ captures the ill-conditioning of $A$ through the inclusion of  inaccurate small singular values in the spectrum of $B_t$.

\subsection{Morozov discrepancy principle}
Although it is well-known that the MDP always leads to an over estimation of the regularization parameter,  e.g. \cite{kilmer:2001}, it is still a widely used method for many applications, and is thus an important baseline for comparison. 
The premise of the  MDP,  \cite{morozov},     to find $\alpha$ is  the assumption that the norm of the residual, $\|\bfrfull(\bfx(\alpha))\|_2^2$ follows a $\chi^2$ distribution with $\delta$ degrees of freedom, 
 $\|\bfrfull(\bfx(\alpha))\|_2^2=\delta$. Heuristically, the rationale for this choice is seen by re-expressing \eqref{fullresidual} 
 \begin{align*}
 \bfrfull(\bfx(\alpha))=A\bfx(\alpha)-\bfb=A(\bfx(\alpha)-\bfxtrue)-\bfeta ,
 \end{align*}
 so that if $\bfx(\alpha)$ has been found as a good estimate for $\bfxtrue$, then the residual \eqref{fullresidual} should be dominated by the whitened error vector $\bfeta $. For white noise $\|\bfeta\|_2^2$ is distributed as a $\chi^2$ distribution with $m$ degrees of freedom, from which  $E(\|\bfeta\|^2)=m$,  with variance $2m$.  Thus we  seek a residual such that $\delta=\upsilon m$   using a Newton root-finding method, see Appendix~\ref{appB}, where we take safety parameter $\upsilon>1$ to handle the well-known over smoothing of the MDP.   Applying the same approach for the projected residual yields the noise term  
$H_{t+1}^T\bfeta$ replacing $\bfeta$. Thus the degrees of freedrom are reduced to   $t+1$ and we seek a residual such that $\delta=\upsilon (t+1)$. 
A number of other suggestions for a projected discrepancy principle  have been presented in the literature, but generally imply using $\delta \approx \upsilon \|\bfeta\|_2^2\approx \upsilon E(\|\bfeta\|_2^2) \approx \upsilon m$ dependent on the noise level of the full problem, e.g. \cite{Hanke,kilmer:2001,ReSgYe}, with $\upsilon>1$.   It is reported in \cite{Hanke}, however, that while the theory predicts choosing $\upsilon>1$, numerical experiments support reducing $\upsilon$. Alternatively this may be seen as reducing the degrees of freedom, instead of reducing $\upsilon$.  We deduce that the interpretation for finding the regularization parameter based on the statistical property of the projected residual in contrast to the full residual should be important in determining $\delta$.  
 
 \begin{remark}[MDP]
For the MDP  the degrees of freedom change from $m$ to $t+1$  when the residual is calculated on the full space as compared to the projected space. Thus $\zetaopt $ is not a good approximation for $ \alphaopt$ when obtained using $\delta_{\mathrm{proj}}$  as a guide for the actual size of the projected residual.  For the \textit{full regularization} the degrees of freedom for the full problem are reduced and again $\zetaopt\approx \alphaopt$. 
 \end{remark}

\subsection{Generalized cross validation}
Unlike the UPRE and MDP, the GCV method for finding the regularization parameter $\alpha$ does not require any information on the noise distribution for $\bfeta$. The optimal parameter $\alpha$ is found as the minimizer of the function 
\begin{align}\label{fullgcv}
G_{\mathrm{full}}(\alpha)= \frac{\|\bfrfull(\bfx(\alpha)) \|_2^2 }{\left(\mathrm{tr}(A(\alpha)-I_m)\right)^2},
\end{align}
ignoring constant scaling of $G_{\mathrm{full}}(\alpha)$ by  $n$, \cite{GoHeWa}. 
The obvious implementation of the GCV for the projected problem is the exact replacement in \eqref{fullgcv} using the projected system
\begin{align}\label{projgcv}
G_{\mathrm{proj}}(\zeta)= \frac{\|\bfrproj(\bfw_t(\zeta)) \|_2^2 }{\left(\mathrm{tr}(B_t(\zeta)-I_{t+1})\right)^2},
\end{align}
as indicated in \cite{kilmer:2001}.  It was recognized in \cite[Section 5.4]{ChNaOl:08}, however, that this formulation, tends to lead to solutions which are over smoothed and as an alternative the WGCV  was introduced,  dependent on  parameter $\omega$,
\begin{align}\nonumber
G_{\mathrm{proj}}(\zeta,\omega)= \frac{\|\bfrproj(\bfw_t(\zeta)) \|_2^2 }{\left(\mathrm{tr}(\omega B_t(\zeta)-I_{t+1})\right)^2}.
\end{align}
Experiments illustrated that $\omega$ should be smaller for high noise cases, but in all cases $0\le \omega \le 1$ is required to avoid the potential of a zero in the denominator. The choice for $\omega$ was argued heuristically and an adaptive algorithm to find $\omega$ was given. 

Consider now the two denominators in \eqref{fullgcv} and \eqref{projgcv}. First of all by \eqref{tracerelation} it is not difficult to show, for $0<\omega<1$,  that  
\begin{align}\nonumber 0>\Tr( B_t(\alpha)-I_{t+1}) > \Tr(A(\alpha)-I_m),
\end{align}
so that $G_{\mathrm{proj}}(\alpha) > G_{\mathrm{full}}(\alpha)$ and  $\alpha$ chosen to minimize the projected GCV will not minimize the full GCV term. For the weighted GCV, however, 
\begin{align} \label{wgcvden}
\Tr(I_{t+1}-\omega B_t(\zeta))&=(1 +t   -\omega t )) +\omega \zeta^2\sum_{i=1}^t\frac{1}{\g^2_i+\zeta^2},\end{align}
and for the full regularization in which we approximate the full TSVD, $\phi_i(\alpha) \approx 0$ for $i>t$,  
\begin{align}\label{gcvden}
\Tr(I_m-A(\alpha))&=(m-m^*)+\alpha^2\sum_{i=1}^{m^*}\frac{1}{\sigma^2_i+\alpha^2}\approx (m-t) +\alpha^2\sum_{i=1}^{t}\frac{1}{\sigma^2_i+\alpha^2}.
\end{align}
 Factoring for $m-t \ne 0$ and $(t+1-\omega t)\ne 0$  in \eqref{gcvden} and \eqref{wgcvden}, respectively, gives the scaled denominators  
\begin{align*}
(m-t) \left(1+\frac{1}{m-t} \sum_{i=1}^t \frac{\alpha^2}{\sigma_i^2+\alpha^2}\right) \quad \text{and} \quad 
(1+t-\omega t) \left(1+\frac{\omega}{1+t-\omega t} \sum_{i=1}^t \frac{\zeta^2}{\g_i^2+\zeta^2}\right).
\end{align*}
Ignoring constant  scaling the denominators  are equilibrated by taking 
\begin{align*}
\frac{1}{m-t}  = \frac{\omega}{1+t-\omega t} \quad \text{yielding} \quad \omega=\frac{1+t}{m}<1. 
\end{align*}
  This result suggests  that we need $(t+1)/m \le \omega \le 1$ in order for $\zetaopt$ to   estimate $\alphaopt$ found with respect to the projected space. 
    \begin{remark}[GCV] 
  Taking $\omega=(t+1)/m$ gives the regularization of the FTSVD of the full problem, in the case of the LSQR with full regularization, namely for moderately and severely ill-posed problems as defined in \cite{HJ:16}.
  \end{remark}

\section{Simulations: One Dimensional Problems}\label{sec:simulationoned}
To illustrate the discussion in \S\ref{sec:parameter estimation} we examine the solution of ill-posed one dimensional problems with known solutions.  In all experiments we use MATLAB\textregistered 2014b  and test problems \texttt{phillips}  and \texttt{gravity} which are discretizations of  Fredholm integral equations of the first kind provided in the Regularization toolbox, \cite{Regtools}. Problem \texttt{gravity}  depends on a parameter $d$  determining the conditioning of the problem, here we use $d=0.75$ yielding a severely ill-posed test problem. In contrast problem \texttt{phillips} is moderately ill-posed and the Picard condition does not hold.\footnote{The continuous and discrete Picard condition are well-described in the literature, e.g. \cite{Hansen:98}. Basically the Picard condition holds if the absolute values of the coefficients of the solution decay on average faster than the singular values.} Simulations for over and under sampled data are  obtained by  straightforward modification of the relevant functions in \cite{Regtools}.  We  discuss representative results obtained for the undersampled case with $m=152$ and $n=304$,  for which the condition number of $A$ is   $4.05e+05$ and $3.38e+17$, for \texttt{phillips} and \texttt{gravity}, respectively. The function \texttt{bidiag\_gk} associated with the software for the paper \cite{HPS}, is used for finding the  factorization \eqref{bidiag} with full reorthogonalization against of all basis vectors (the default). 

\begin{figure}[!htb]
\begin{center}
\subfigure[\texttt{phillips}]{\includegraphics[width=.40\textwidth]{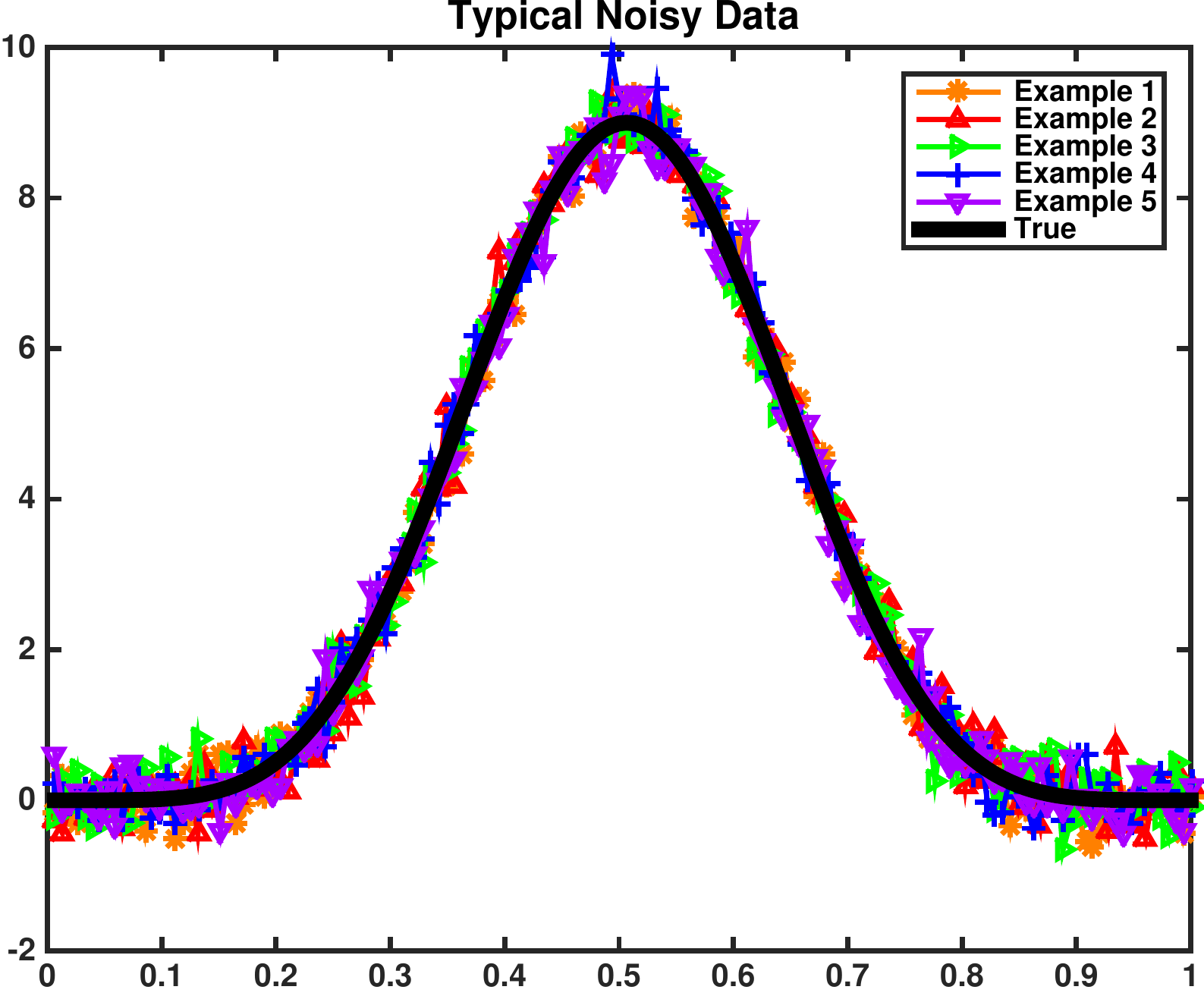}}\hspace{1cm}
\subfigure[\texttt{gravity}, $d=.75$]{\includegraphics[width=.40\textwidth]{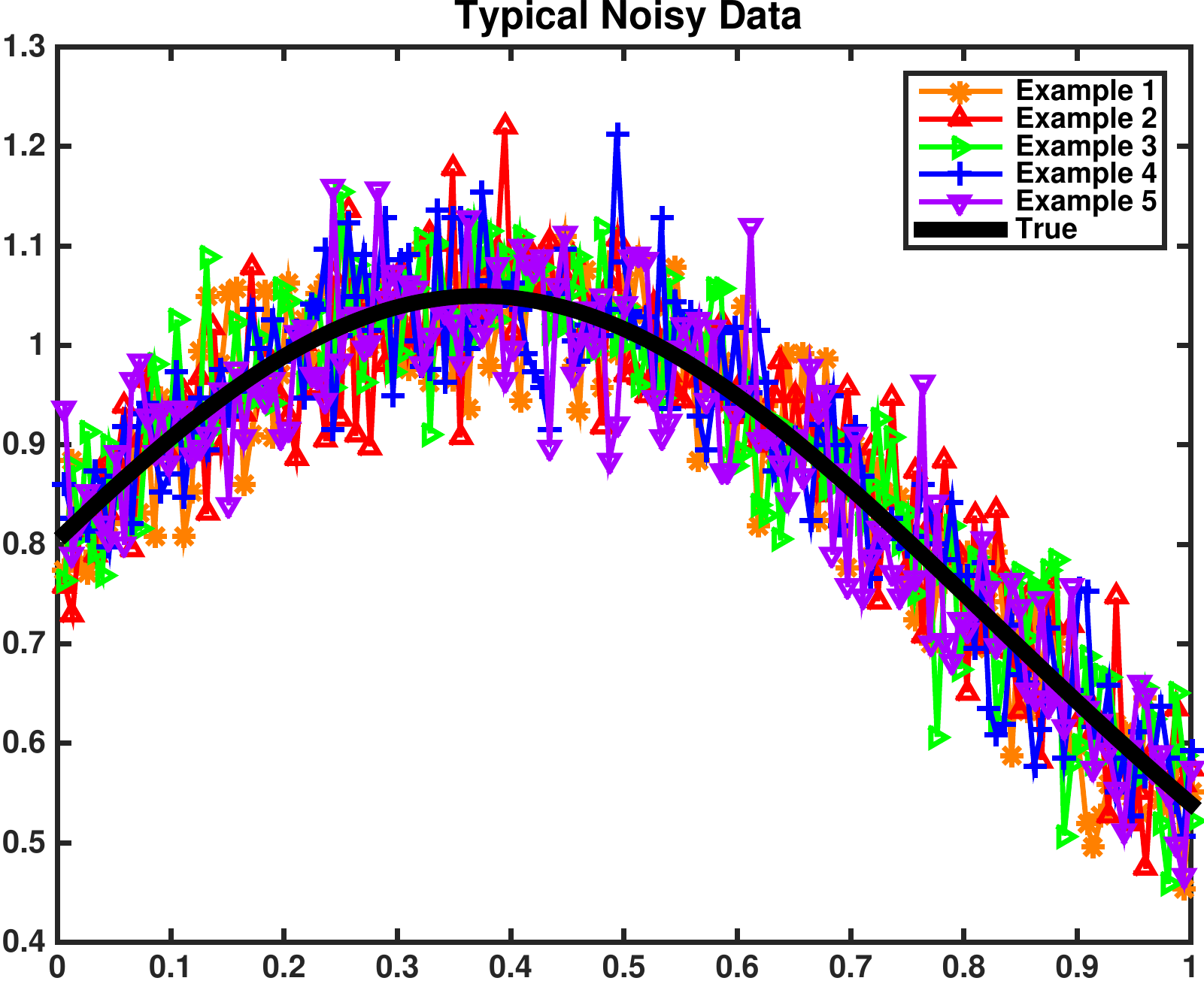}}
\end{center}\caption{Illustrative test data  for noise level $\eta=.005$ for sample right hand side data $\bfb^c$, $c=1\colon5$, with $m=152$ and $n=304$. Exact solutions are given by the solid lines in each plot. \label{fig:data}}
\end{figure}

For a given problem defined by \eqref{linear equation} without noise, noisy data are obtained as 
\begin{equation}\label{noiselevel}\bfb^c=\bfbtrue+\bfeta^c = \bfbtrue+ \eta \|\bfbtrue\|_2 \bfeps^c,
\end{equation}
for noise level $\eta$ and with $\bfeps^c$ the $c^{\text{th}}$ column of error matrix $\bfEps$ that has columns sampled from a random normal distribution using the MATLAB function \texttt{randn}$(m,nc)$.     The signal to noise ratio for the data  given by 
\begin{align}\label{BSNR}
\text{BSNR}(\eta,m)=20\, \mathrm{log}10\left(\frac{\|\bfbtrue\|}{\|\bfb^c-\bfbtrue\|}\right) \approx -20 \, \mathrm{log}10 (\eta \sqrt{m}),
\end{align}
is independent of the test problem. In particular  $\text{BSNR}(.005,152)$ $\approx 24.2$. Example simulation data are shown in Figure~\ref{fig:data} for noise levels $\eta=.005$ for each test problem with $m=152$ and in each case for $5$ samples of the noise, $\bfb^c$, $c=1\colon5$. In all simulations  the matrices and right hand side data are weighted by the diagonal inverse square root of the covariance matrix, assuming colored noise.

 \subsection{Spectra of $A$ and $B_t$}\label{spectra}
\begin{figure}[!htb]
\begin{center}
\subfigure[\texttt{phillips}]{\includegraphics[width=.4\textwidth]{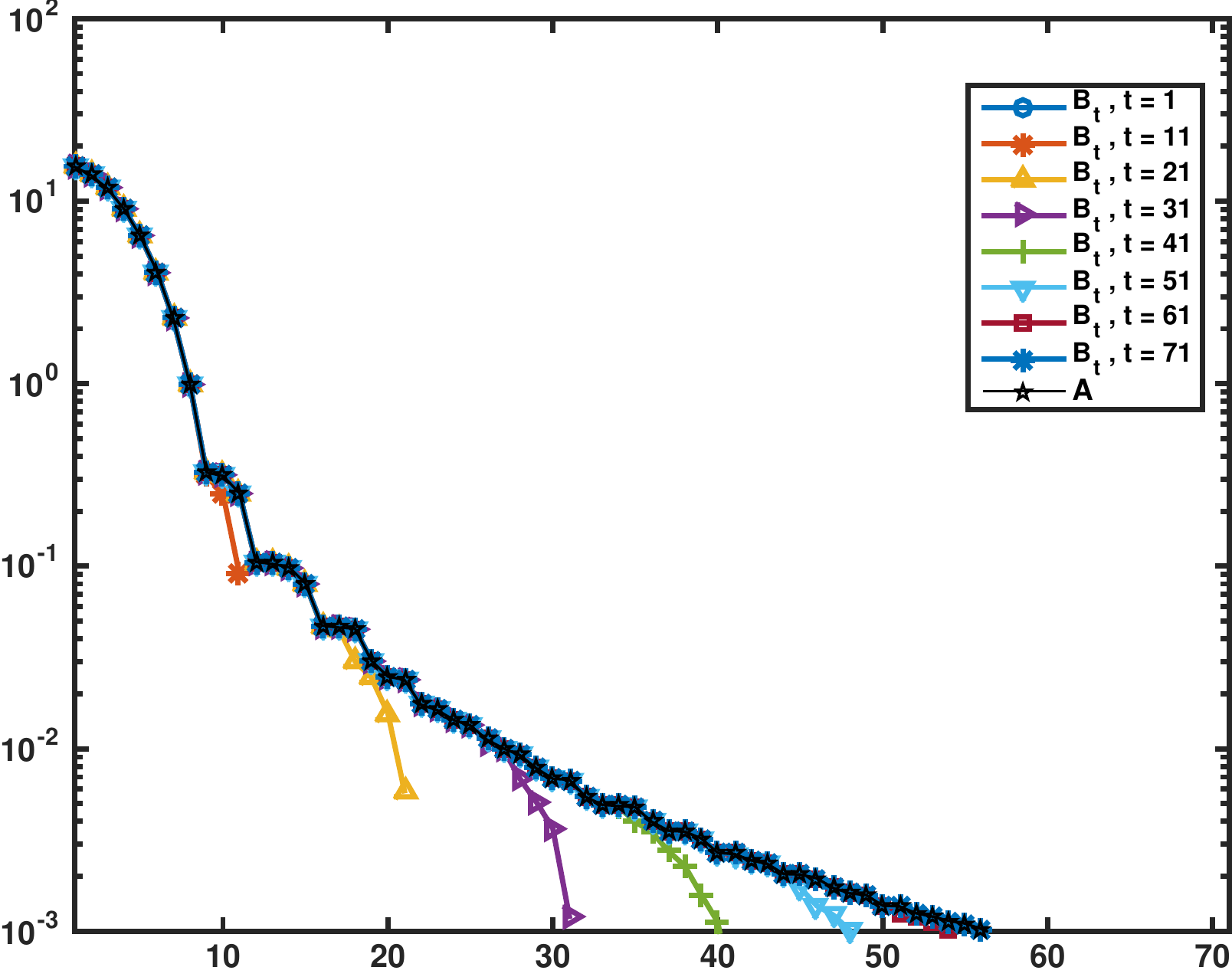}}\hspace{1cm}
\subfigure[\texttt{gravity}   $d=.75$]{\includegraphics[width=.4\textwidth]{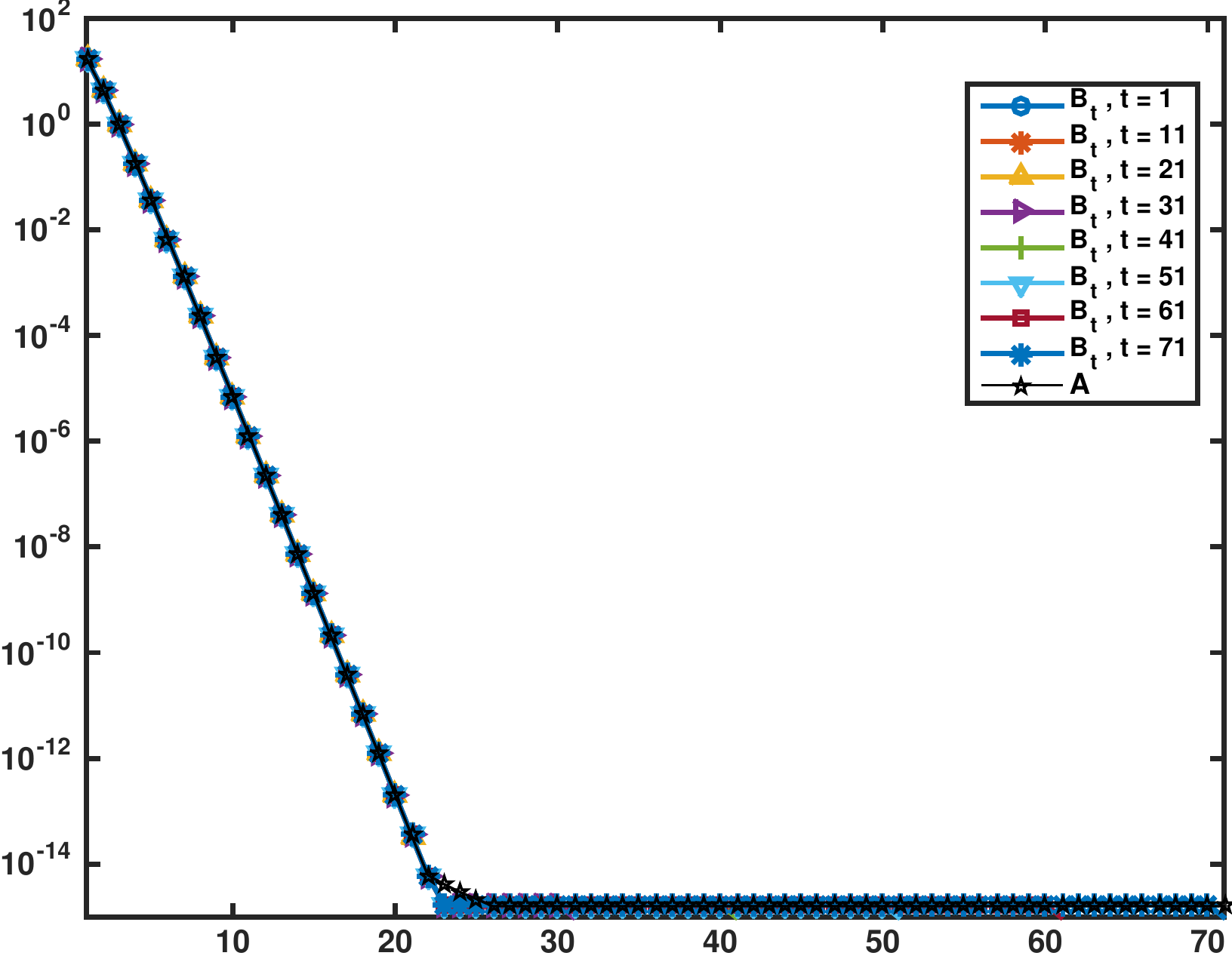}}
\end{center}\caption{Plots of singular values against the index of the singular value for matrix $B_t$ for increasing $t$, $t=1:10:71$ compared to the first $71$ singular values of matrix $A$  for underdetermined cases $m=152$ and $n=304$. \label{figspectra}}
\end{figure}

Figure~\ref{figspectra}  illustrates the spectra of the matrices $A$ and $B_t$ for each test problem for the pairs $(m,n)=(152, 304)$, for $t=1\colon10\colon91$, and with noise level $\eta=.005$ used in the calculation of $B_t$. 
For  \texttt{phillips} there are clear steps in the singular values at indices $9$, $12$, $16$, $19$ and $22$, and the problem is only moderately ill-posed. In contrast,  \texttt{gravity} is severely ill-posed, the singular values decay continuously and exponentially to machine precision. Because \texttt{gravity} is severely ill-posed the LSQR iteration quickly captures the dominant right singular subspace for small $t$. On the other hand, for \texttt{phillips}, the slower decay of the spectrum and the generation of small Ritz values introduces inaccurate small singular values into the spectrum of $B_t$, as is seen by the departure of the spectrum of $B_t$ from the spectrum of $A$. This may present difficulty for estimating the regularization parameter using the presented approaches, unless $t$ is small.

 \subsection{Estimating $t$: the subproblem size}
\begin{figure}[!htb]
\begin{center}
\subfigure[\texttt{phillips}]{\includegraphics[width=.4\textwidth]{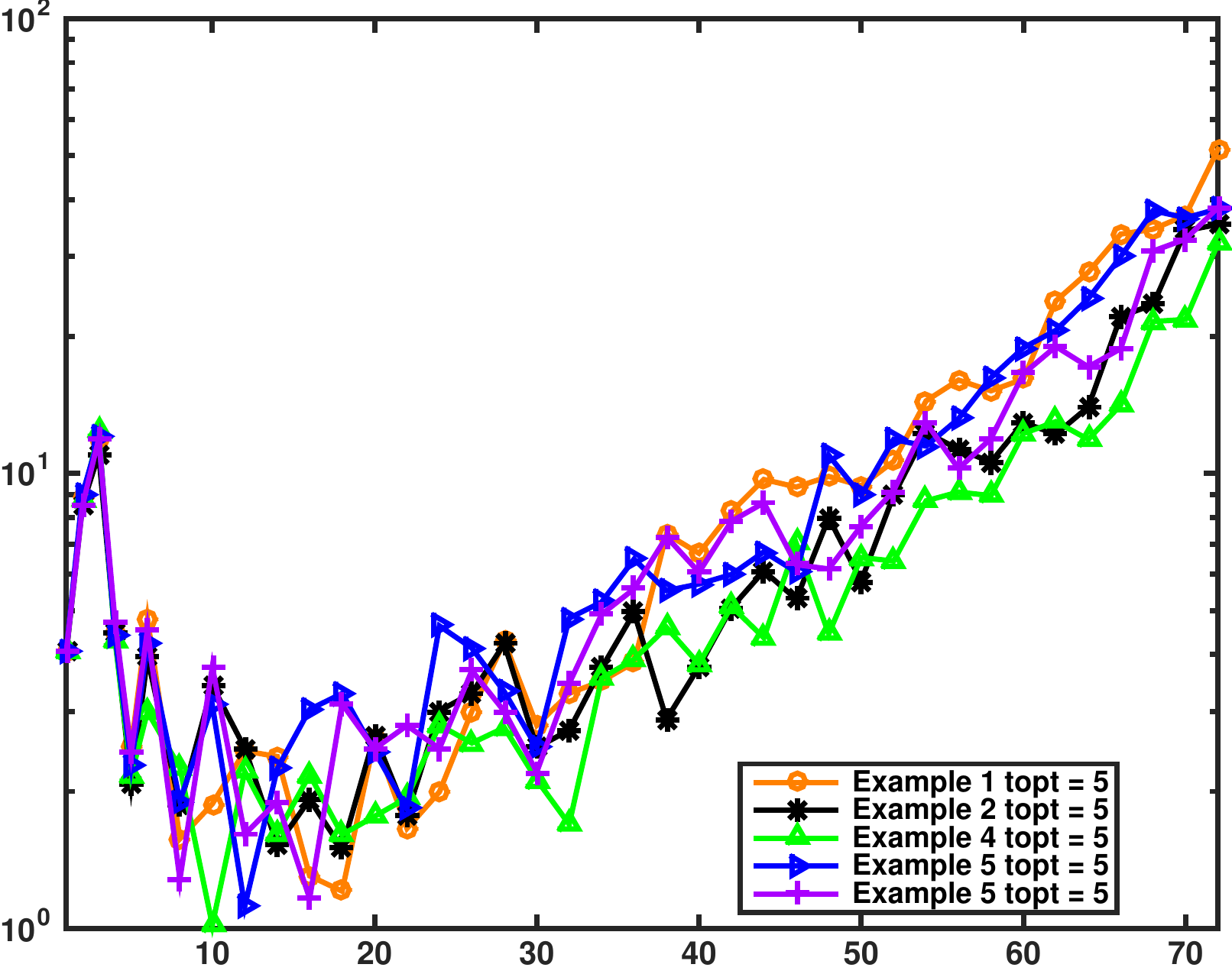}}\hspace{1cm}
\subfigure[\texttt{gravity}  $d=.75$]{\includegraphics[width=.4\textwidth]{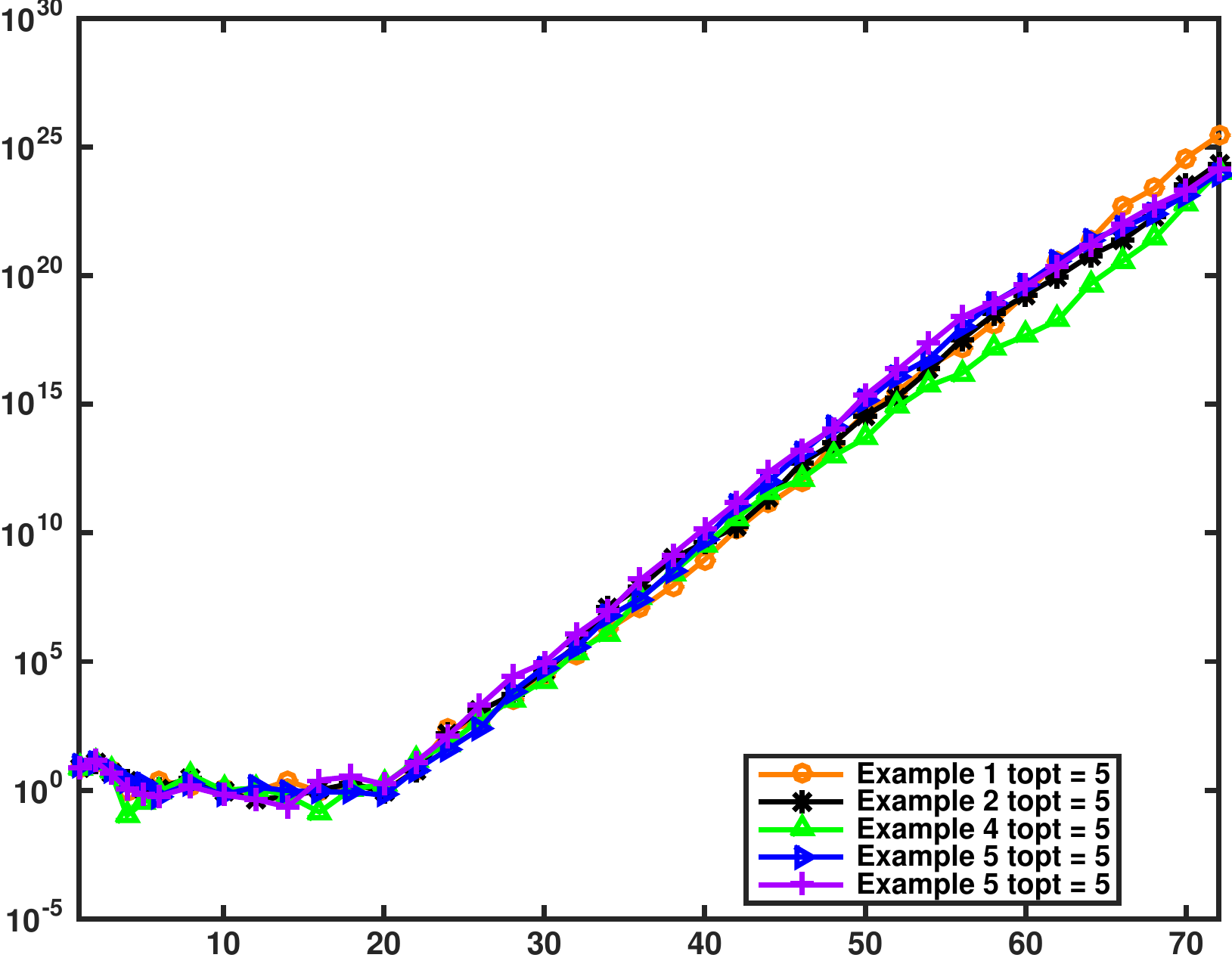}}
\end{center}\caption{Noise revealing function $\rho(t)$ for the two test problems for underdetermined cases with noise level $\eta=.005$  as in Figure~\ref{figspectra}. \label{fig:rho}}
\end{figure}

 To determine the size of the subproblem we examine the approach suggested in \cite[(3.9)]{HPS} for determining the appearance of noise in the subspace. Denoting  the diagonal components of $B_t$  by $\theta_j$, $j=1\colon\!t$ and the subdiagonal entries by $\beta_j$, $j=2\colon\!t+1$,  the cumulative ratio, $\rho(t)=\prod_{j=1}^t (\theta_j/\beta_{j+1})$, shows the impact of $\beta_{j+1}$ approaching the precision of the algorithm as $t$ increases. For large enough $t$ and exact arithmetic $\beta_{t+1}\ll\theta_t$, \cite{HPS}. Noise is identified as entering at optimal iteration $\toptp$,  given by
 \begin{align}\label{tnoise}
\toptp= \min\{ \argmax{{t>\tmin}}(\rho(t))\} +2.
 \end{align}
Here $\tmin$ is problem dependent and chosen to ensure that noise has entered the problem and we  take $2$ steps beyond $\tmin$. There is little difference in the characteristic oscillation of $\rho(t)$  considered in our examples;  noise enters quickly and $\tmin=3$ is already sufficient.  Figure~\ref{fig:rho} illustrates $\rho(t)$ for the  test problems with $5$ samples of the noisy data. $\rho(t)$ correctly identifies the  point when the singular values reach machine precision  for problem \texttt{gravity}. We point out that the noise levels used in these examples, and the subsequent simulations, are substantially larger than the noise levels used in \cite{HPS}, for which the optimal choice of $t$ is thus correspondingly larger. If  we run our examples with less significant noise we do obtain results consistent with \cite{HPS}. It was already noted in \cite{HPS},  that \eqref{tnoise} cannot be used  when the  discrete Picard condition does not hold, e.g. for \texttt{phillips}. An alternative method for identifying the subspace size  is to minimize the GCV function for the TSVD, \cite[3.12]{ChKiOl:15}, 
\begin{equation}\label{gcvtsvd}
\G(t) = \frac{\tmax}{(\tmax-t)^2} \sum_{t+1}^{\tmax} | \bfu_i^T \bfb|^2.
\end{equation}
$\G(t)$ depends on the choice of $\tmax$, i.e. the size of the largest subspace considered, in contrast $\rho(t)$ depends on the selection of $\tmin$ but is independent of $\tmax$. 
 To assess the impact of the \textit{correct} choice of $t$ we also examine solutions obtained with larger subspaces. 
 
 \subsection{Regularization Parameter Estimation}\label{sec:regest}
 In implementing the parameter regularization we test the GCV with $\omega=1$ as compared to WGCV with $\omega=(t+1)/m$, and the MDP with $\delta=m$ as compared to  projected MDP (PMDP) $\delta=(t+1)$.  Estimates using GCV, WGCV, MDP, PMDP and UPRE are obtained by calculating the relevant functions for the same set of regularization parameters, and then minimizing within the region of the optimum as used in Regularization Tools for the GCV. For all tests the regularization parameter yielding a minimum error (denoted in results by MIN) for the projected space is also found by calculating the regularization parameter at $1000$ logarithmically sampled points between $\gamma_1$ and $\max(10^{-14}\gamma_1,\gamma_t)$.\footnote{In practice one would not take such a large selection of regularization parameters, but in tests we found that the minimization step in the UPRE may yield reduced error if the optimum is not found by sampling over a sufficiently fine distribution for the regularization parameter.}

\subsection{Evaluating the results}
\begin{figure}[!htb]
\begin{center}
\subfigure[\texttt{phillips}]{\includegraphics[width=.4\textwidth]  {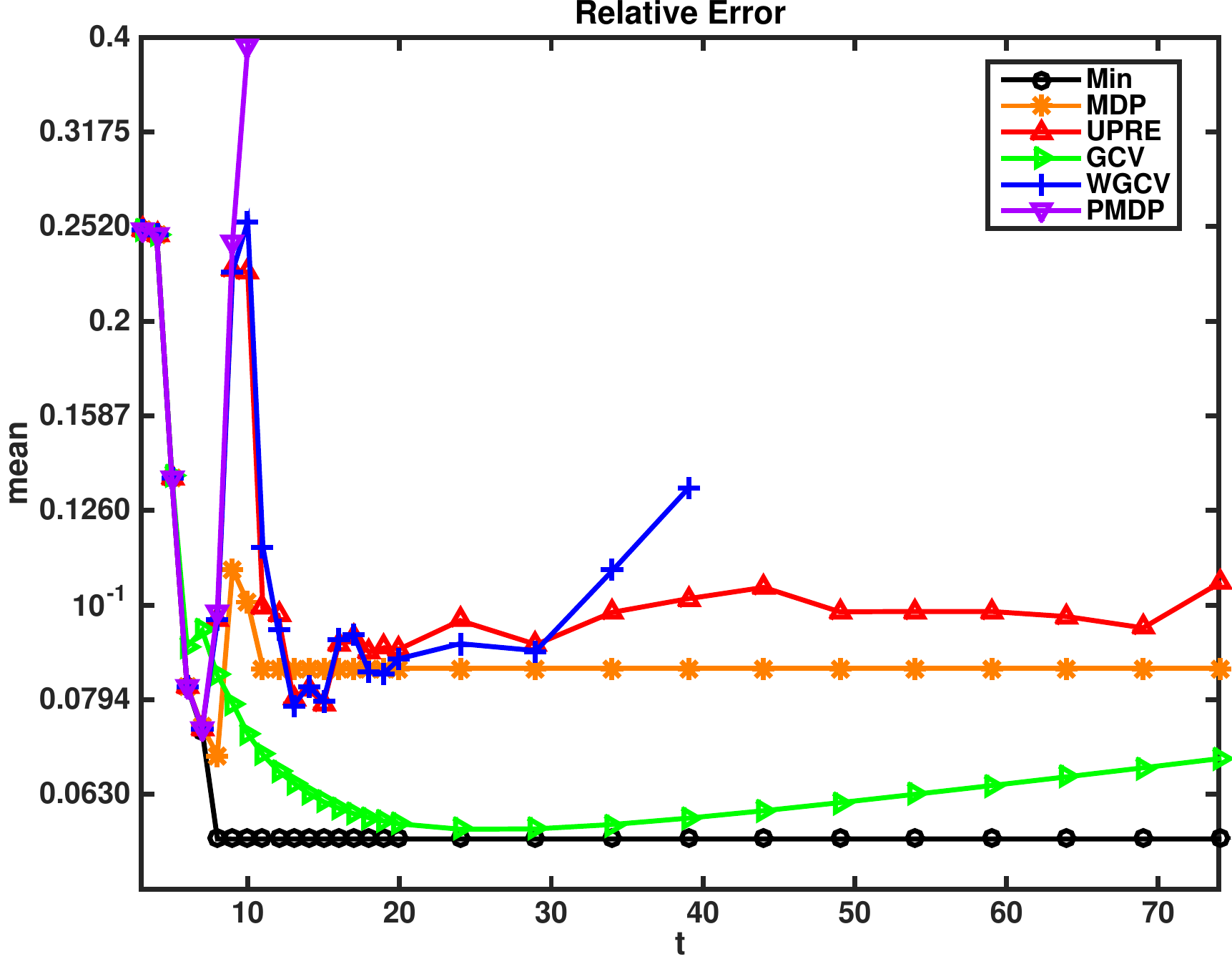}}\hspace{.3cm}
\subfigure[\texttt{gravity}, $d=.75$]{\includegraphics[width=.4\textwidth]{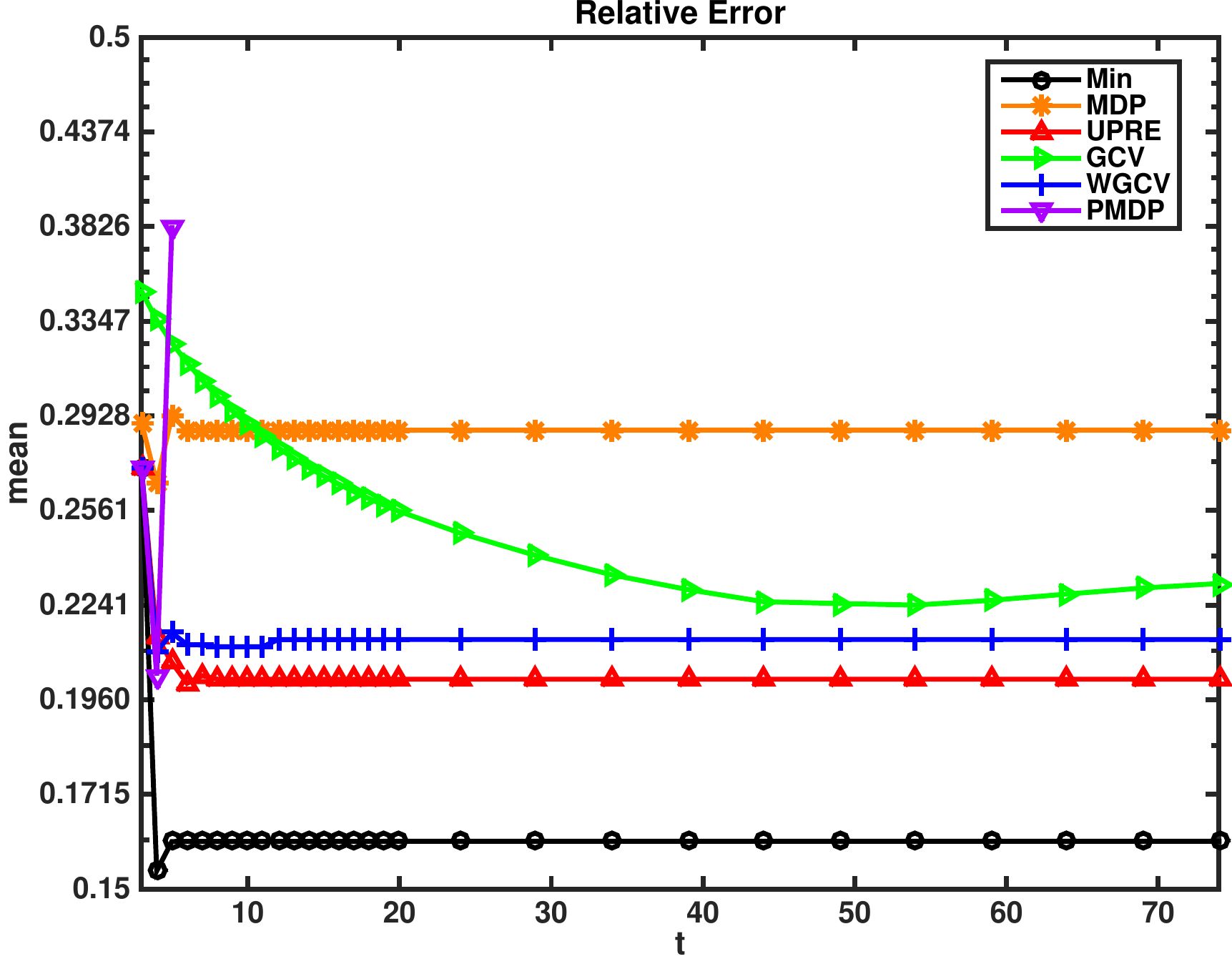}}\hspace{.3cm}
\end{center}
\caption{Average RE over all samples against $t$ for the underdetermined case: $m=152$ and $n=304$ : Noise level $\eta=.005$. \label{figrelerrs}}
\end{figure}

Contaminated data  $\bfb^{c}$ are generated using \eqref{noiselevel} for $c=1\colon50$ yielding solutions $\bfx(t,c)$ for problem size $t$ with $t=[3\colon20,24\colon5\colon74]$. The relative error (RE)  of the solution with respect to the known true solution is given by
\begin{equation}\nonumber
 \text{RE}(t,c) = \|\bfx(t,c)-\bfxtrue\|_2/\|\bfxtrue\|_2. \end{equation}
 The results in Table~\ref{errorsmint} are the average  RE  over all $50$ samples at the reported average $\toptp$ and the minimum RE over all samples and all $t$. These results summarize the graphs of the errors for the same cases in Figure~\ref{figrelerrs}. Immediately we observe that these one dimensional problems generally need rather small subspaces  to obtain optimal solutions with respect to the full problem. They also show that the estimators, other than the PMDP and WGCV  are quite robust to the choice of  $t$,  away from the optimal value.   As suggested by the spectra of $B_t$ and $A$ shown in \S\ref{spectra} the WGCV is robust for \texttt{gravity}, for which the LSQR gives the full regularization,  but not for \texttt{phillips}. Parameter dependent PMDP is not robust in either case. Also the GCV has minimal error for a larger subspace, and thus estimating $\topt$ via $\toptp$ is not effective. On the other hand, using \eqref{gcvtsvd}  $\toptG$ is larger, $6$ and $74$ respectively for the two problems, for which a smaller error is achieved by GCV, $.13$ and $.23$, respectively.   \begin{table}[!htb]
 \begin{center}
 \caption{Average RE over $50$ samples for problem size $m=152$ and $n=304$.  $\tmin=3$ and  with average ${\toptp}$ as given. Results are illustrated in Figure~\ref{figrelerrs}.\label{errorsmint}}
 \begin{tabular}{|*{7}{c|}} \hline
&MIN&MDP&UPRE&GCV&WGCV&PMDP\\ \cline{2-7}
& \multicolumn{6}{|c|}{Average RE: Average $\toptp=5$}\\ \hline
\texttt{phillips}    &     $0.16    $&$0.16    $&$0.16    $&$0.17    $&$0.16    $&$0.16$\\ \hline
\texttt{gravity} $d=.75$    &    $0.17    $&$0.66    $&$0.52    $&$0.35    $&$0.49$     &   $>1$\\ \hline
& \multicolumn{6}{|c|}{Minimum RE (average $t$ for the minimum)}\\ \hline
\texttt{phillips}                  &$0.06$ $(9)$ &$0.07    $ $(8)$ &$0.07    $ $(7)$ &$0.06$ $(24)$ &$0.07$ $(7)$ &$0.07$ $(7)$ \\ \hline
\texttt{gravity} $d=.75$     &$0.15$ $(4)$&$0.27$ $(4)$&$0.21$ $(4)$&$0.22$ $(54)$&$0.22$ $(4)$&$0.22$ $(4)$   \\ \hline
\end{tabular}
\end{center}\end{table}

The results confirm the expectation of the analysis on the performance of the techniques dependent on the degree of ill-posedness of the problem, \texttt{gravity} is more severely ill-conditioned and the WGCV and UPRE estimators give robust solutions independent of $t$, improving on the GCV and MDP solutions. 

\section{Two Dimensional Simulations}\label{sec:simulationtwod}
\subsection{Image Deblurring}
We consider image deblurring problems, \texttt{grain} and \texttt{satellite}, of size $256 \times 256$  from RestoreTools \cite{Restoretools}. Our main aim here is to first demonstrate the use of the regularization techniques PMDP, WGCV and UPRE for increasing $t$,  
and then to examine a stabilizing technique using an IRR, \S\ref{sec:irr}.  Results without IRR are presented for completeness in \S\ref{sec:LSQR} and with IRR in \S\ref{sec:resirr}. 

\begin{figure}[!htb]
\centering
\subfigure[\texttt{True grain}]{\includegraphics[width=.2\textwidth]{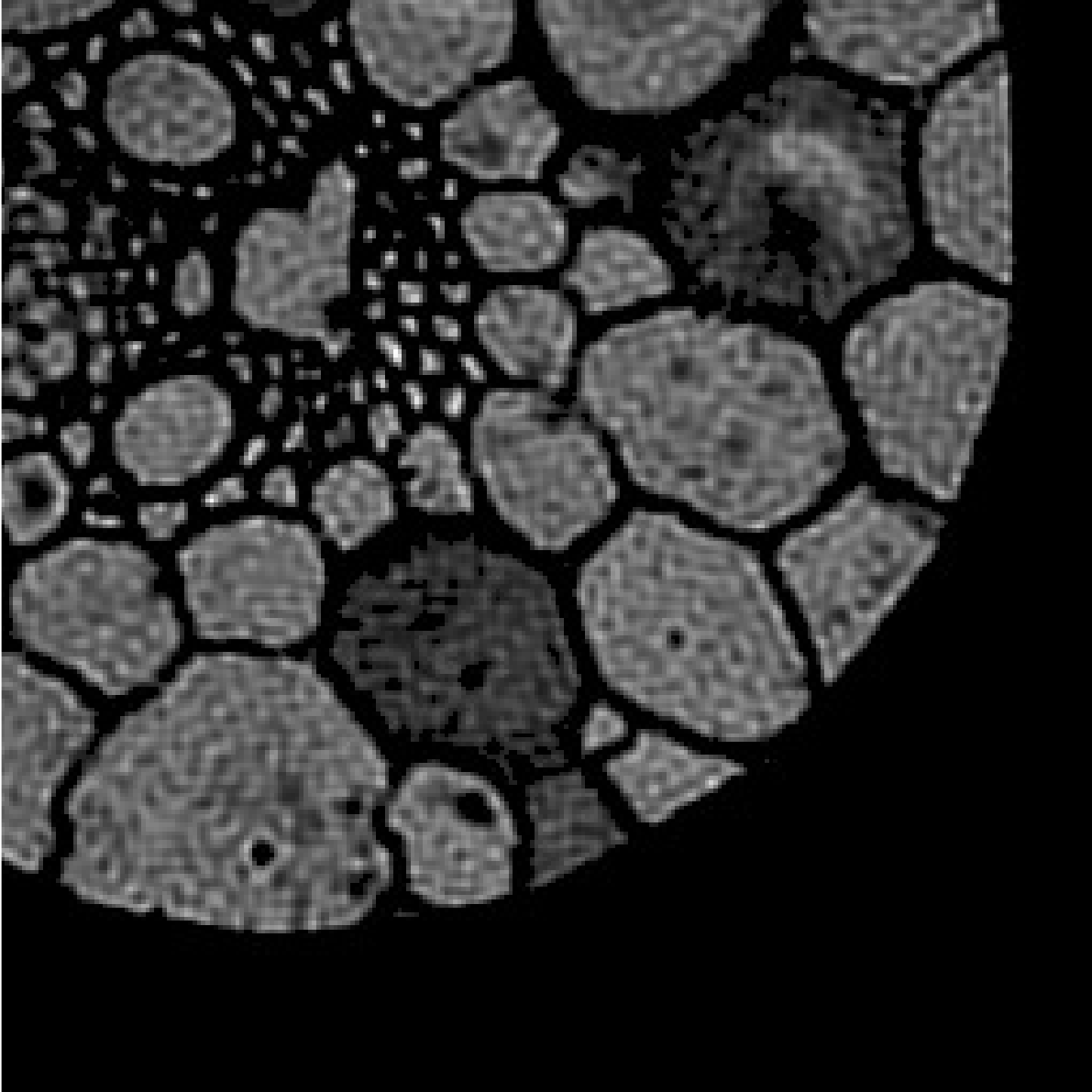}}\hspace{1cm}
\subfigure[\texttt{Contaminated}, $\nu=.1$ \label{badgrain}]{\includegraphics[width=.2\textwidth]{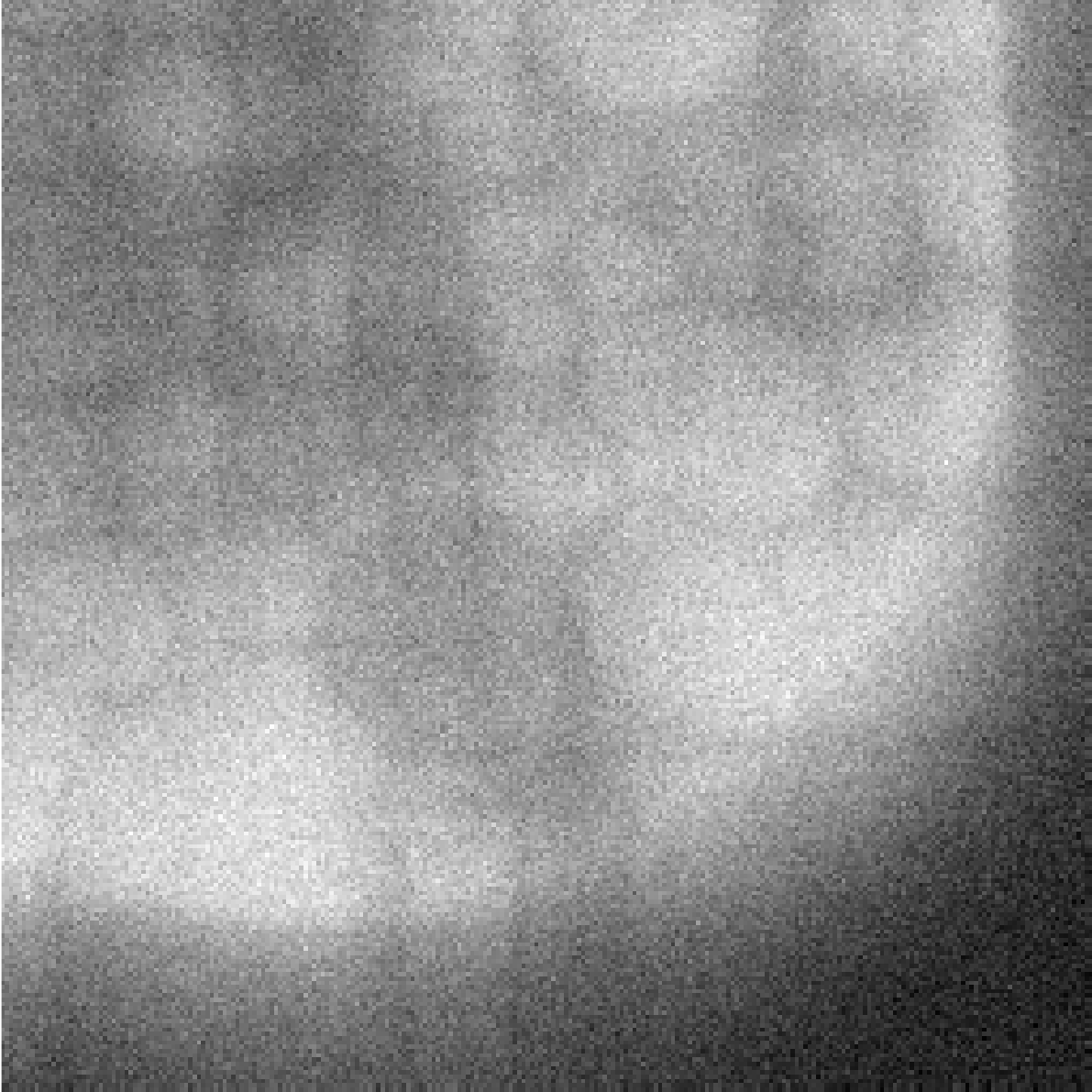}}\hspace{1cm}
\subfigure[\texttt{Point spread function}\label{PSFgrain}]{\includegraphics[width=.2\textwidth]{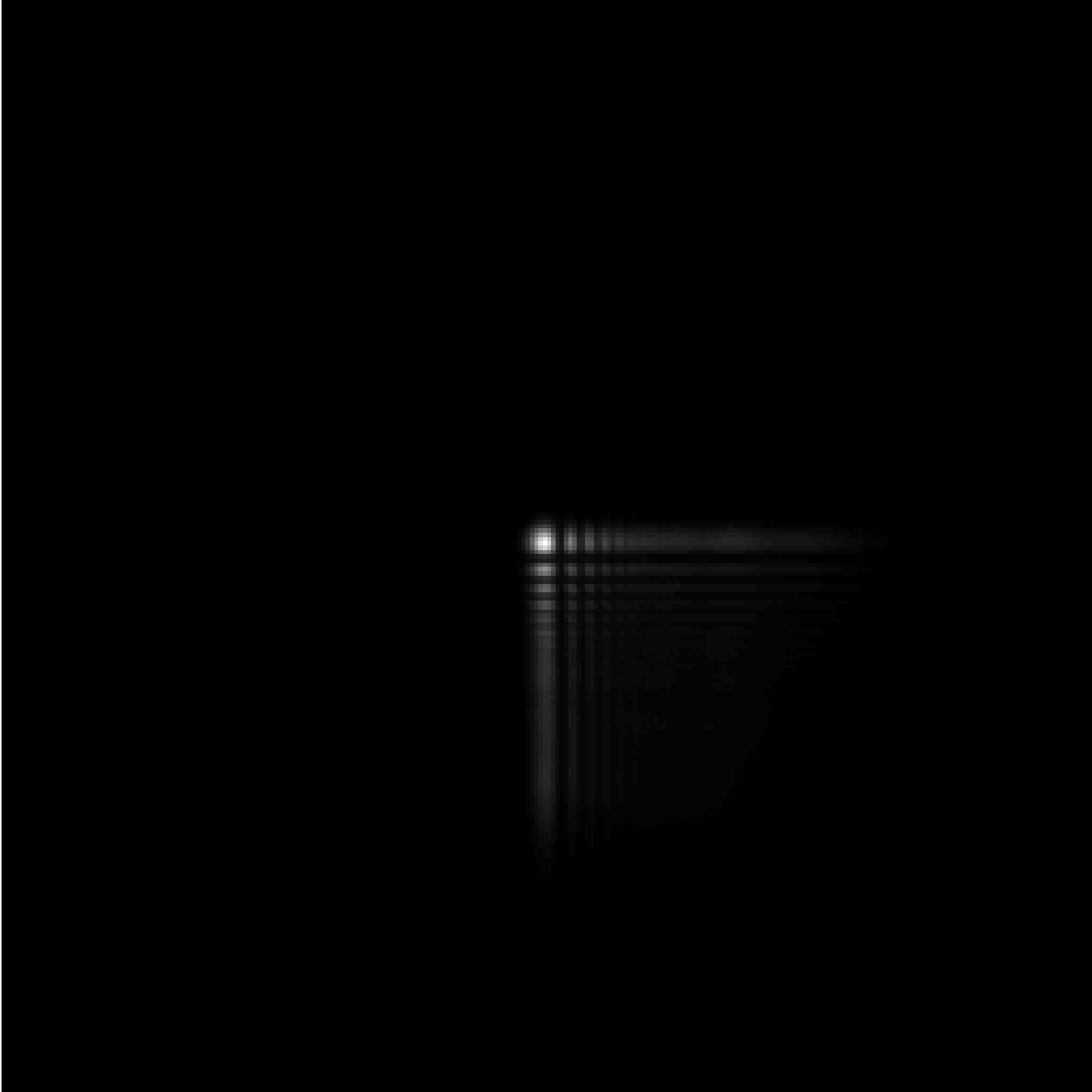}}\\
\subfigure[\texttt{True satellite}]{\includegraphics[width=.2\textwidth]{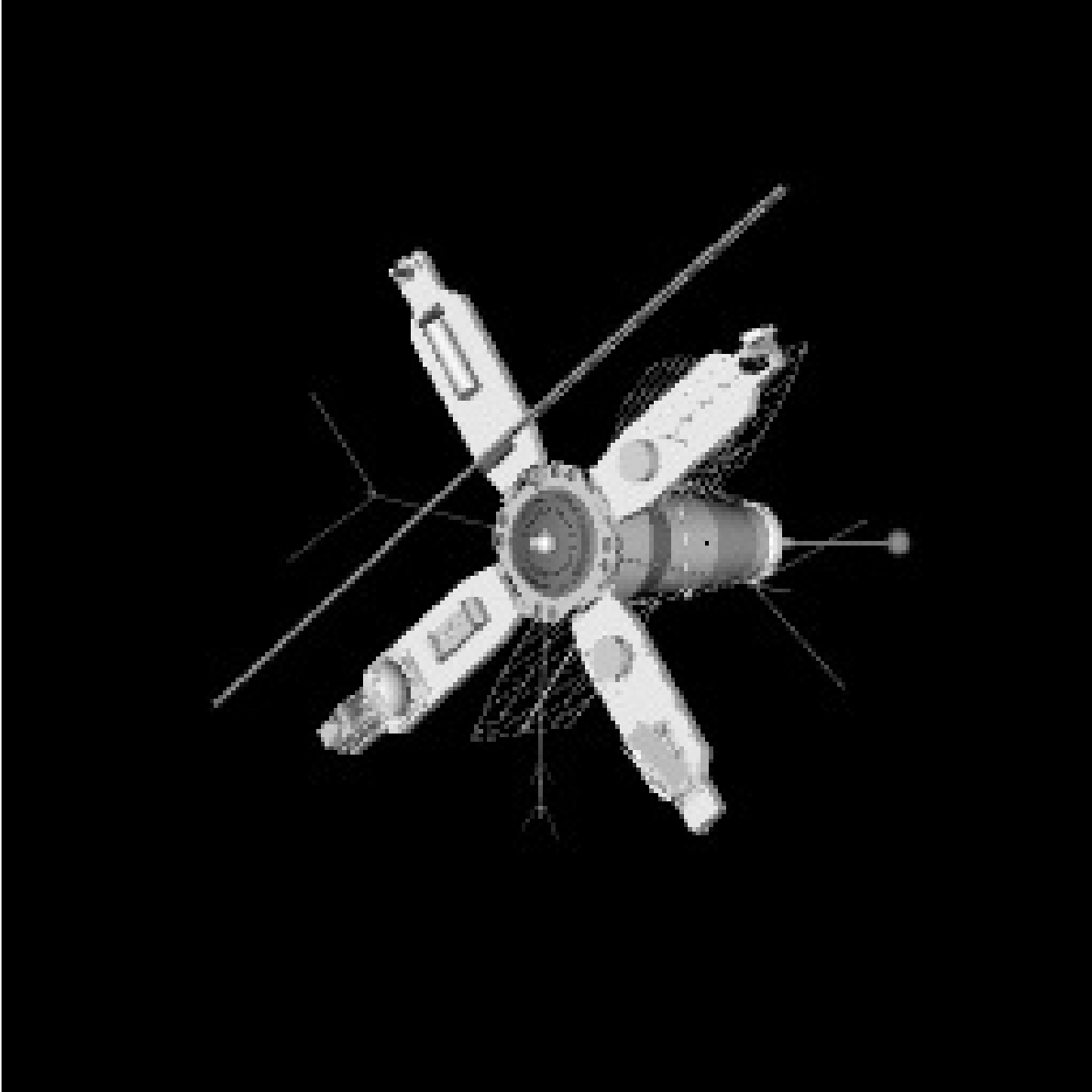}}\hspace{1cm}
\subfigure[\texttt{Contaminated}, $\nu=.1$\label{badsatellite}]{\includegraphics[width=.2\textwidth]{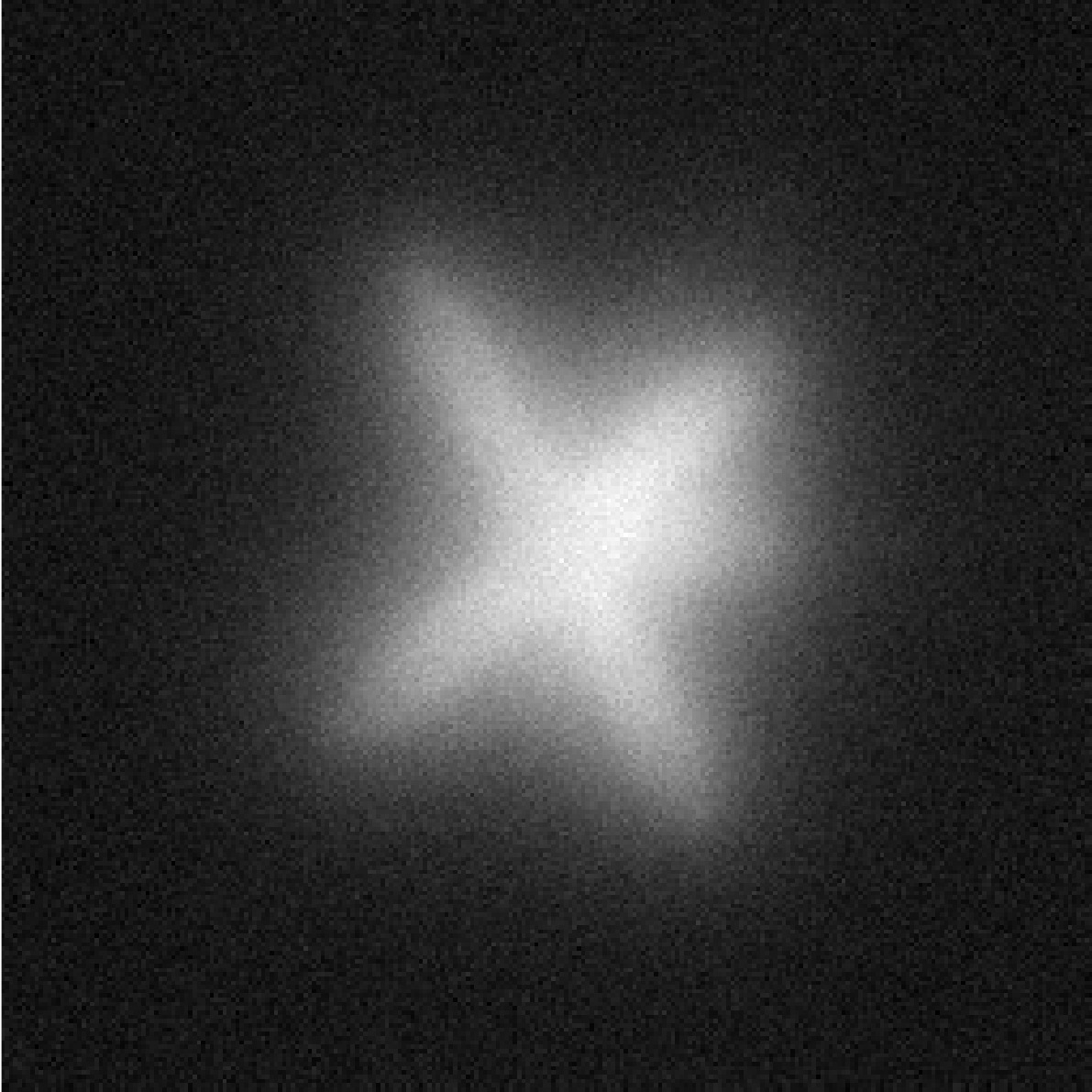}}\hspace{1cm}
\subfigure[\texttt{Point spread function\label{PSFsatellite}}]{\includegraphics[width=.2\textwidth]{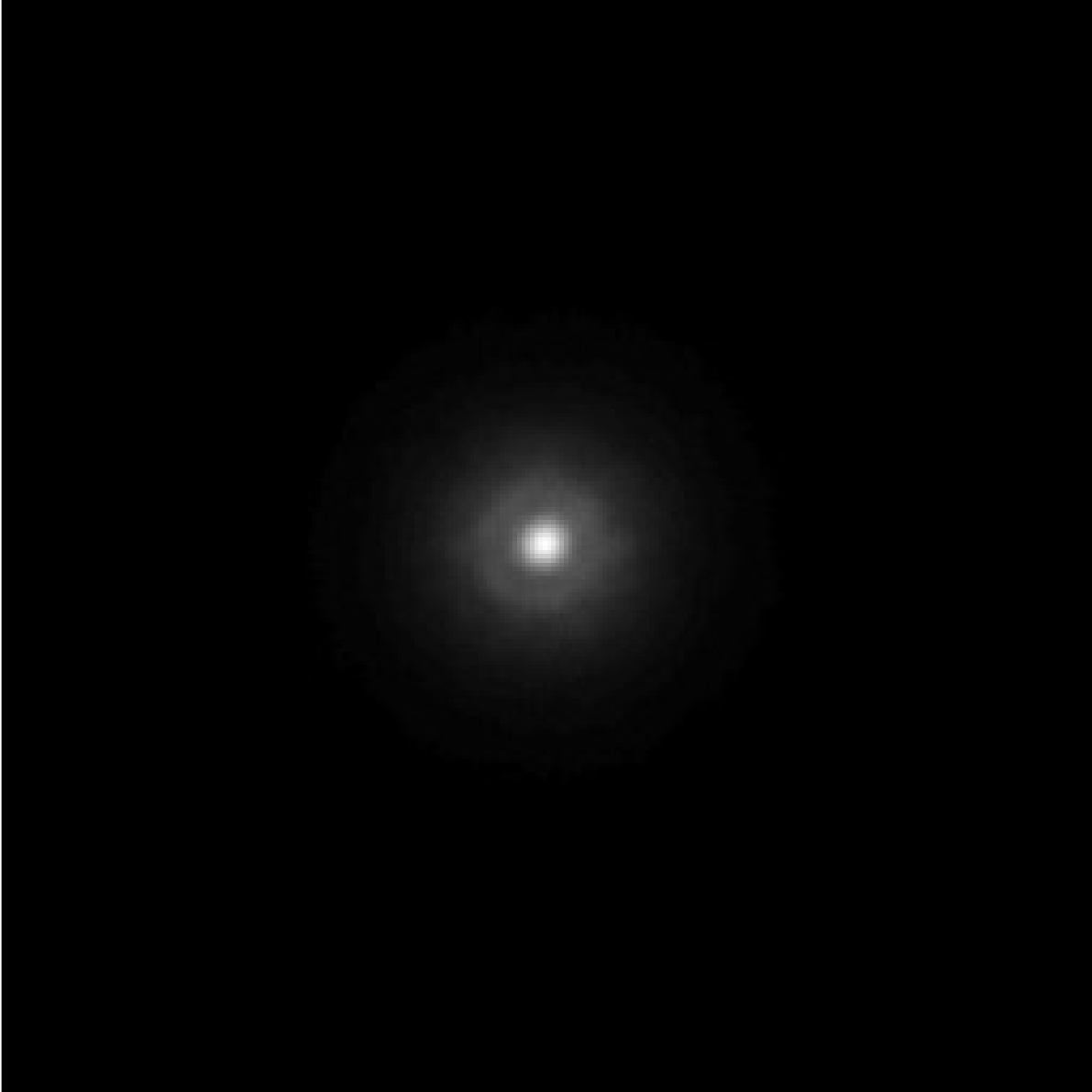}}
\caption{Data for grain and satellite images with blur by the given point spread function  and noise level $10\%$, corresponding to $\nu=.1$. \label{fig2dexample}}
\end{figure}

For contrast with the results presented in \cite{Restoretools} we use noise levels $\eta=.00039$ and $.00019$ in \eqref{noiselevel} which corresponds to noise levels $\nu=10\%$ and  $5\%$, respectively, in \cite{Restoretools} with $\|\bfeta^c\|_2=\nu \|\bfbtrue\|_2$, yielding  $\nu=\eta \sqrt{m}$. These correspond to BSNR $20$DB and $26$DB  as calculated by \eqref{BSNR}.\footnote{For comparison with the results in \cite{ChKiOl:15} we note that there the BSNR is calculated using a definition in \cite{NeChBa}, and their results with $10$DB and $25$DB according to that definition yield   $19$DB and $34$DB with \eqref{BSNR}, resp. i.e. approximately $11\%$ and $2\%$ noise. Hence results here for $10\%$ may be approximately compared with the $10$DB results in \cite{ChKiOl:15}.} For immediate comparison with \cite{Restoretools} we indicate the results using the noise level $\nu$ rather than $\eta$. 
 In Figure~\ref{fig2dexample} we give the true solution, blurred and noisy data and the point spread function.

\subsection{Algorithm Details}
\begin{figure}[!htb]
\begin{center}
\subfigure[\texttt{grain} $\rho(t)$]{\includegraphics[width=.4\textwidth]{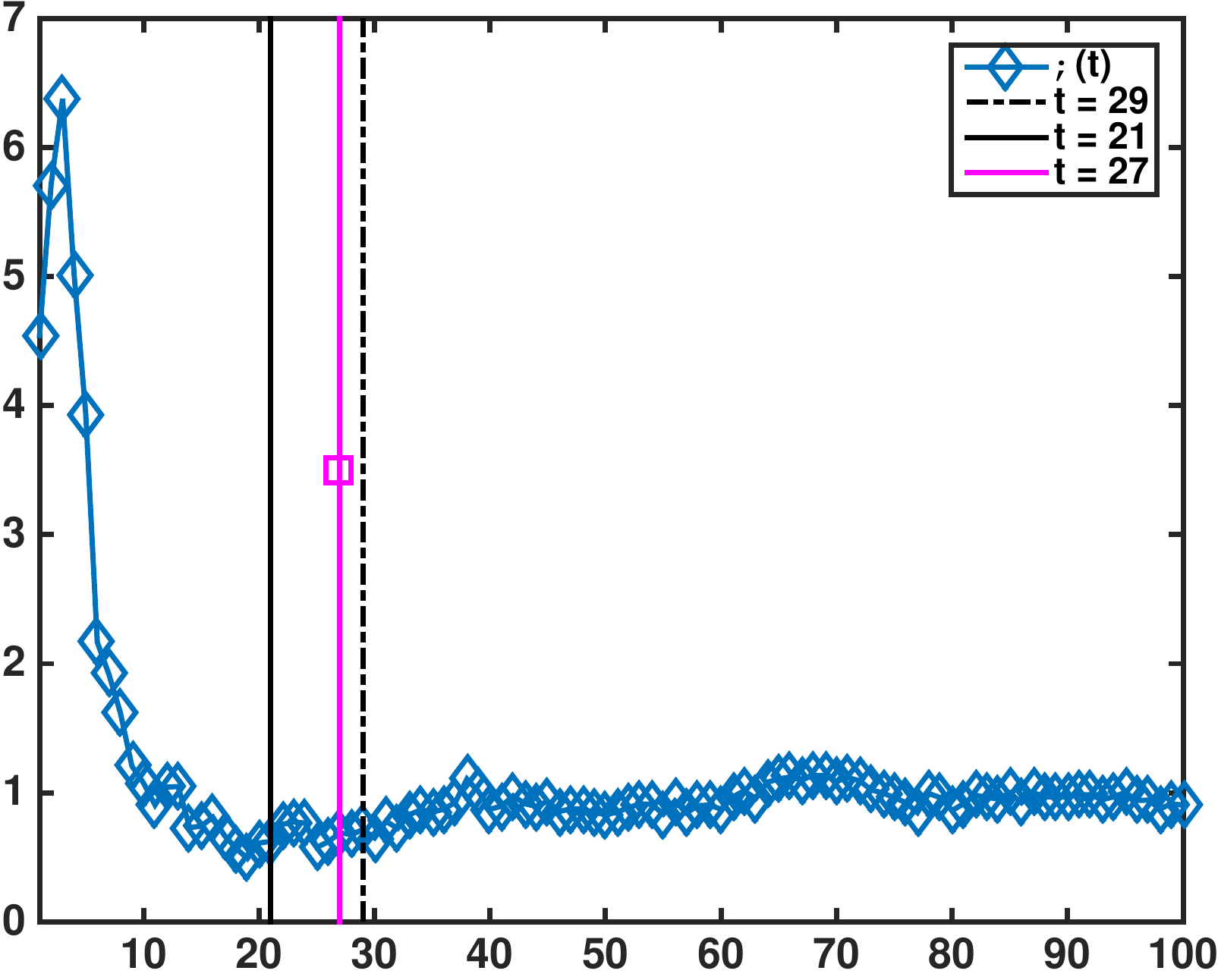}}\hspace{.4cm}
\subfigure[\texttt{satellite} $\rho(t)$]{\includegraphics[width=.4\textwidth]{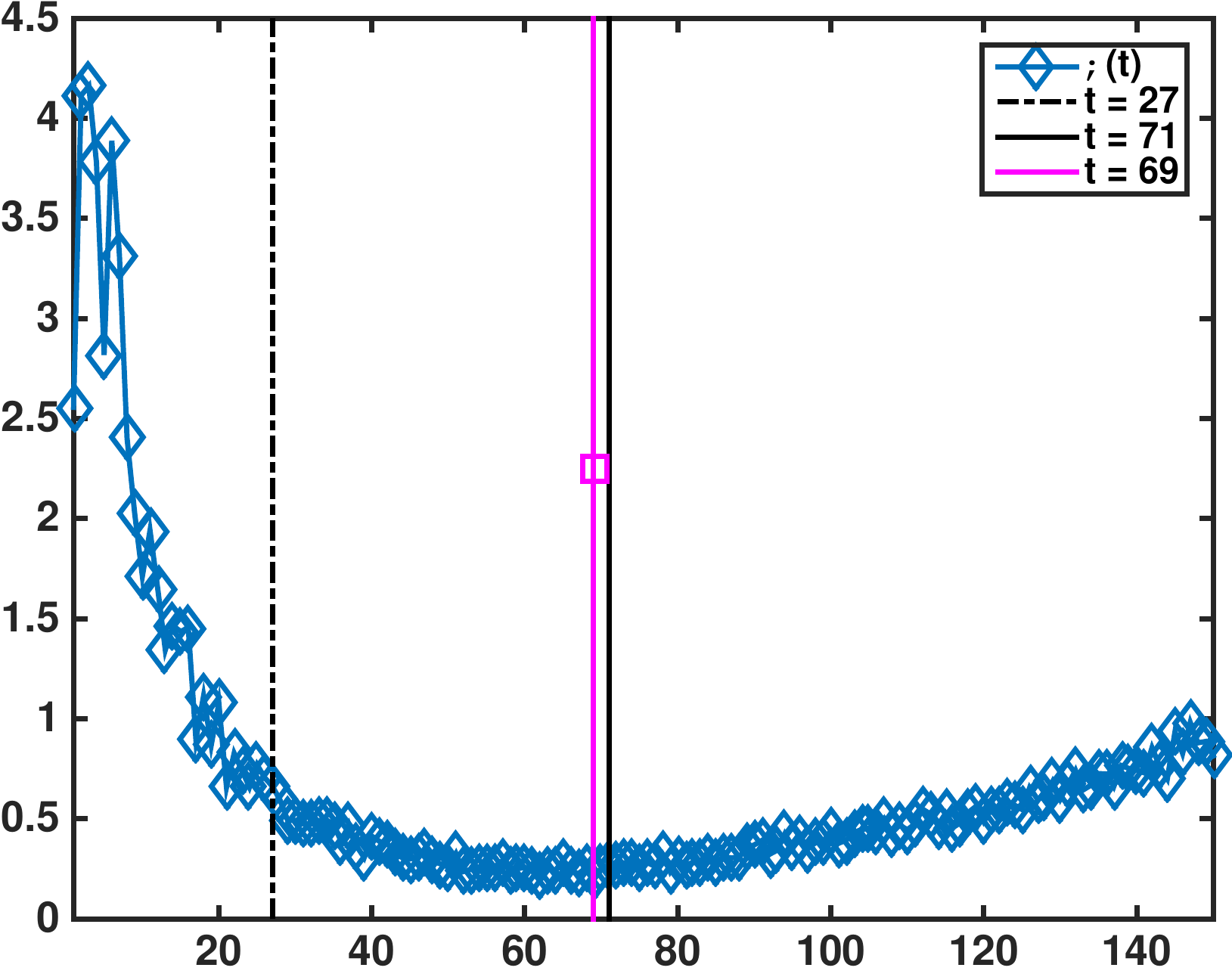}}\hspace{.4cm}
\end{center}
\caption{Noise revealing function $\rho(t)$ for the data illustrated in Figures~\ref{badgrain} and \ref{badsatellite} with $\tmin=25$. Here the dashed-dot vertical line corresponds to the location of $\toptp$, the solid line with symbol to $\toptG$ and the solid line to $\toptmin$. 
 \label{rhotwoD}}
\end{figure}

\begin{figure}[!htb]
\begin{center}
\subfigure[\texttt{grain}]{\includegraphics[width=.4\textwidth]{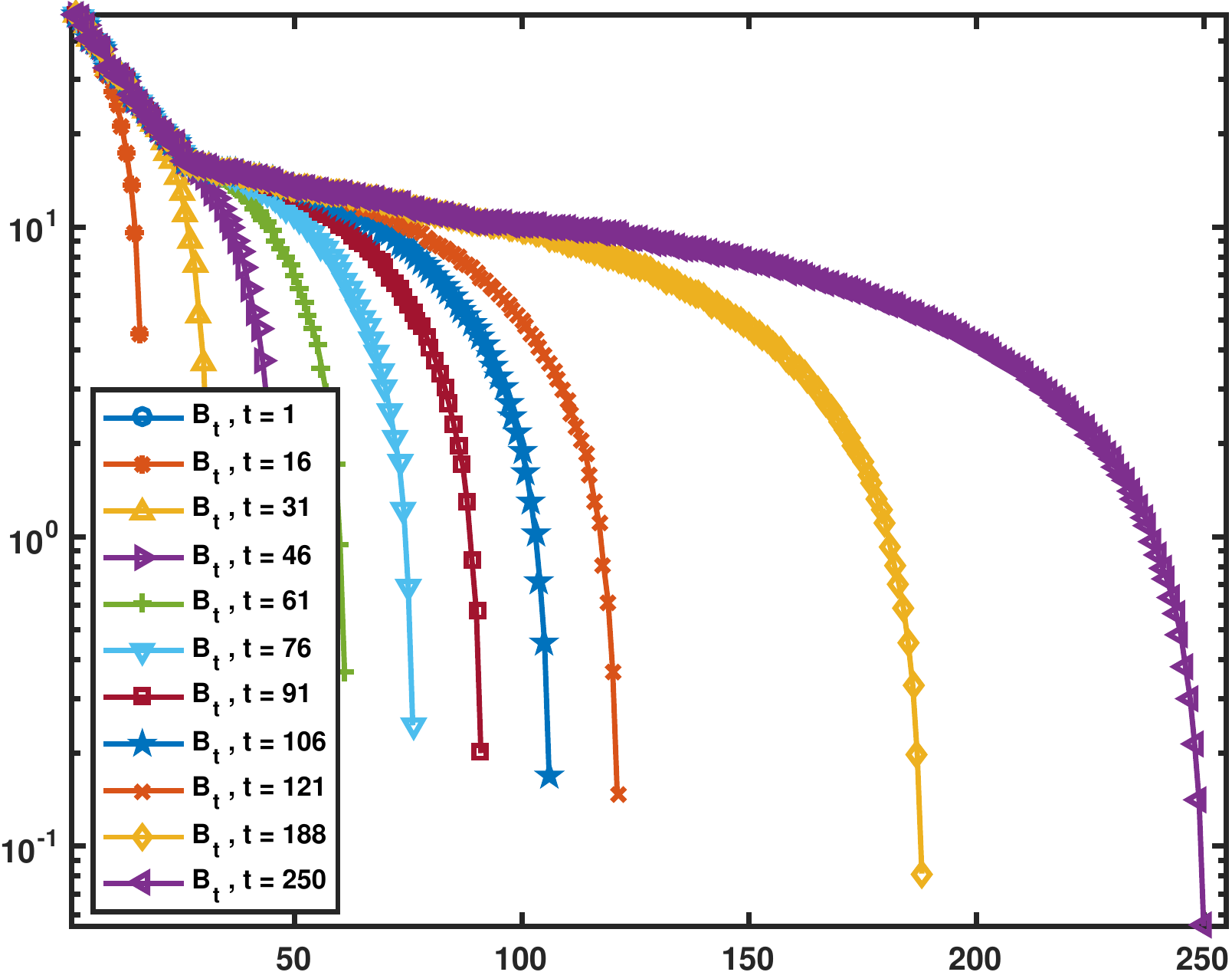}}\hspace{1cm}
\subfigure[\texttt{satellite}]{\includegraphics[width=.4\textwidth]{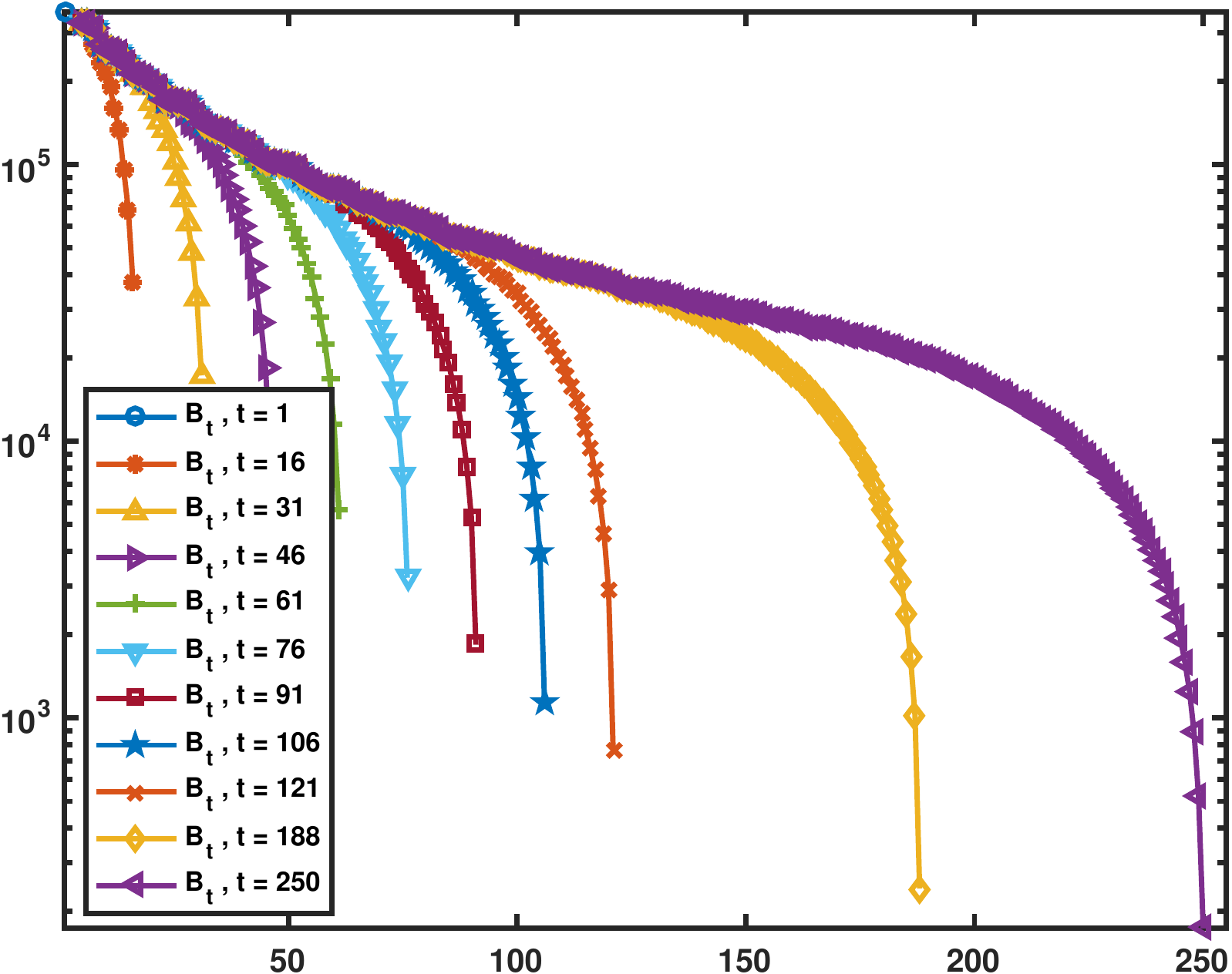}}
\end{center}\caption{Plots of singular values against the index of the singular value for matrix $B_t$ for increasing $t$, $t=1:15:121$, $t=188$ and $t=250$. \label{figspectratwod}}
\end{figure}

In finding the restorations for the data indicated in Figures~\ref{badgrain} and \ref{badsatellite} we note first that the matrices $A$ for the PSFs indicated in Figures~\ref{PSFgrain} and \ref{PSFsatellite} do not satisfy the Picard condition. As illustrated in Figure~\ref{rhotwoD},  $\rho(t)$ does not show the increase within the shown range of $t$ as is clear for $\rho(t)$ with large $t$ obtained for the one dimensional examples in Figure~\ref{fig:rho}. The spectra for increasing $t$, shown in Figure~\ref{figspectratwod} also demonstrate that these problems are only mildly to moderately ill-posed so that the LSQR iterate does not adequately capture the right singular subspace. Still, $\rho(t)$  attains a minimum within the shown  range for $t$ and then exhibits a gradual increase. This suggests that noise is entering the data after the minimum and that one may  use 
\begin{equation*}
\toptmin=\argmin{}({\rho(t)})+2, \end{equation*}
where again we advance $2$ steps under the assumption that noise enters after the minimum. 
In Figure~\ref{rhotwoD} the vertical lines indicate the positions of $\toptp$, $\toptG$ and $\toptmin$. For \texttt{grain}, $\toptG$ becomes quickly independent of the number of terms used, already stabilizing at $\toptG=27$ with just $\tmax=50$ terms, with no change even out to a maximum size of $\tmax=250$ in the calculation. 
For \texttt{satellite} $\toptG$ is less stable  only reaching $32$ when $\tmax=100$ terms are used in the estimation, but increasing to $96$ if $\tmax=250$ terms are imposed. 
Stability in the choice of $\toptp$ with respect to $\tmin$ also follows lack of stability in choice of $\toptG$, suggesting that it is preferable to use $\toptmin$. 
In our experiments we have deduced that it is important to examine the characteristic shape of $\rho(t)$ in determining the optimal choice for the size of the subspace, and will show results using $\toptmin$, $\toptp$ and $\toptG$.

The range for the regularization parameter is also important as is indicated through the windowing approach based on \eqref{gcvtsvd}, \cite{ChEaOl:11}.  From  Figure~\ref{figspectratwod}  it is clear that LSQR iteration only provides a \textit{partial regularization} for either problem and that $B_t$ only captures a portion of the spectrum. Thus we use a single window defined by $t^*$ and apply a FTSVD for the solution which is dominant for the first $t^*$ terms, i.e. with  filter factors $\phi_i(\zeta) \approx 1$ for $i <t^*$ and $\phi_i(\zeta) \rightarrow 0$ for $i>t^*$.    With $\zeta=\tau \gamma_{t^*}$, $\phi_{t^*}(\zeta)=1/(1+\tau^2)<1$.  In our results $\tau=.1$, $t^*=\max({\toptp}, {\toptG})$, and we impose  $\tau \gamma_{t^*} < \zeta <\gamma_1$, for the range of $\zeta$. For the \textit{minimal} (MIN) solutions, the range is adjusted to $10^{-1.5} \le  \zeta \le 1$, consistent with the range for the regularization parameter used in \cite{ChKiOl:15}. This range is scaled by the mean of the standard deviation of the noise in the data, consistent with the inverse covariance scaling of the problems. In finding the MIN solution the range is sampled at $100$ logarithmically sampled points. 

\subsection{Results}\label{sec:LSQR}
\begin{figure}[!htb]
\begin{center}
\subfigure[\texttt{grain} RE]{\includegraphics[width=.4\textwidth]{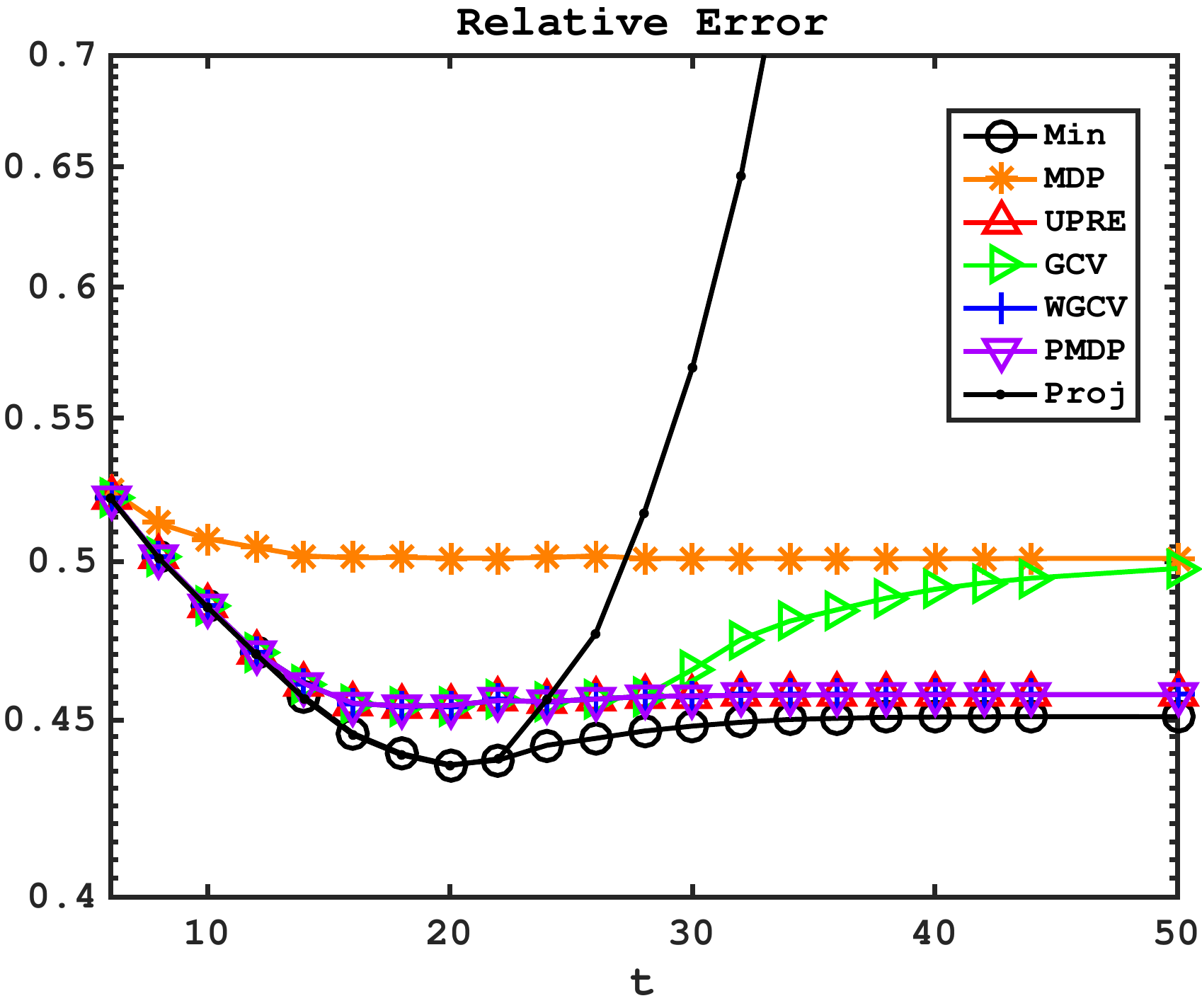}}\hspace{1cm}
\subfigure[\texttt{satellite} RE]{\includegraphics[width=.4\textwidth]{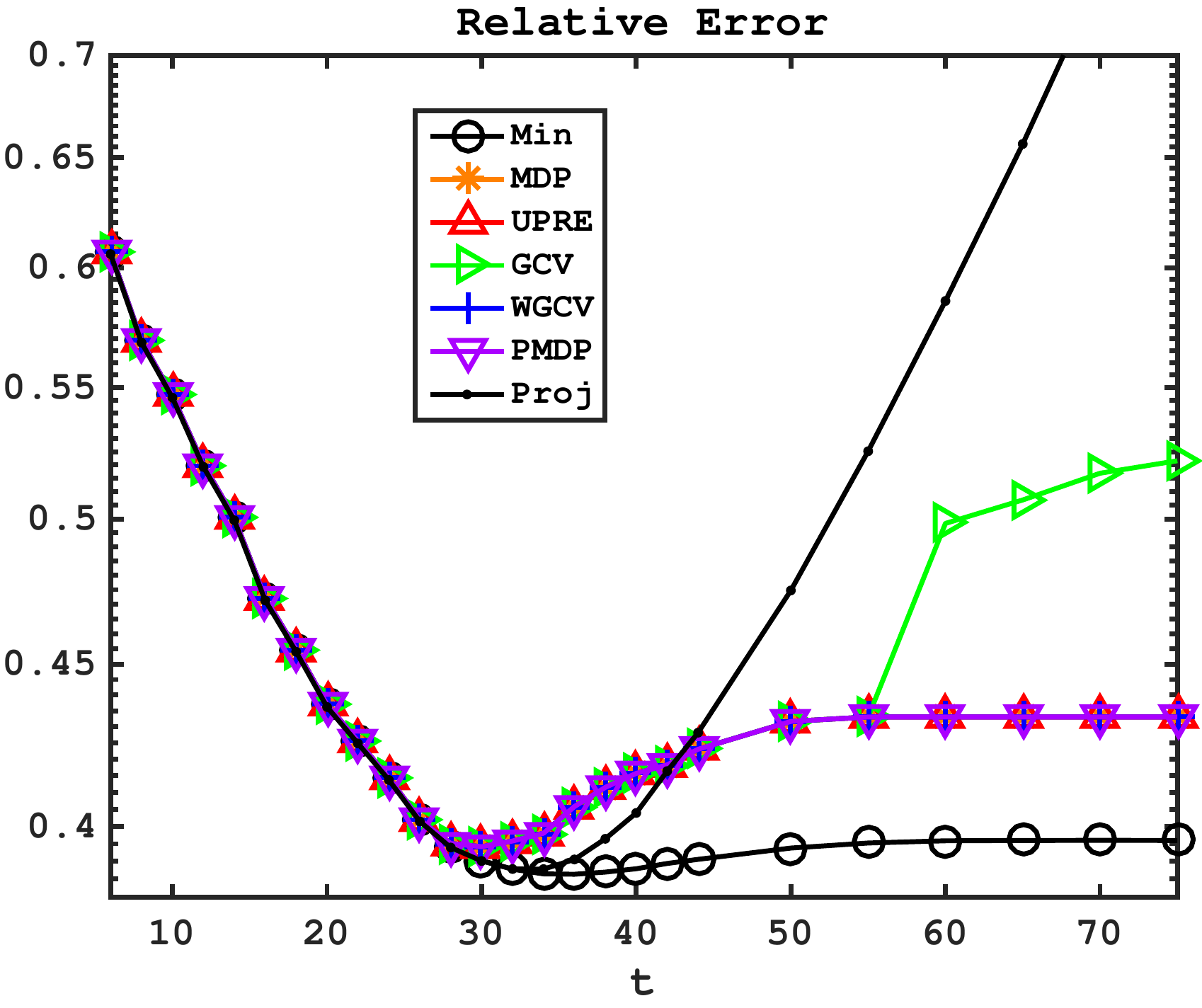}}
\end{center}
\caption{RE  with increasing $t$ with regularization parameter calculated by the different regularization techniques for  examples illustrated in Figures~\ref{badgrain} and \ref{badsatellite}. The solid line in each case is the solution with projection and without regularization. \label{fig2derrorsstep1}}
\end{figure}

\begin{figure}[!htb]
\centering
\subfigure[MIN $\topt=22$]{\includegraphics[width=.17\textwidth]{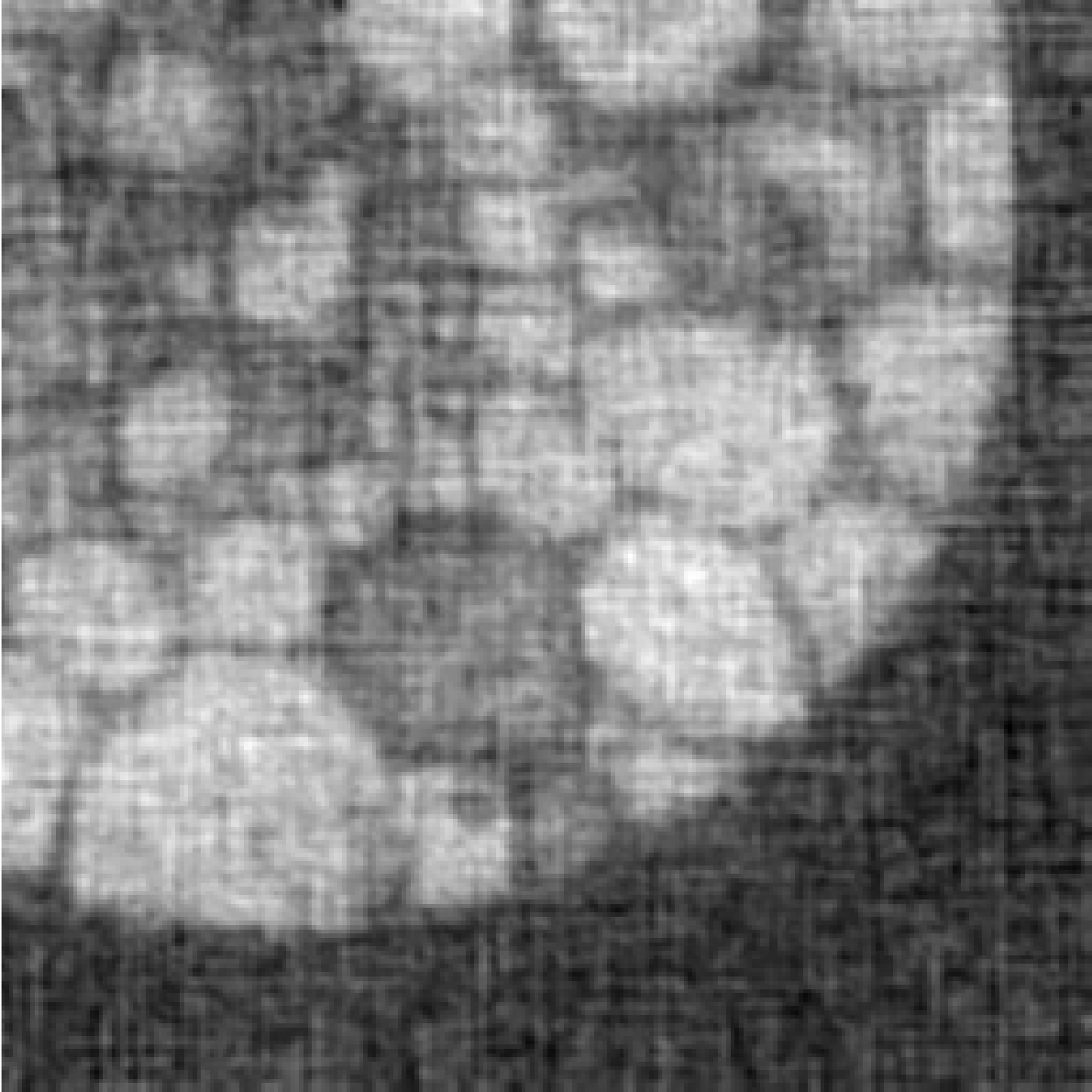}}\hspace{1cm}
\subfigure[$\toptmin=21$]{\includegraphics[width=.17\textwidth]{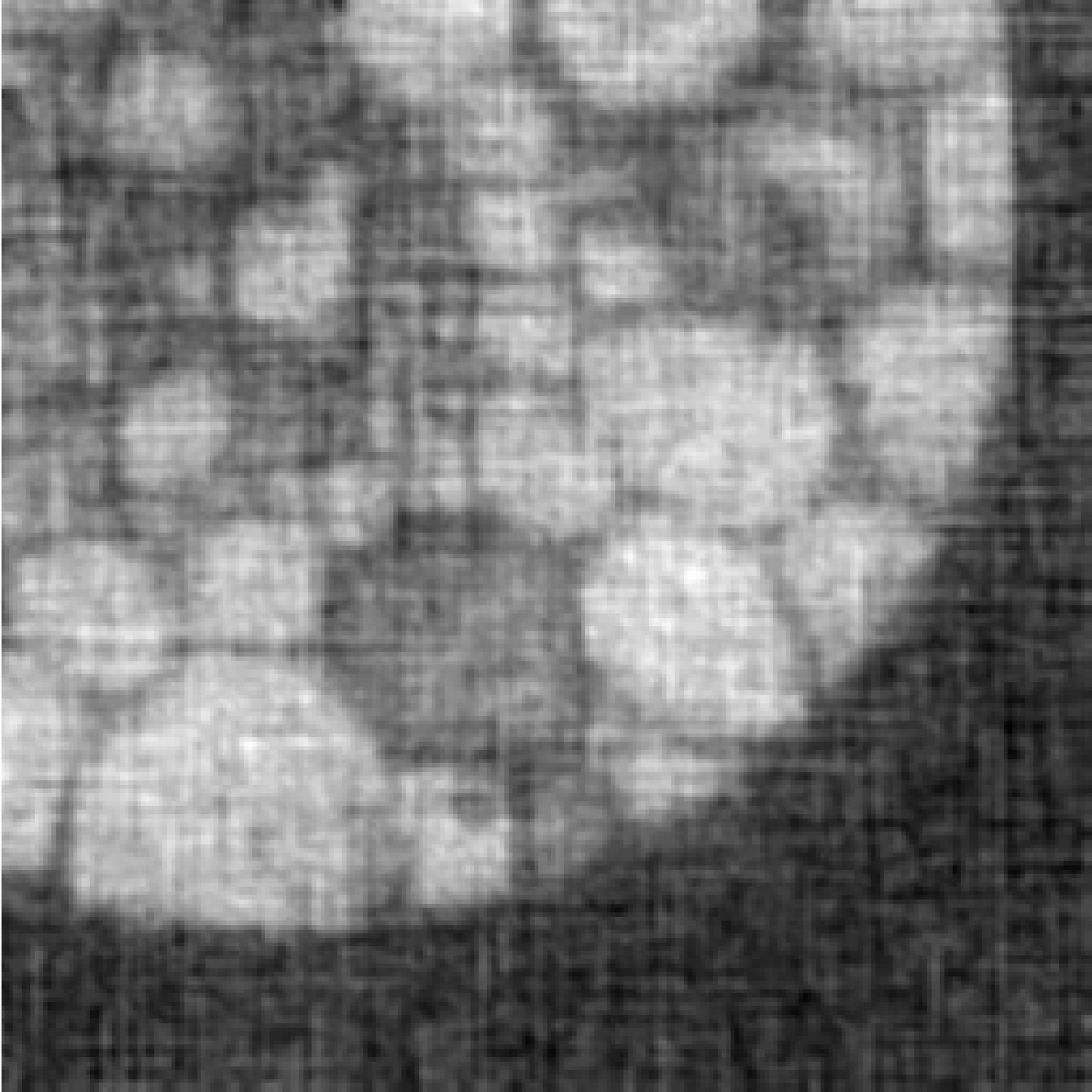}} \hspace{1cm}
\subfigure[$\toptG=27$]{\includegraphics[width=.17\textwidth]{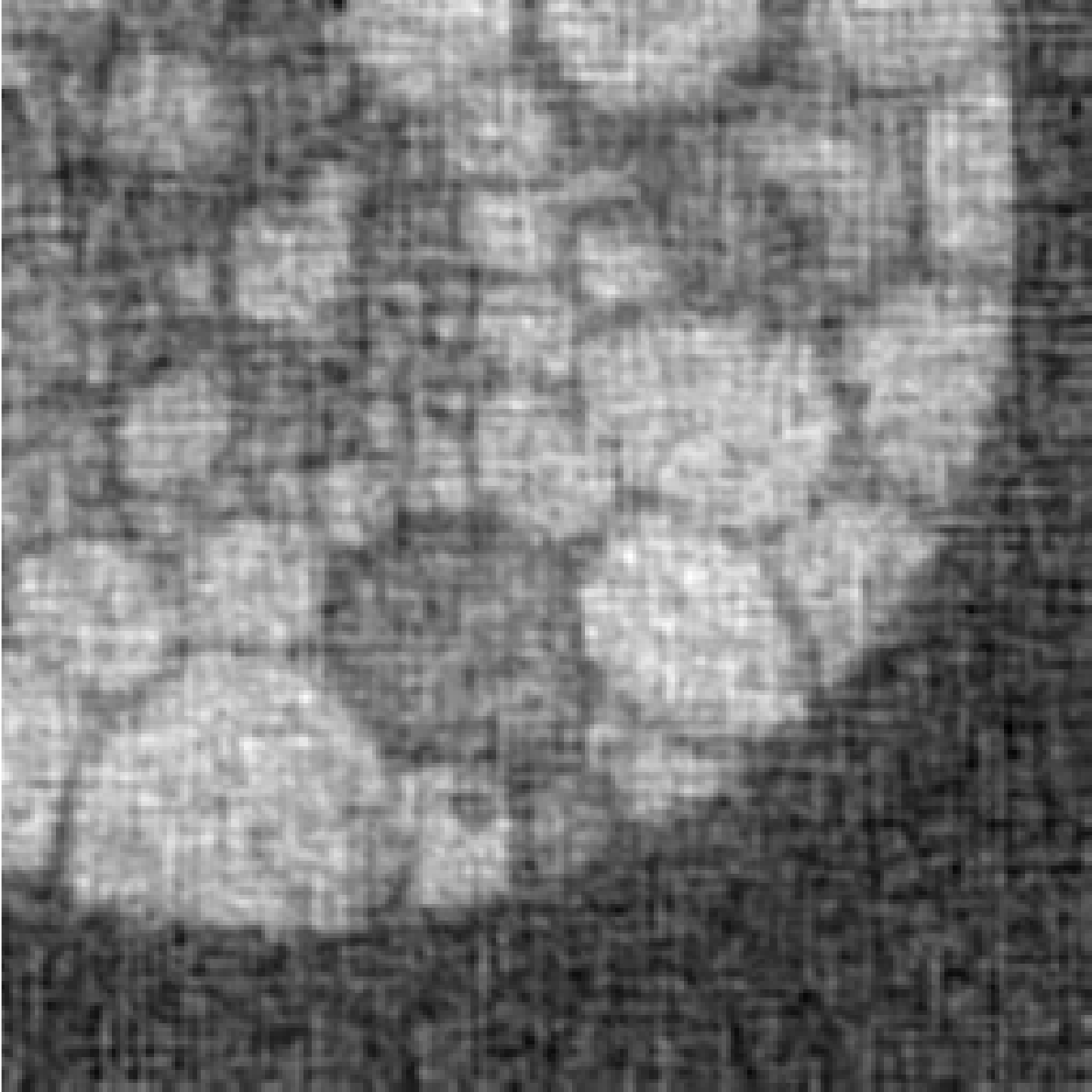}}\hspace{1cm}
\subfigure[$\toptp=29$]{\includegraphics[width=.17\textwidth]{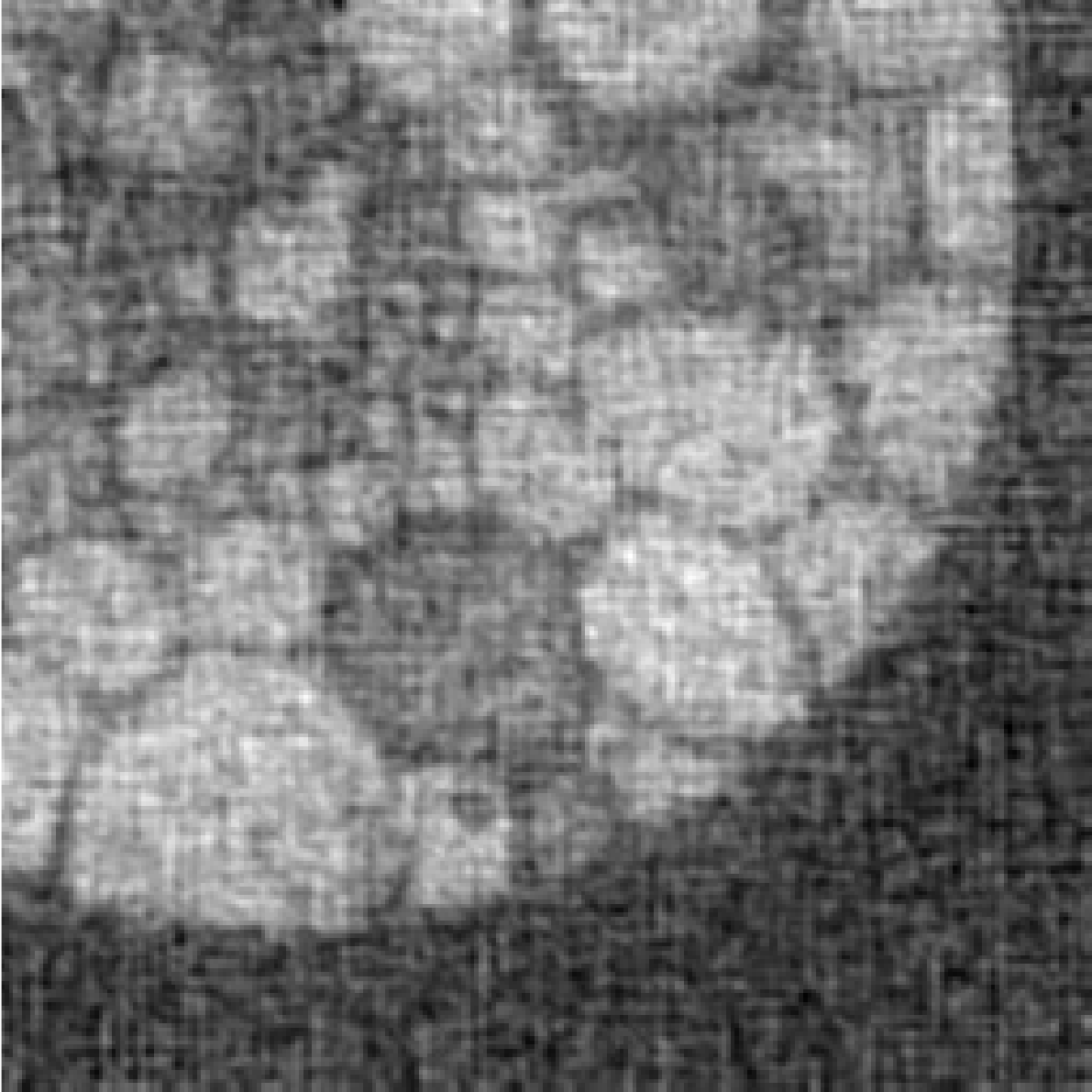}}\\
\caption{Solutions for noise level $10\%$ for  \texttt{grain}  using UPRE to find the regularization parameters and comparing the solutions obtained for $\toptp$, $\toptmin$  and $\toptG$ as compared to the solution with minimum error, MIN. \label{fig2dsolutionsgrain_step1}}
\end{figure}

For quantitative measurement of a given solution as compared to the known solution we again use RE. Other possibilities include using the signal to noise ratio, which is directly related to  RE, and the mean structural similarity index (MSSM) suggested in \cite{WaBoShSi}. We found in our experiments that high MSSM corresponds to low RE and thus providing these results delivers little in terms of further assessment of the algorithms for image deblurring.  The REs using the regularization parameter estimators in contrast to MIN and PROJ are illustrated in Figure~\ref{fig2derrorsstep1} for restoration of the images in Figure~\ref{badgrain} and Figure~\ref{badsatellite}. The  results are consistent with the literature in terms of the semi convergence behavior of the LSQR and also confirm the increase in error seen with GCV without the weighting parameter. Results obtained with  WGCV, PMDP and UPRE are consistent and verify the analysis in Section~\ref{sec:regest}, providing a stable solution for increasing $t$. Solutions  found at the noted $\topt$ for UPRE as compared to the  \textit{optimal} solution with minimum error  are illustrated in Figure~\ref{fig2dsolutionsgrain_step1} for problem \texttt{grain}. Results for the  \texttt{satellite} image are similar. Overall the results  demonstrate that the restorations are inadequate at this level of noise.

\subsection{Iteratively Reweighted Regularization}\label{sec:irr}
Iteratively reweighted regularization provides a cost effective approach for sharpening images e.g. \cite{WoRo:07}, and has been introduced and applied for focusing geophysical inversion, in this context denoted as minimum support regularization,   \cite{PoZh:99,vatan:2014,VAR:2014b,VRA:2014,Zhd:2002}. Regularization operator $D$ is replaced by a solution dependent operator $D^{(k)}$, initialized with $D^{(0)} =I$ and  $\bfx^{(0)}=0$, yields iterative solution $\bfx^{(k+1)}(\alpha)$. For $k>0$, with  
\begin{equation*}
(D^{(k)})_{ii}= ((\bfx^{(k)}_i-\bfx^{(k-1)}_{i})^2 + \beta^2)^{-1/2},
\end{equation*}
where $\beta>0$ is a focusing parameter which assures that $D^{(k)}$ is invertible. Immediately
\begin{equation*}
(D^{(k)})^{-1}_{ii}= ((\bfx^{(k)}_i-\bfx^{(k-1)}_{i})^2 + \beta^2)^{1/2}.
\end{equation*} 
Thus we can use \eqref{tikhonov3} with system matrix $\tilde{\tilde{A}}^{(k)}=\tilde{A}(D^{(k)})^{-1}$,  to obtain the iterative solution $\bfx^{(k+1)}$, $k>0$. Furthermore, it is straightforward to modify the algorithm for calculation of the factorization \eqref{bidiag}  for the left and right preconditioned matrix $\tilde{\tilde{A}}^{(k)}$,  also noting for the specific preconditioners that operations with the diagonal matrix are simple component-wise products. 

\subsubsection{Comments on the parameter $\beta$}
Suppose that $\bfx^{(k)}_i=\bfx^{(k-1)}_{i}$, for $i \in \mathcal{I}$ and $\beta=0$, then $(D^{-1})_{ii} = 0$ for $i \in \mathcal{I}$. Then rather than solving $\tilde{A} D^{-1} \bfz \approx \bfr$, we solve the reduced system $
\hat{A}  \hat{\bfz} \approx \bfr$, where $\hat{A}$ is  $\tilde{A}$ but with column $i$ removed for $i \in \mathcal{I}$ and all other columns scaled by the relevant diagonal entries from $(D^{-1})_{ii}$. Matrix $\hat{A}$ is of size $m \times \hat{n}$ where $\hat{n}=n-|\mathcal{I}|$ and  vector  $\hat{\bfz}$ is vector $\bfz$ with entries $i \in \mathcal{I}$ removed. With regularization  $\hat{\bfz}$  is obtained as the solution of $(\hat{A}^T\hat{A}+\alpha^2 I_{\hat{n}}) \hat{\bfz} = \hat{A}^T \bfr$. Thus $\hat{\bfy}=\hat{D}^{-1}\hat{\bfz}$, where $\hat{D}$ is obtained from $D$ with the same diagonal entries $i \in \mathcal{I}$ removed. The update for $\bfx$ is therefore obtained using \eqref{xupdate} with entries $\bfx^{(k)}_i(\alpha)=\bfx^{(k-1)}_i + \hat{\bfy}_i$, for $i\notin \mathcal{I}$ and $\bfx^{(k)}_i(\alpha)=\bfx^{(k-1)}_i$, for $i \in\mathcal{I}$. Forthwith we use $\beta=0$ and factorize  the reduced system with system matrix $\hat{A}$.

\subsubsection{Algorithmic Details for IRR}\label{sec:resirr}
The approach for the iteration requires some explanation as to how the range of $t$ is obtained at IRR iterations $k>0$. Noise revealing function $\rho^{(k)}(t^{(k-1)}, t)$  depends on the subspace size $t^{(k-1)}$ from the prior step $k$, and current subspace size $t$. Further, $\toptp$, $\toptmin$ and $\toptG$ are all dependent on $t^{(k-1)}$ as well as $\tmin^{(k)}(t^{(k-1)})$ and $\tmax^{(k)}(t^{(k-1)})$, i.e. given a specific subspace size at $t^{(k-1)}$ the minimum and maximum sizes to use at step $k$ need to be specified. Because we anticipate that further noise enters with increasing $k$, we expect $\tmin^{(k)}<\tmin^{(k-1)}$ and that $\tmax^{(k)}(t)<\topt^{(k-1)}$. With these constraints the cost of an IRR step will be less than the first step $k=0$. Examination of $\rho^{(k)}(t^{(k-1)},t)$ is useful in identifying the constraints on $t^{(k)}$.  For each iteration the range for $\zeta$ is constrained using the current singular values, by $\tau \gamma_t^* \le \zeta \le \gamma_1$, where at step $0$, $t^*=\max(\toptp, {\toptG})$ and $t^*=\toptp$ for the IRR updates.  We will examine the choices for the case with $5\%$ noise.

\subsubsection{Results with IRR for $5\%$ noise}
\begin{figure}[!htb]
\begin{center}
\subfigure[$k>0$ $\rho(t)$\label{grainrhostepk}]{\includegraphics[width=.4\textwidth]{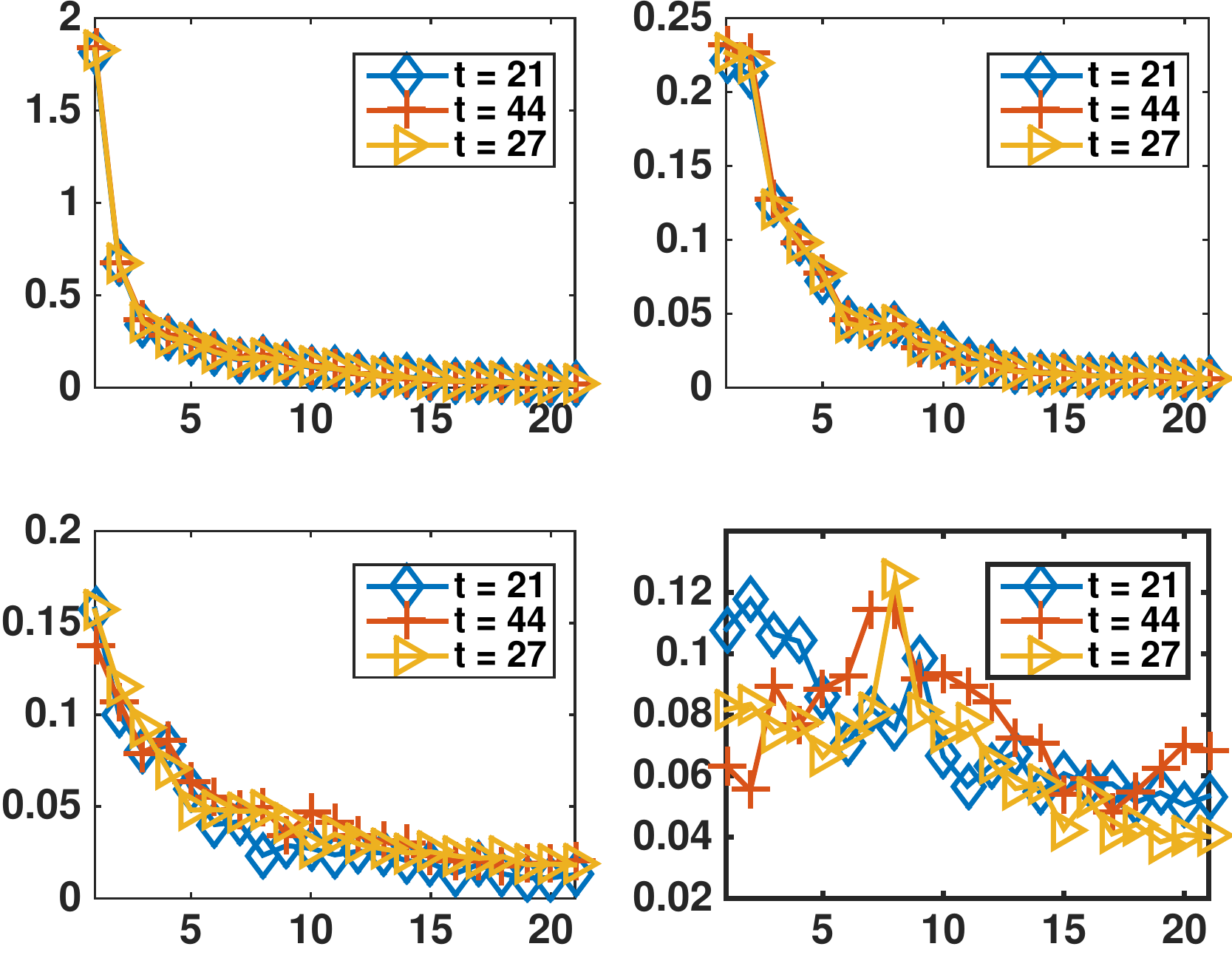}}\hspace{.6cm}
\subfigure[RE\label{grainRE}]{\includegraphics[width=.4\textwidth]{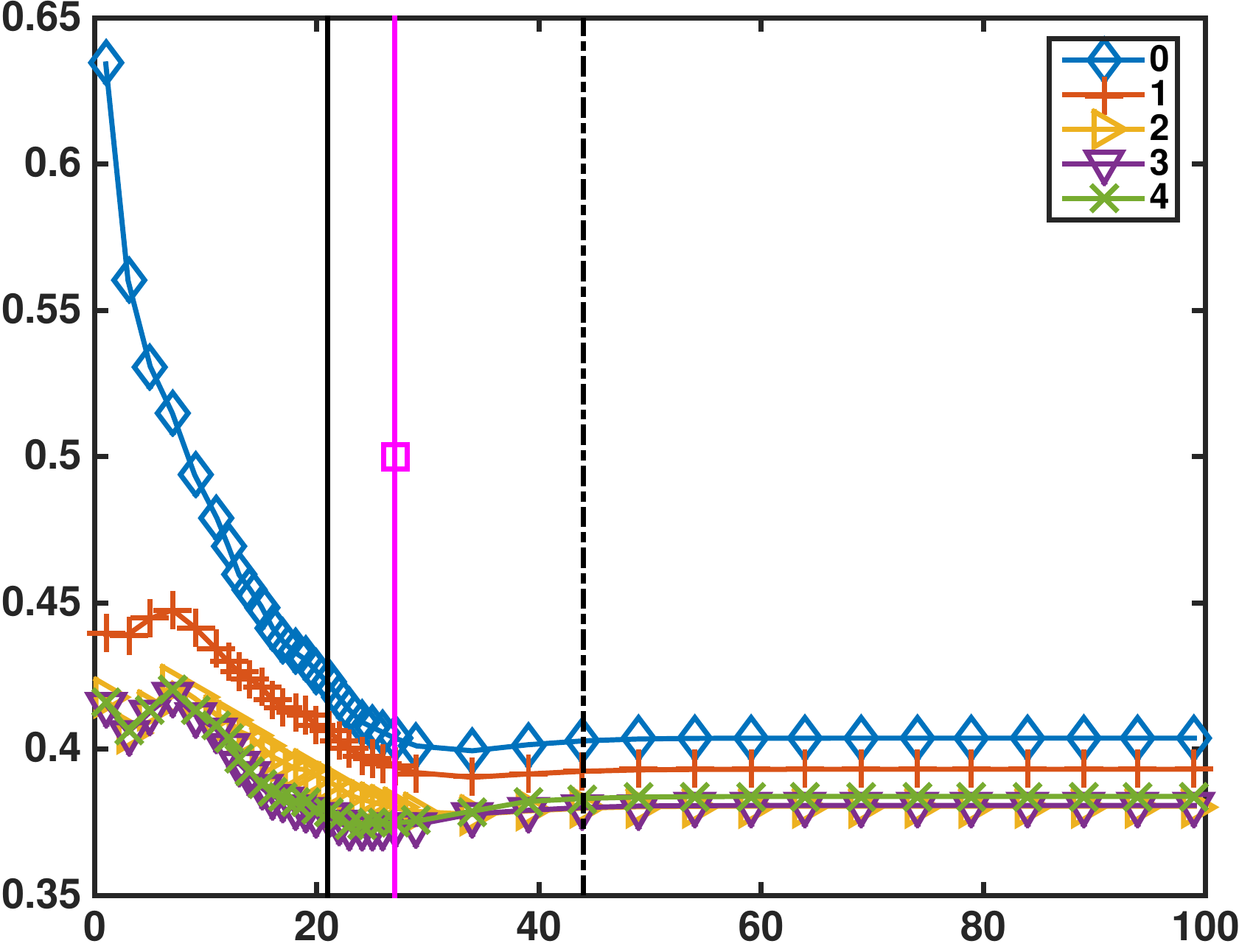}}\hspace{1cm}
\end{center}
\caption{Demonstration of determining the projected problem size with IRR iterations for problem \texttt{grain} with $5\%$ noise.In Figure~\ref{grainrhostepk} $\rho^{(k)}(t)$ for increasing $k$, for $k=1\colon4$ left to right and then above and below.  In Figure~\ref{grainRE} the dashed-dot vertical line corresponds to the location of $\toptp$, the solid line with symbol to $\toptG$ and the solid line to $\toptmin$.}
\end{figure}

Here we only report results obtained using the UPRE as compared to the \textit{optimal} solutions. These results are indicative of IRR implemented with the other methods for finding the regularization parameter. 
To examine the process carefully in one case, we focus on problem \texttt{grain} with $5\%$ noise. Function $\rho(t)$ at the first step $k=0$ does not differ significantly from the case with $10\%$ noise, shown in Figure~\ref{rhotwoD}.  We use a maximum subspace size with $t=100$ for the calculation of  ${\toptG}$.  Figure~\ref{grainrhostepk} shows $\rho^{(k)}(t^{(k-1)},t)$ for the choices of $t$, $\toptp$, $\toptmin$ and  $\toptG$.   It is clear that  $\rho^{(k)}(t^{(k-1)},t)$ is almost independent of $t^{(k-1)}$ for the first steps, but that noise enters for $k=4$. 
This is also reflected in the RE in Figure~\ref{grainRE}, the RE stabilizes for increasing $t$, and decreases for the first three steps of IRR, but increases at step $4$. The impact of the IRR is similar for problem \texttt{satellite}.

The RE for the two simulations, contrasted with MIN are detailed in Table~\ref{tab:errortwodfive}.   The IRR stabilizes the solutions leading to results comparable to those of MIN. Moreover, it is also clear that one may not conclude that finding  $\topt$   using $\toptmin$ is preferable to using  $\toptG$ or $\toptp$. The results in Table~\ref{tab:errortwodfive} indicate that using $\toptp$ leads to best solutions in one case, and $\toptG$ in the other, although effectively  the quality is comparable.  Provided that the solutions are stabilized with the IRR, improvements in the solutions are obtained in a limited number of steps using IRR with relatively small subspaces for the iterative updates. 
\begin{table}[!htb]\begin{center}
\caption{RE for problem \texttt{grain} and \texttt{satellite} with $5\%$ noise corresponding to $\nu=.05$, for solutions found using different selection of $\topt$ as compared to the optimum found with UPRE and overall optimum selected over the range for $\zeta$.\label{tab:errortwodfive}}
\begin{tabular}{|*{5}{c|}} \hline
&\multicolumn{4}{|c|}{\texttt{grain} with ${\toptmin} = 21$ and ${\toptG}=27$} \\ \hline
Iteration& ${\toptmin}$ &$\toptG$ &Min for UPRE&Overall Min\\  \hline
$1$&$0.4248$&    $0.4074 $ &  $0.4030$&   $0.3993$ \\ \hline
$2$&$0.4073  $&    $  0.3962   $&    $ 0.3929  $&    $  0.3903$ \\ \hline
$3$&$0.3899 $&    $   0.3811 $&    $   0.3788$&    $   0.3776$ \\ \hline
$4$&$0.3775 $&    $   0.3731 $&    $   0.3728$&    $   0.3731$ \\ \hline
&\multicolumn{4}{|c|}{\texttt{satellite} with ${\toptp} = 42$ and ${\toptG}=69$}\\ \hline
Iteration& ${\toptp}$ &${\toptG}$ &Min for UPRE&Overall Min\\  \hline
$1$& $0.3566$&    $    0.3520$&    $   0.3517  $&    $  0.3863$ \\ \hline
$2$&$0.3375 $&    $   0.3430  $&    $  0.3373  $&    $  0.3382$ \\ \hline
$3$&$0.3374  $&    $  0.3431 $&    $   0.3372 $&    $   0.3352$ \\ \hline
$4$&$0.3385  $&    $  0.3485 $&    $   0.3371 $&    $   0.3352$ \\ \hline
\end{tabular}
\end{center}\end{table}

\subsubsection{Terminating the IRR iteration}
The graphs of $\rho(t)$  with increasing $k$  in Figure~\ref{grainrhostepk} indicate that  the properties of $\rho^{(k)}(t^{(k-1)},t)$  can be used to determine effective termination of the IRR, based on the iteration $k$ when noise enters into $\rho^{(k)}(t^{(k-1)},t)$. Our experience has shown that the optimal solution in terms of image quality is achieved not at the step before noise enters in $\rho^{(k)}(t^{(k-1)},t)$  but two steps before. 

\subsubsection{Results for $10\%$ noise}
\begin{figure}[!htb]
\begin{center}
\subfigure[RE \texttt{grain} $10\%$ noise \label{grainreten}]{\includegraphics[width=.42\textwidth]{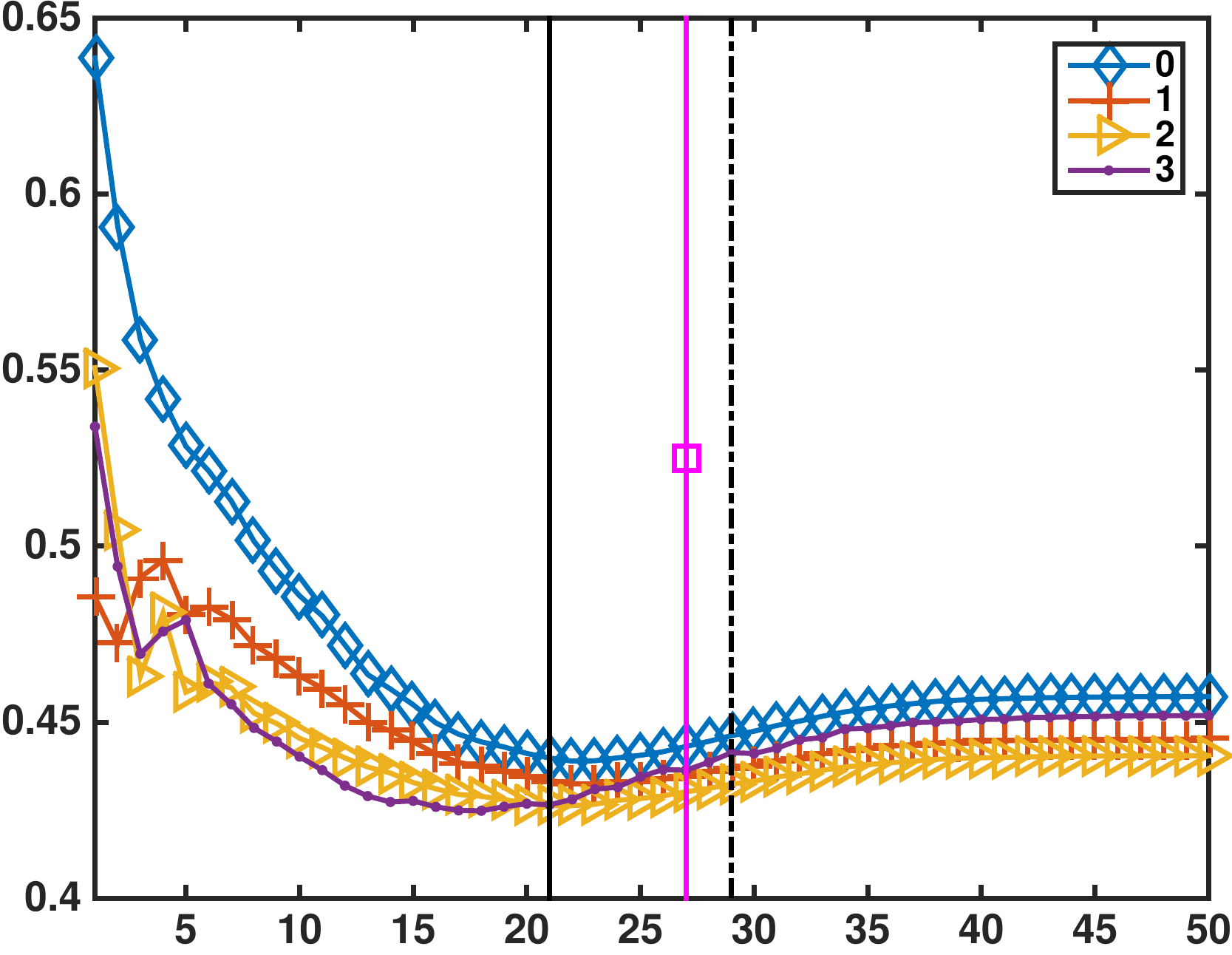}}\hspace{.3cm}
\subfigure[RE \texttt{satellite} $10\%$ noise \label{sattelitere}]{\includegraphics[width=.42\textwidth]{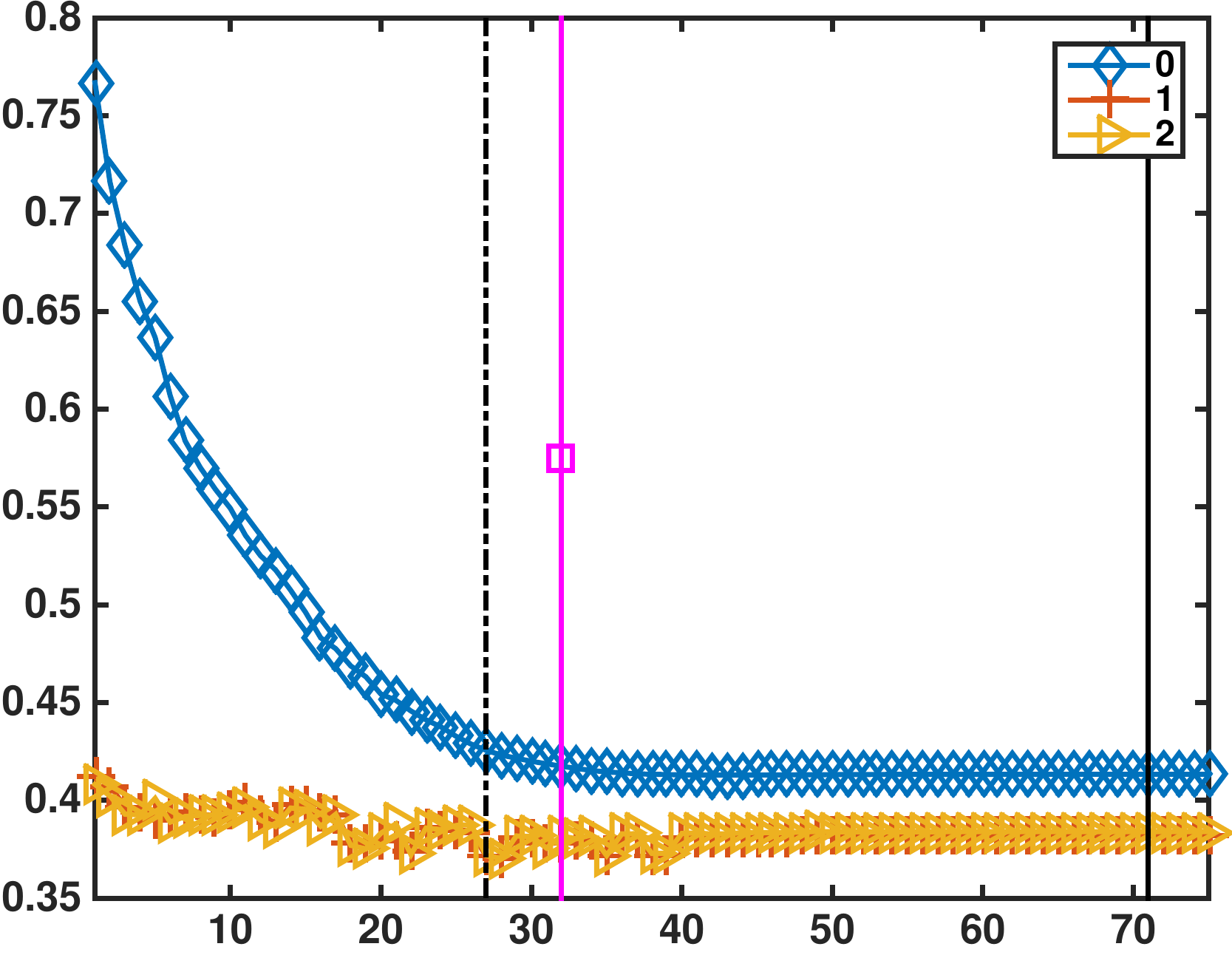}}\hspace{.3cm}
\end{center}
\caption{RE  for problem \texttt{grain}, (\texttt{satellite})  with noise level $10\%$, in Figures~\ref{grainreten} and \ref{sattelitere}, respectively. In each case with increasing iteration  for solutions calculated using UPRE. The solutions are stable out to $k=2$ iterations of IRR. Here the dashed-dot vertical line corresponds to the location of $\toptp$, the solid line with symbol to $\toptG$ and the solid line to $\toptmin$. 
\label{fig2dmeasuresgrain}}
\end{figure}

\begin{figure}[!htb]
\begin{center}
\subfigure[$k=3$, MIN $\topt=19$]{\includegraphics[width=.17\textwidth]{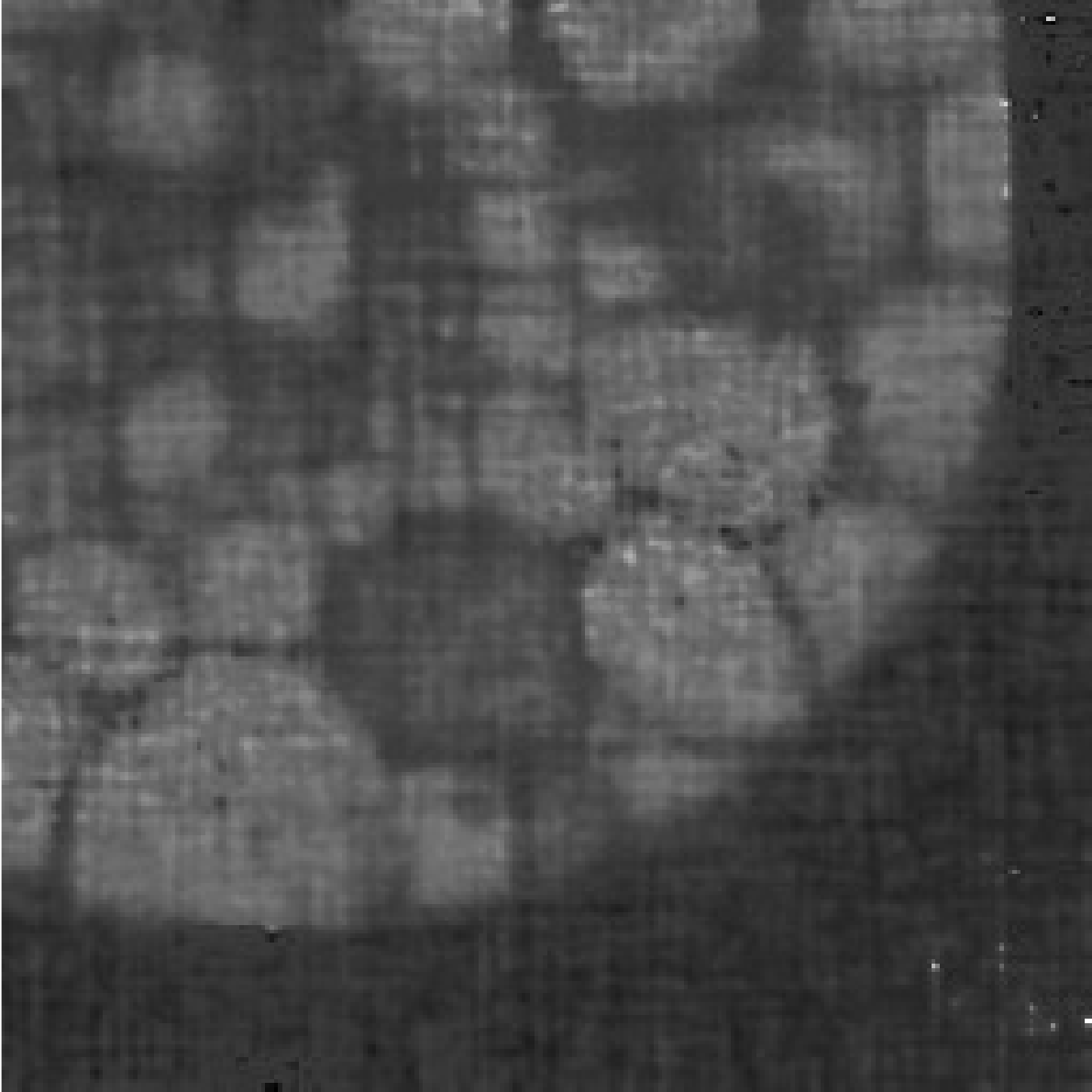}}\hspace{1cm}
\subfigure[$k=2, \toptmin=21$]{\includegraphics[width=.17\textwidth]{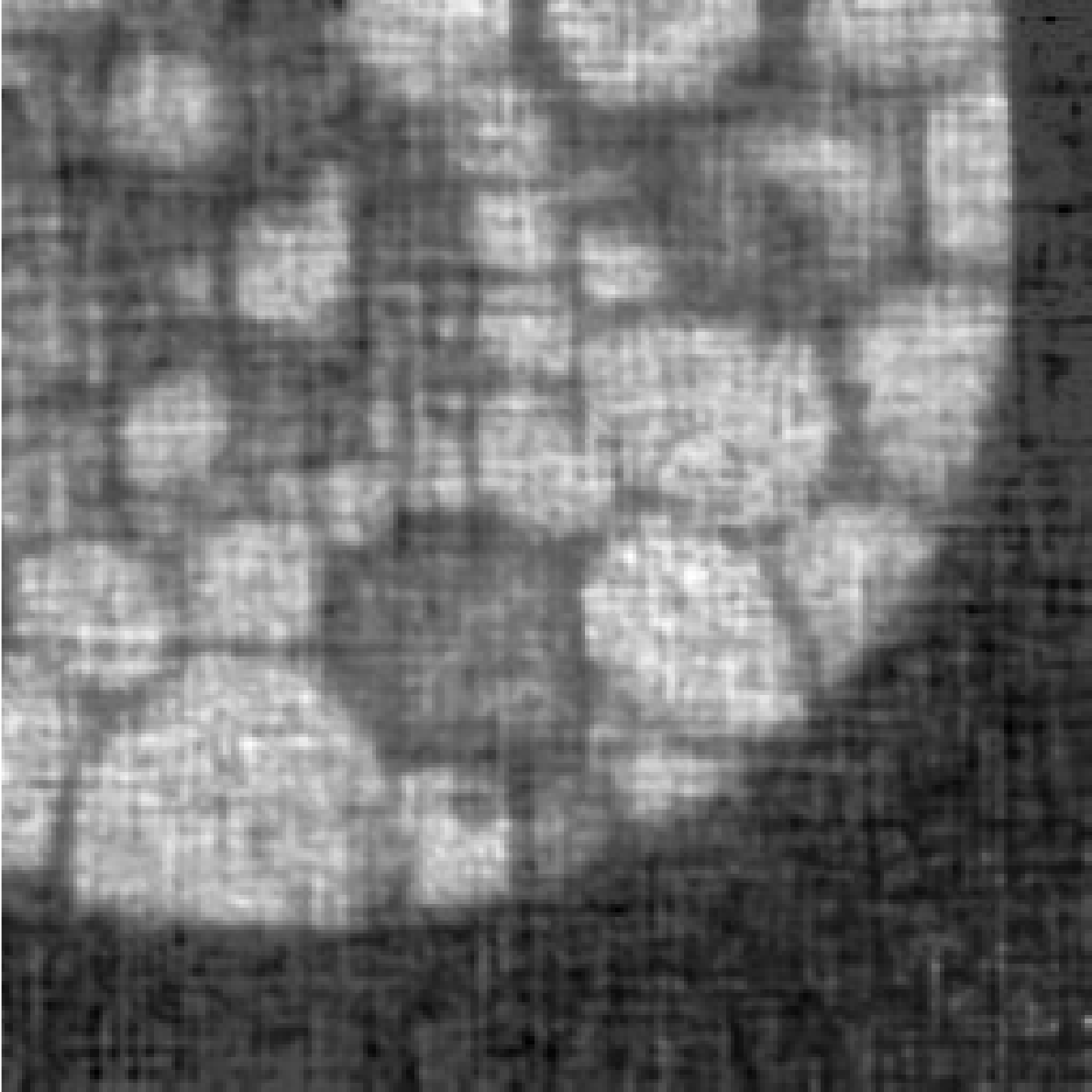}}\hspace{1cm}
\subfigure[$k=2, \toptG=27$]{\includegraphics[width=.17\textwidth]{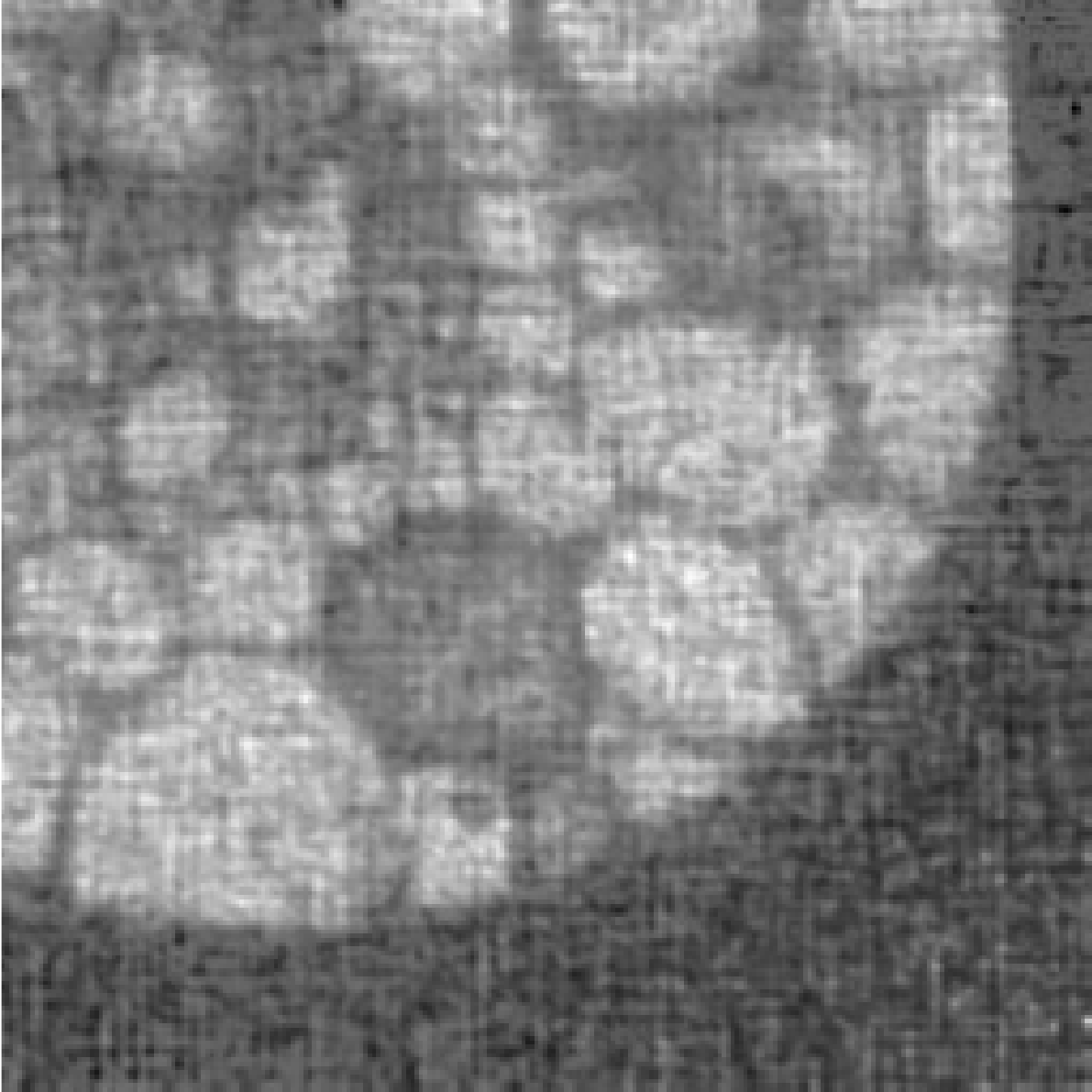}}\hspace{1cm}
\subfigure[$k=2, \toptp=29$]{\includegraphics[width=.17\textwidth]{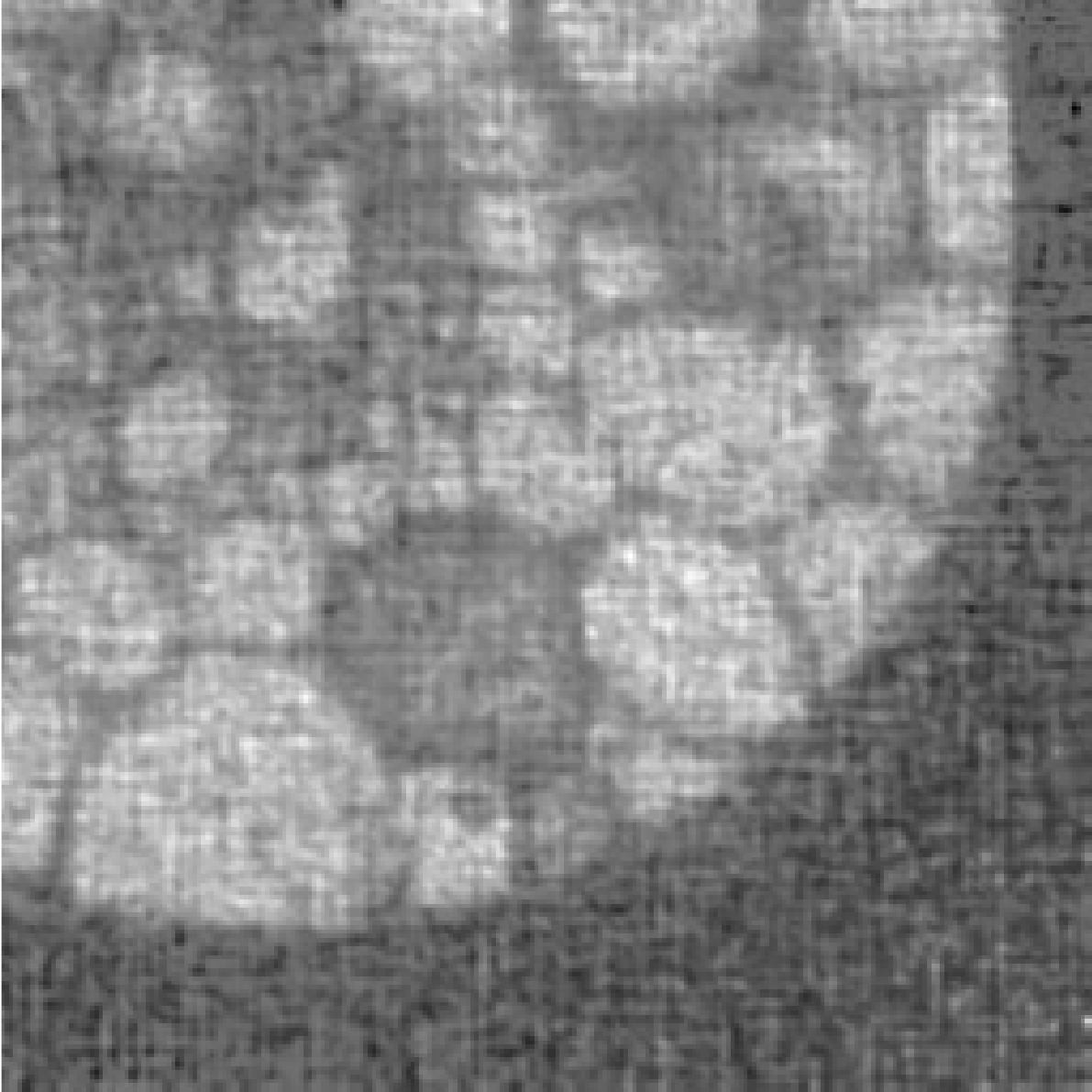}}\\
\end{center}
\caption{Solutions for noise level $10\%$ for   \texttt{grain} using UPRE to find the regularization parameters and comparing the solutions obtained for different $\topt$ as compared to the solution with minimum error, MIN, for IRR at the indicated step applied to the solutions in Figure~\ref{fig2dsolutionsgrain_step1}.   \label{fig2dsolutionsgrain_irrstep}}
\end{figure}

In Figure~\ref{fig2dmeasuresgrain}  we illustrate how the RE changes with the iteration count. 
The vertical lines in Figure~\ref{grainreten} demonstrate that  using $\toptp$ or $\toptG$ makes little difference to the quality of the solution when measured with respect to the RE. Example solutions for \texttt{grain} are given in Figure~\ref{fig2dsolutionsgrain_irrstep}, for contrast with the solutions without IRR shown in Figure~\ref{fig2dsolutionsgrain_step1}.

\subsection{Sparse tomographic reconstruction of a walnut}\label{sec:walnut}
\begin{figure}[!htb]
\begin{center}
\subfigure[$\rho(t)$ \label{figrhowalnut}]{\includegraphics[width=.4\textwidth]{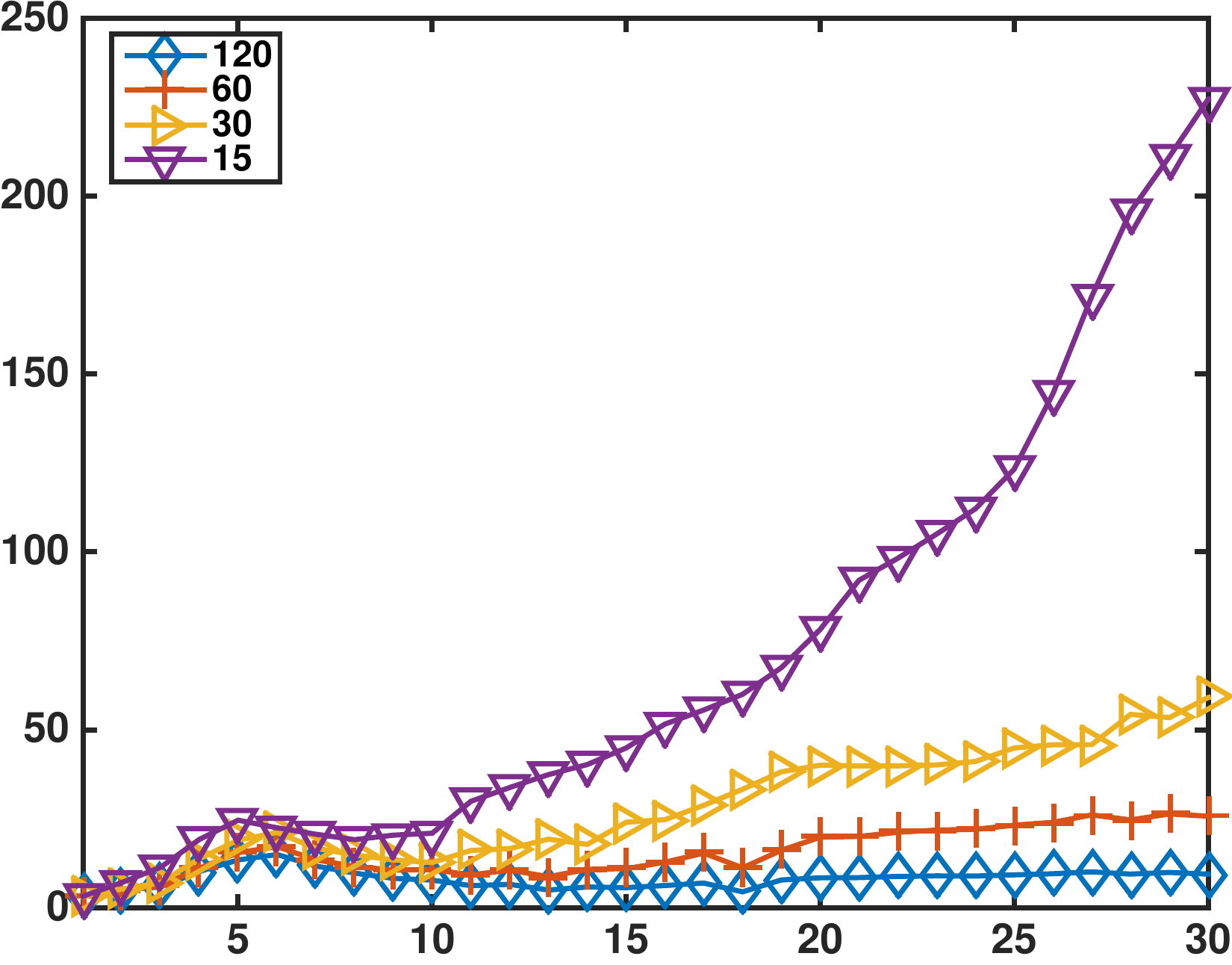}} \hspace{1cm} 
\subfigure[$60$ projections  \label{figspectradN2}]{\includegraphics[width=.4\textwidth]{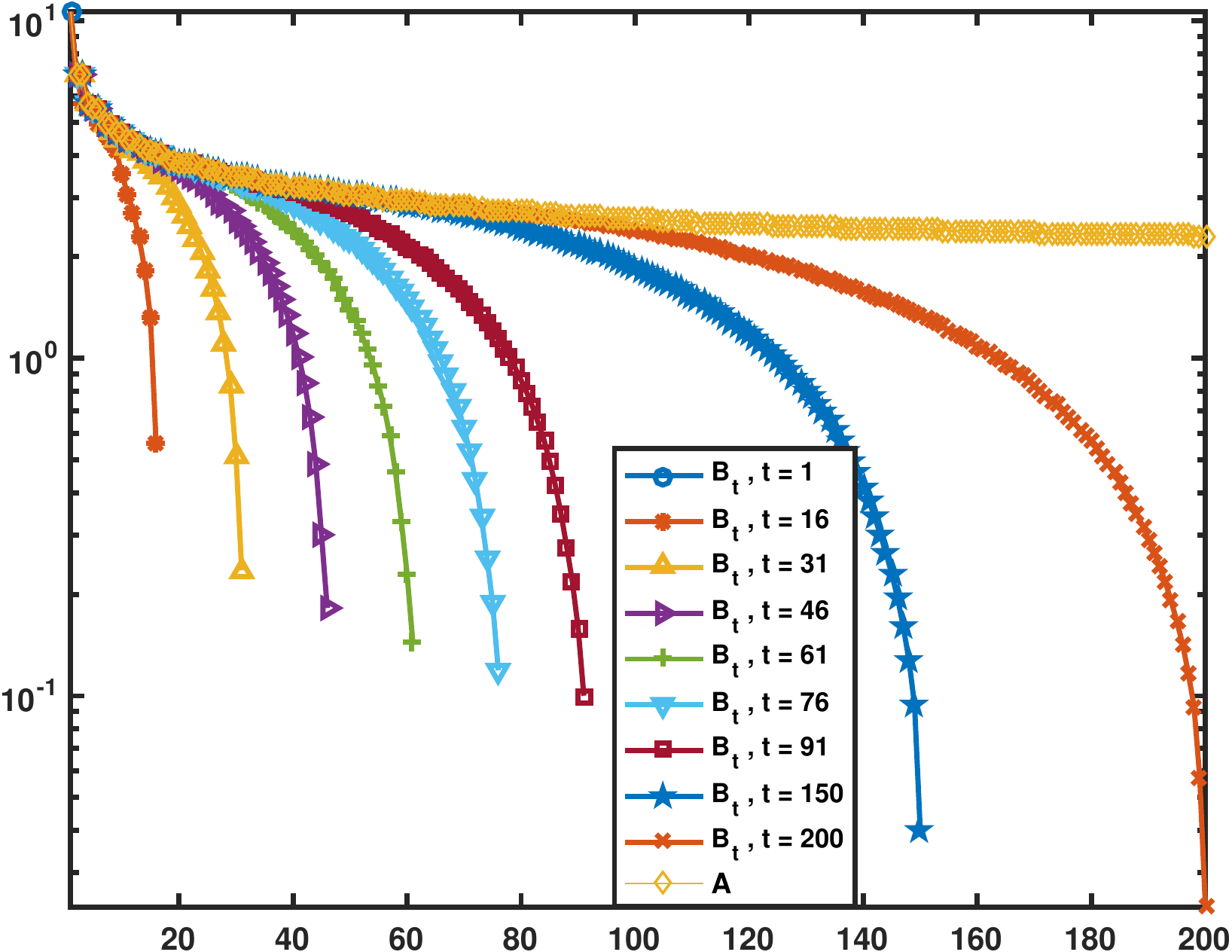}} \\
\subfigure[$30$ projections ]{\includegraphics[width=.4\textwidth]{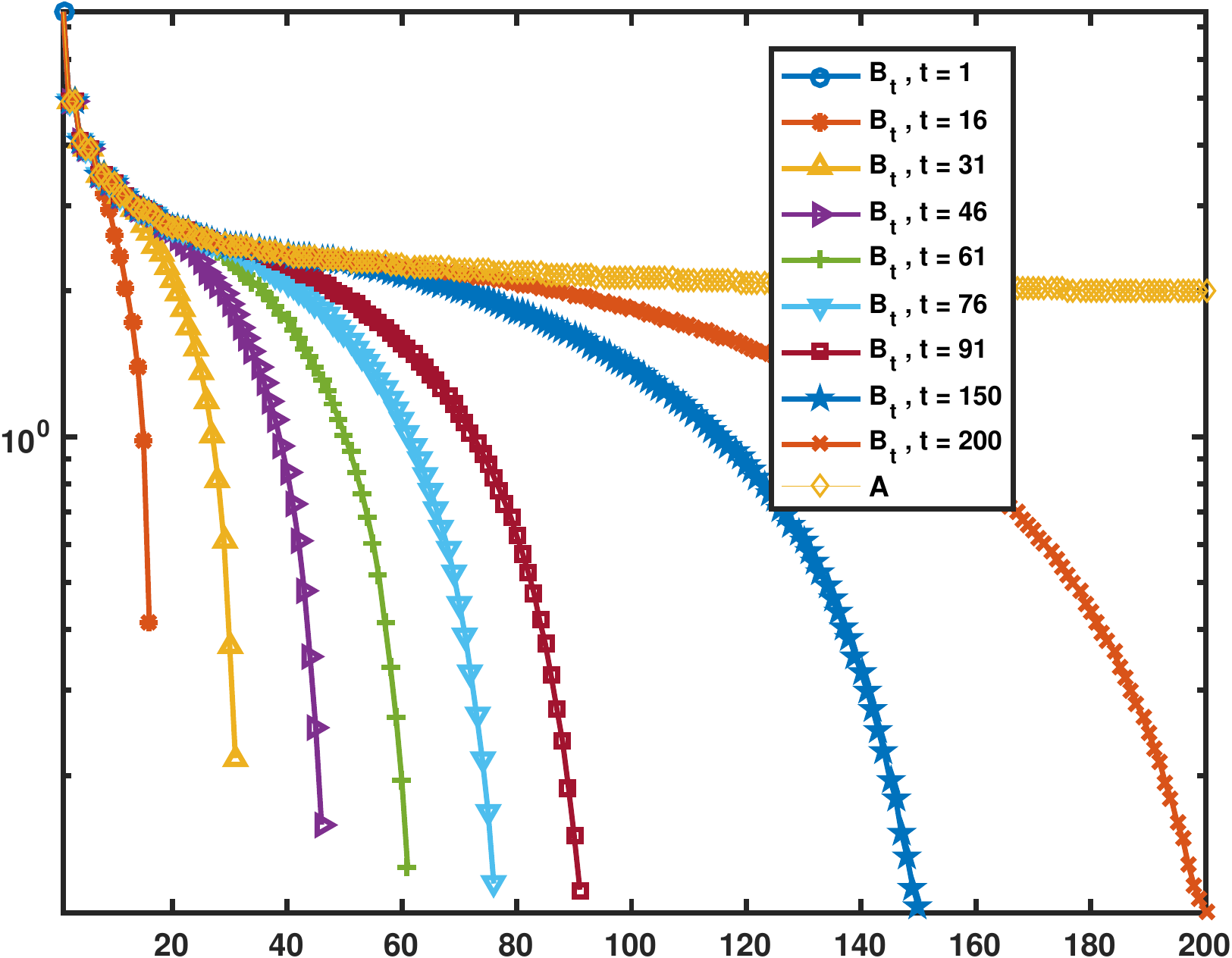}} \hspace{1cm} 
\subfigure[$15$ projections  \label{figspectradN8}]{\includegraphics[width=.4\textwidth]{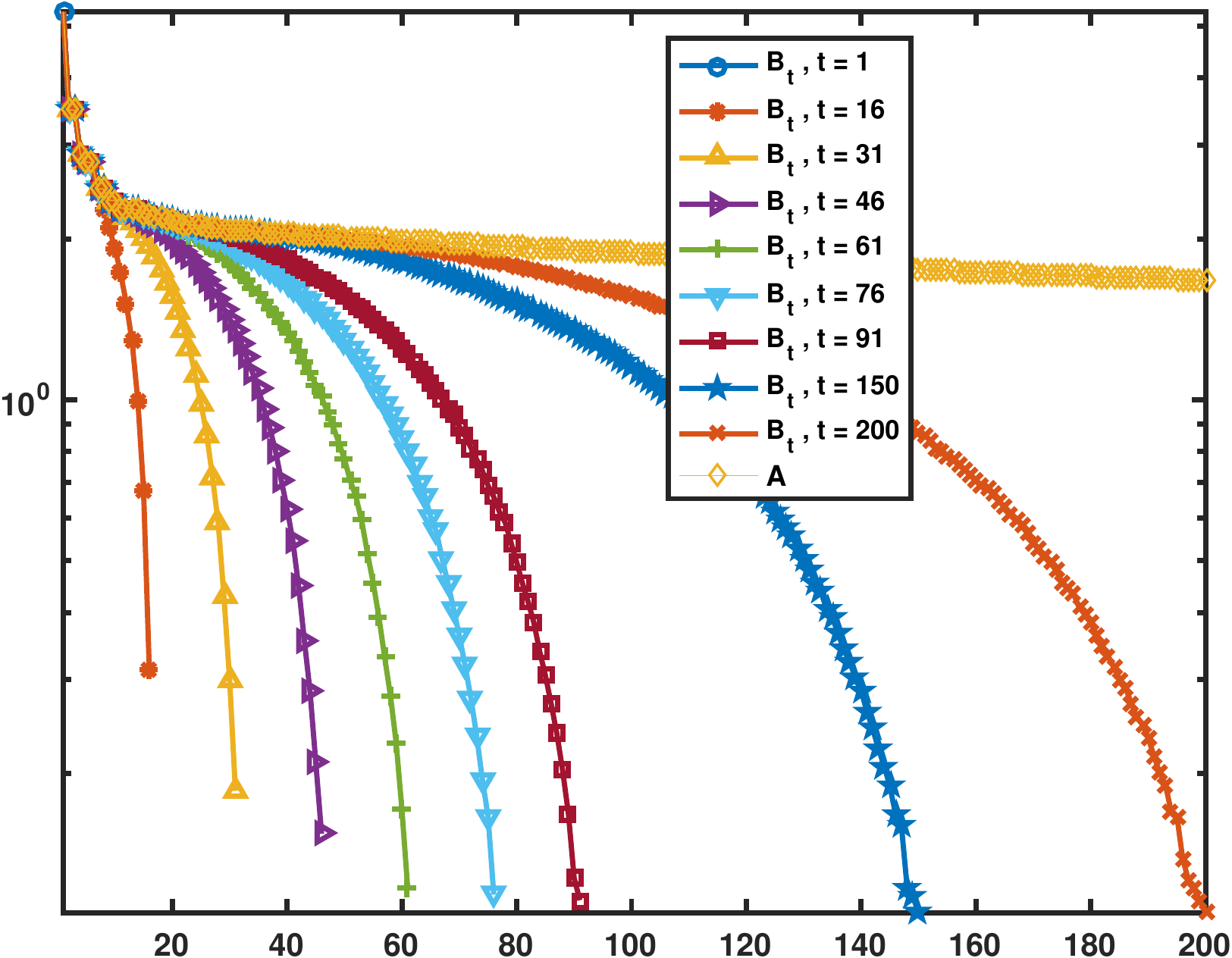}} \\
\end{center}
\caption{\label{walnutphi} $\rho(t)$ for increasing sparsity for resolution $164$ for the walnut data in Figure~\ref{figrhowalnut} and the spectra for increasing $t$, $t=1:15:91$, $150$, $200$, with increasing sparsity in Figures~\ref{figspectradN2}-\ref{figspectradN8}.}
\end{figure}

\begin{figure}[!htb]
\begin{center}
\subfigure[UPRE:      $k=0$]{\includegraphics[width=.24\textwidth]{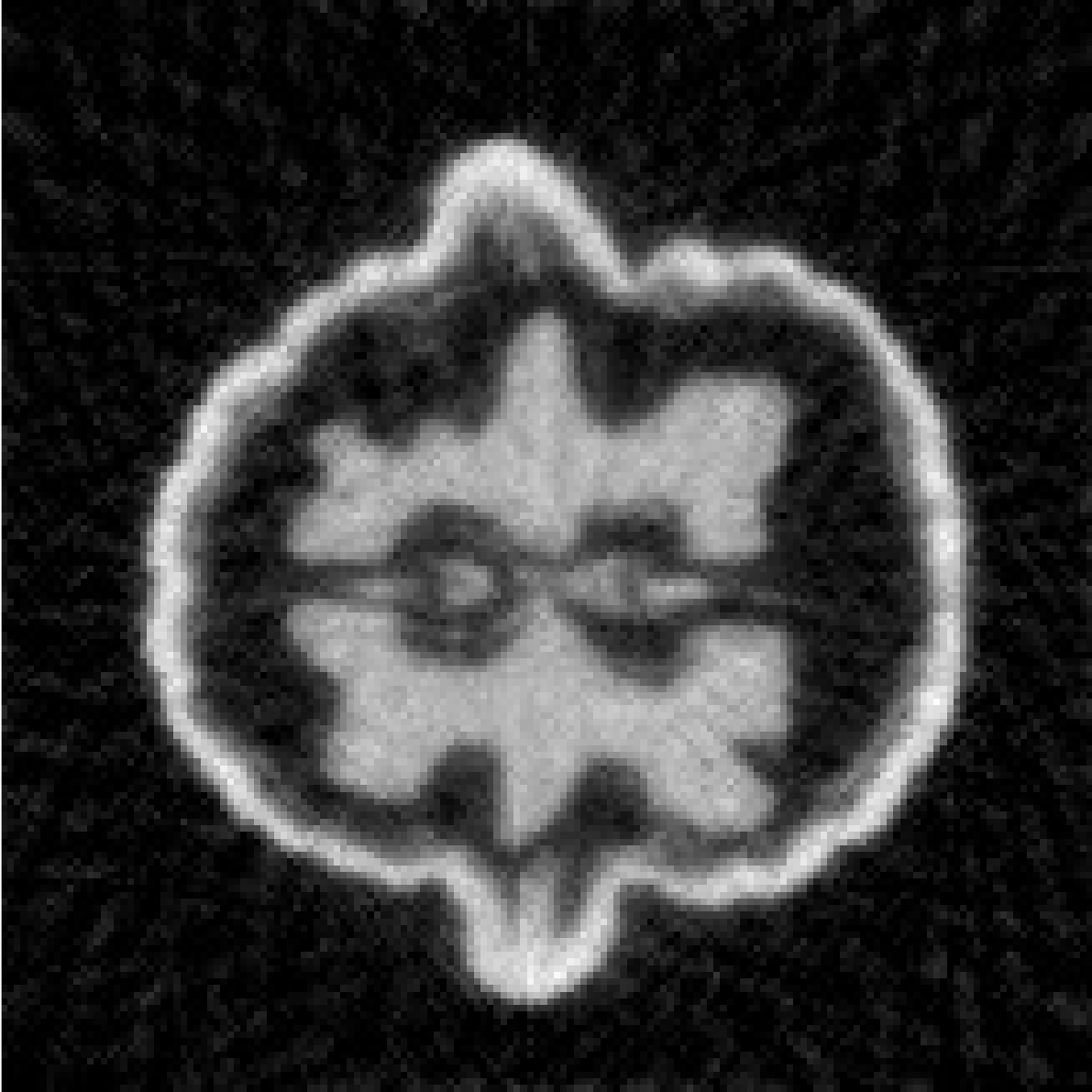}}\hspace{1.1cm}
\subfigure[UPRE:     $k=1$]{\includegraphics[width=.24\textwidth]{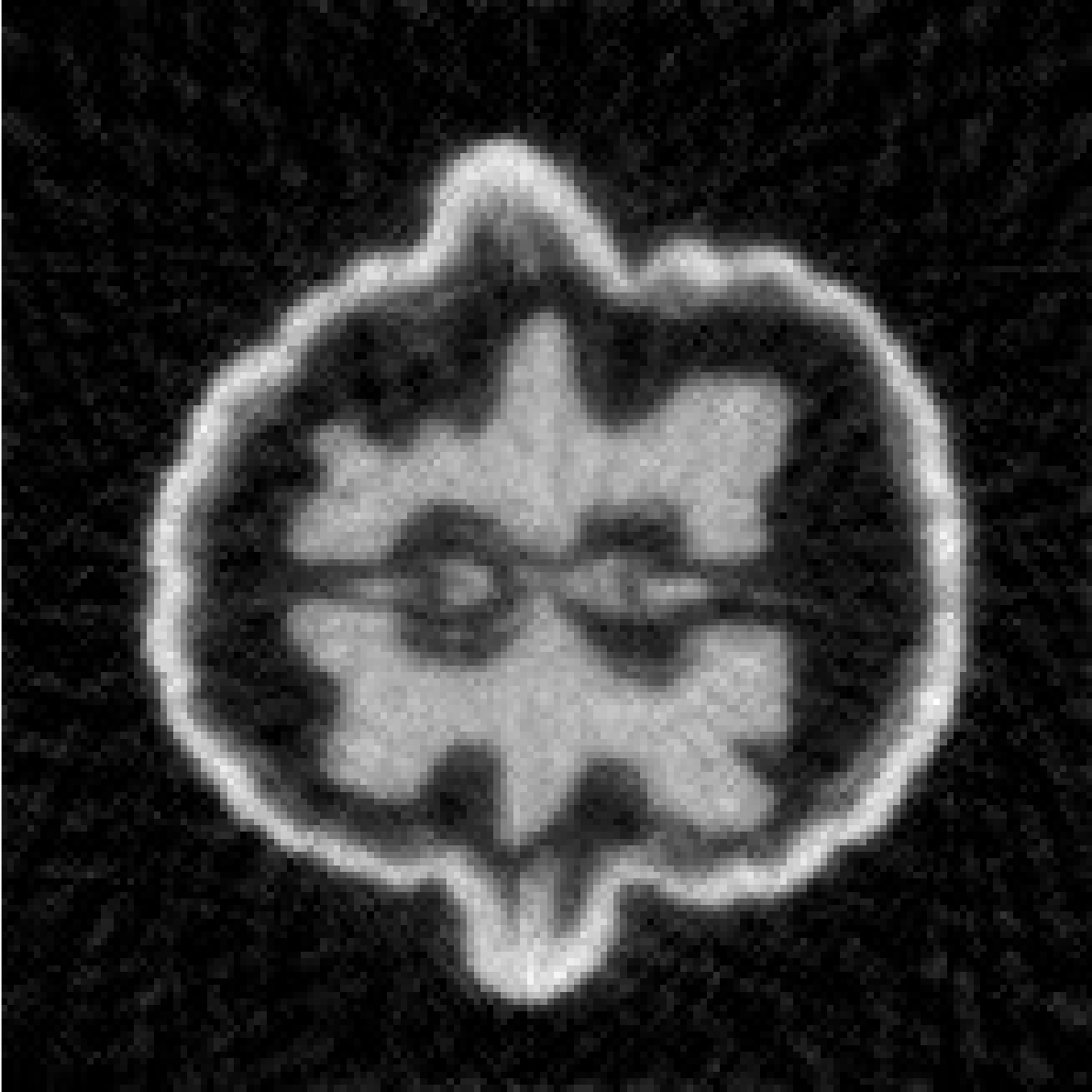}}\hspace{1.1cm}
\subfigure[UPRE:      $k=2$]{\includegraphics[width=.24\textwidth]{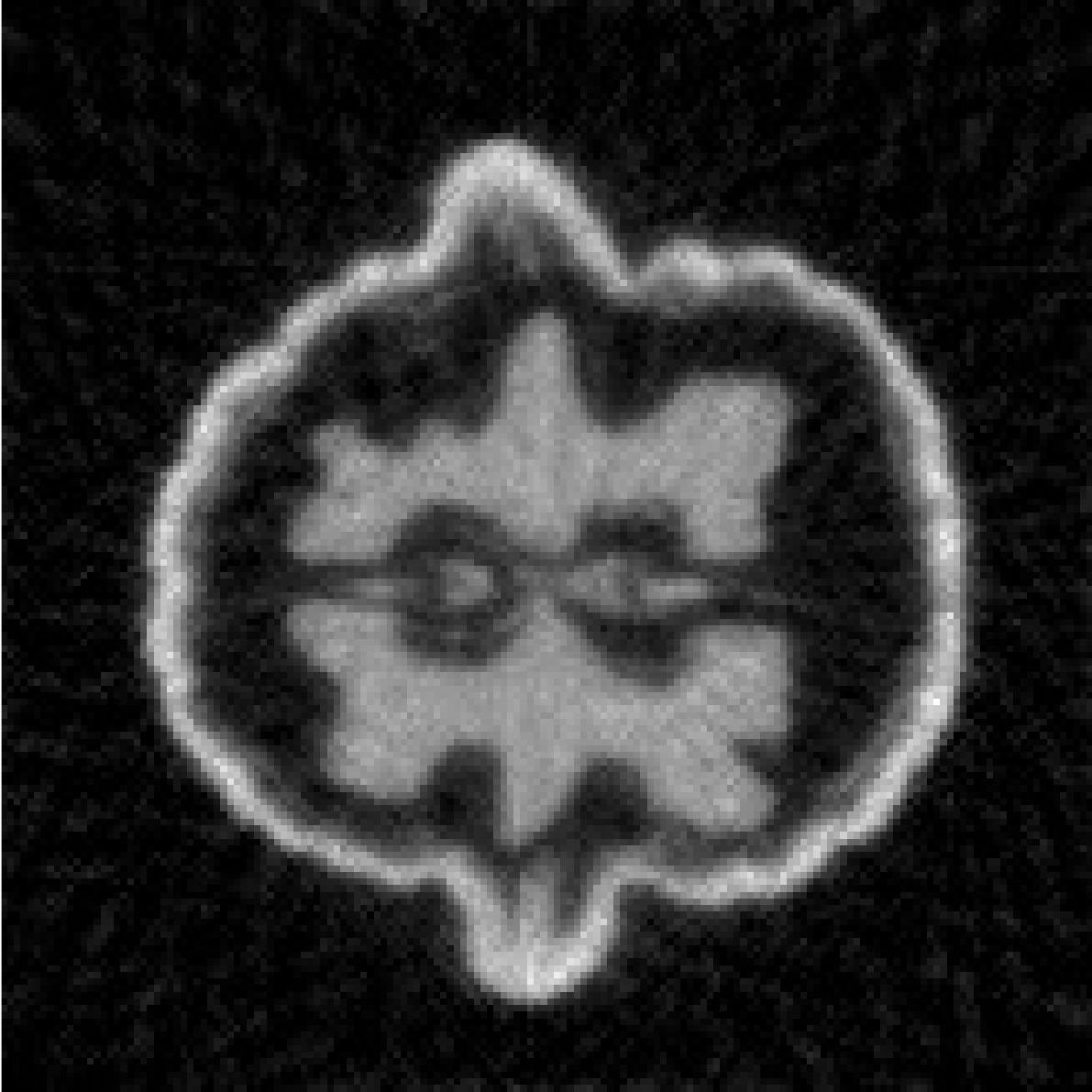}}\\
\end{center}
\caption{Solutions at increasing iterations for walnut with resolution $164\times 164$,  $\tmin=5$, $\toptp=8$, positivity constraint  and sampling at $6^\circ$ intervals, $60$ projections. \label{walnut16460}}
\end{figure}

\begin{figure}[!htb]
\begin{center}
\subfigure[UPRE:      $k=0$]{\includegraphics[width=.24\textwidth]{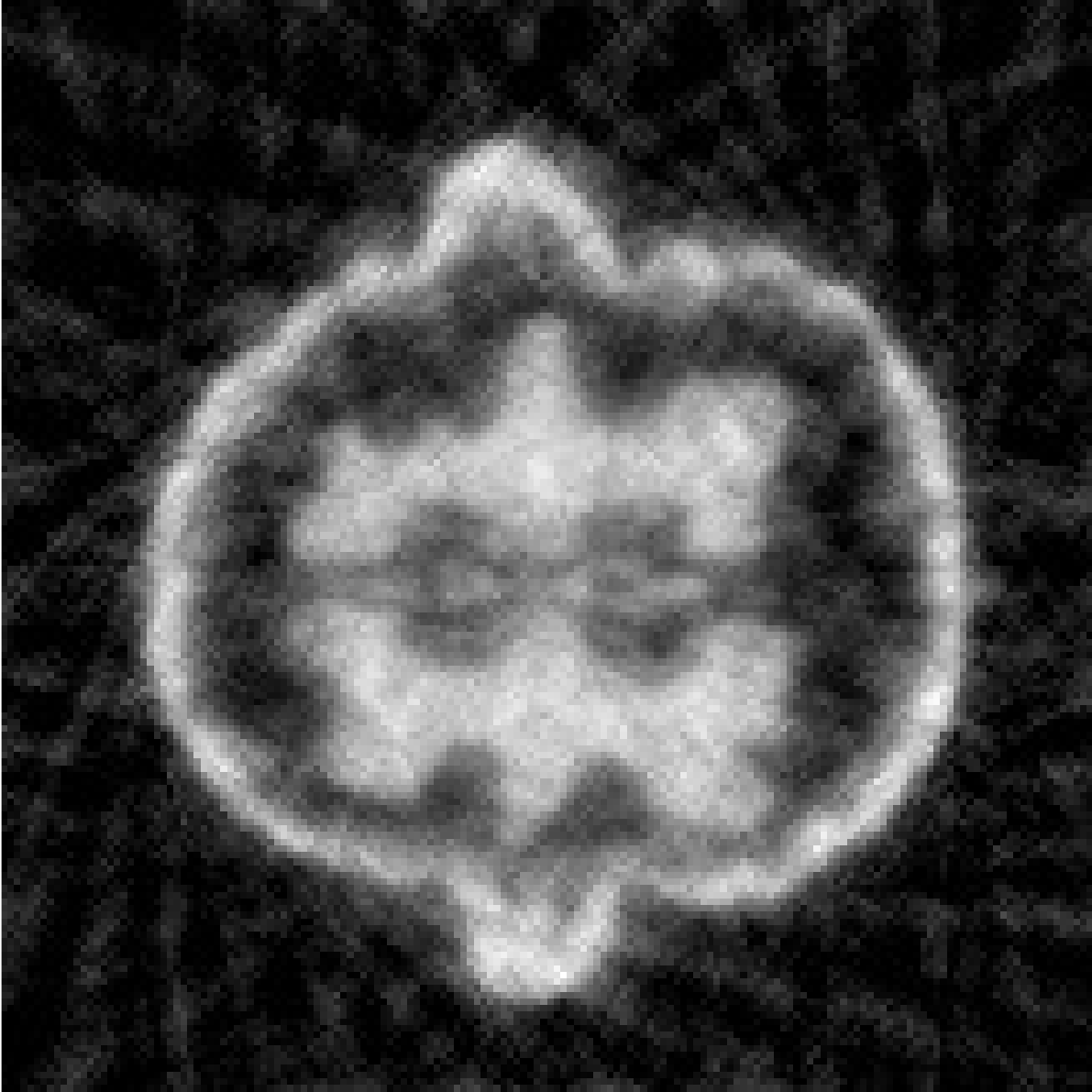}}\hspace{1.1cm}
\subfigure[UPRE:     $k=1$]{\includegraphics[width=.24\textwidth]{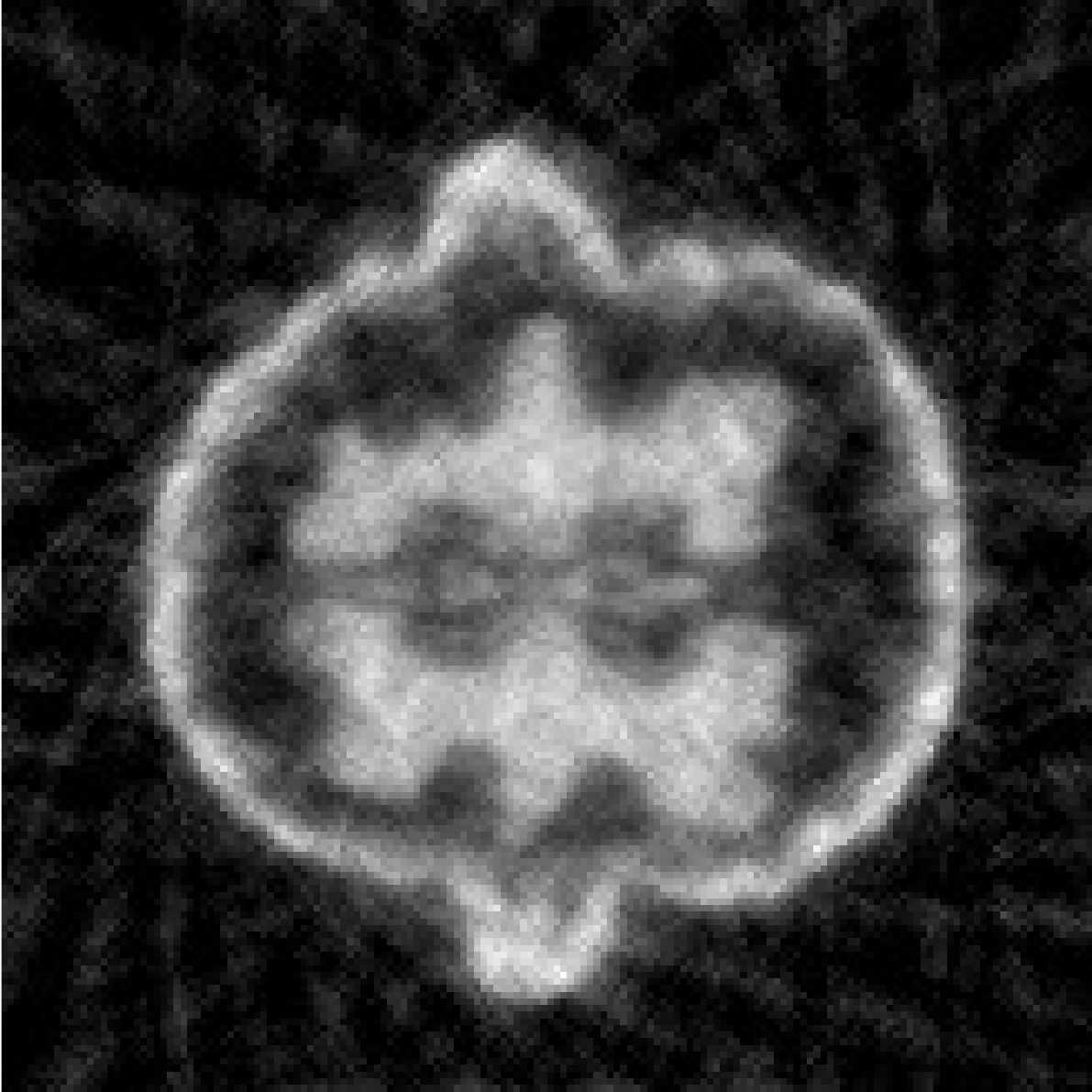}}\hspace{1.1cm}
\subfigure[UPRE:      $k=2$]{\includegraphics[width=.24\textwidth]{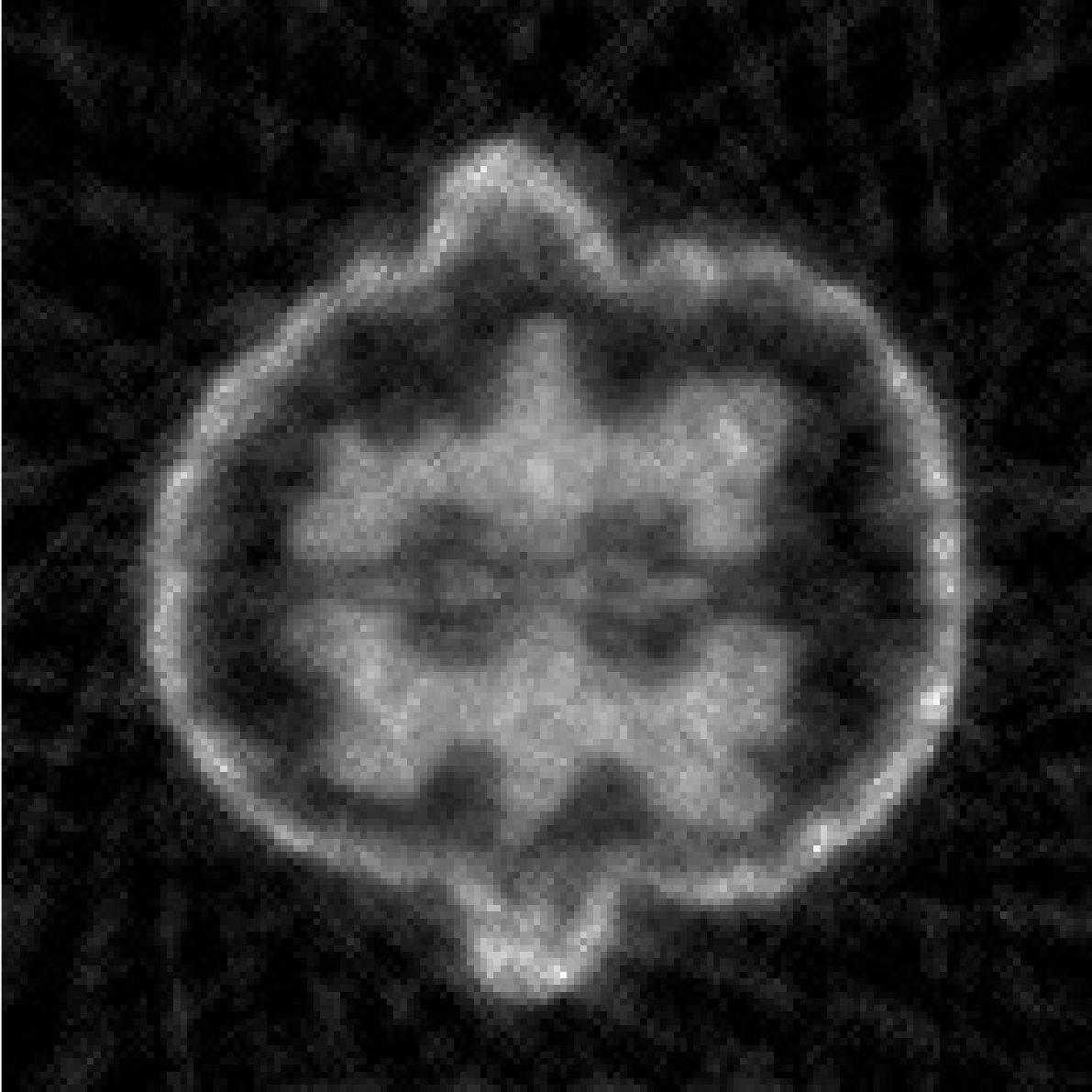}}\\
\end{center}
\caption{Solutions at increasing iterations for walnut with resolution $164\times 164$,  $\tmin=5$, $\toptp=8$, positivity constraint  and sampling at $12^\circ$ intervals, $30$ projections. \label{walnut16430}}
\end{figure}

To contrast the success of the regularization parameter estimation techniques in the context of a $2$D projection problem, we present results for the reconstruction of projection data obtained from tomographic x ray data of a walnut, used for edge preserving reconstruction in \cite{walnut2} and with the description of the data described in \cite{walnut1}. Datasets \texttt{DataN} correspond to resolution $N \times N$ in the image, and use $120$ projections, corresponding to $3^\circ$ sampling. Data are provided with $N=82$, $164$ and $328$. We use  resolution $164$ with $120$ projections, and then downsampled to $60$,  $30$ and $15$ projections,  i.e. angles $3^\circ$, $6^\circ$, $12^\circ$ and $24^\circ$.  Results with resolutions $82$ and $328$ are comparable. Results are presented using the solution at $\toptp$  regularized using UPRE, Figures~\ref{walnut16460}-\ref{walnut16430}. Results for $120$ projections are almost perfect due to the apparent limited noise in the provided data, while the results with $15$ projections clearly show the projection data.  In all the presented results the parameters $\toptp$ are determined automatically, after manually picking $\tmin=5$  from manual  consideration of the plot for $\rho(t)$, see Figure~\ref{walnutphi}.  Further, from Figures~\ref{figspectradN2}-\ref{figspectradN8}, it is immediate that in this case the LSQR iteration  only offers a partial regularization independent of the sparsity and thus again it is necessary to identify the window for the usable spectrum. All  parameters are then estimated in the same way as for the image restoration cases.

\begin{figure}[!htb]
\begin{center}
\subfigure[UPRE: $60$]{\includegraphics[width=.24\textwidth]{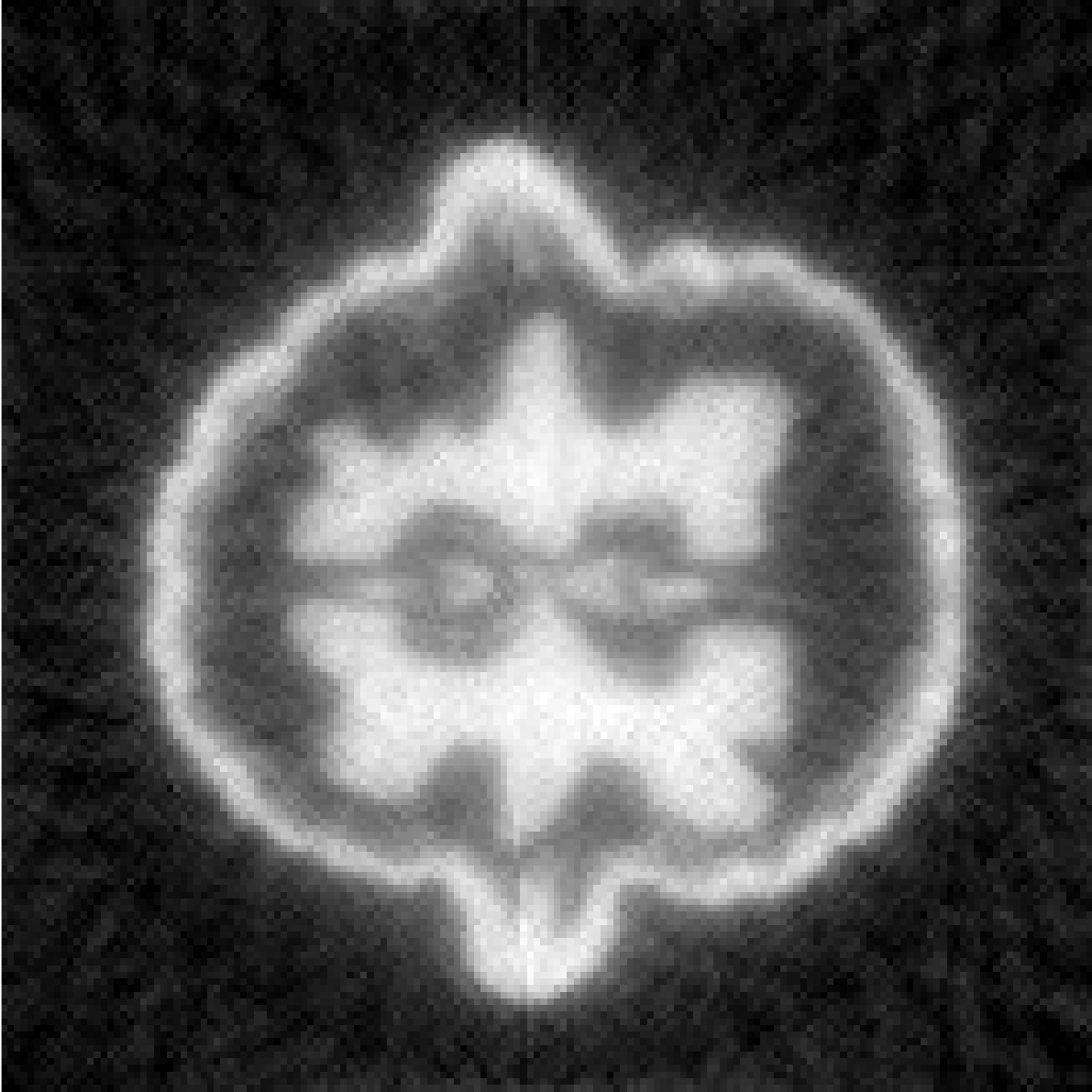}}\hspace{1.1cm}
\subfigure[UPRE: $30$]{\includegraphics[width=.24\textwidth]{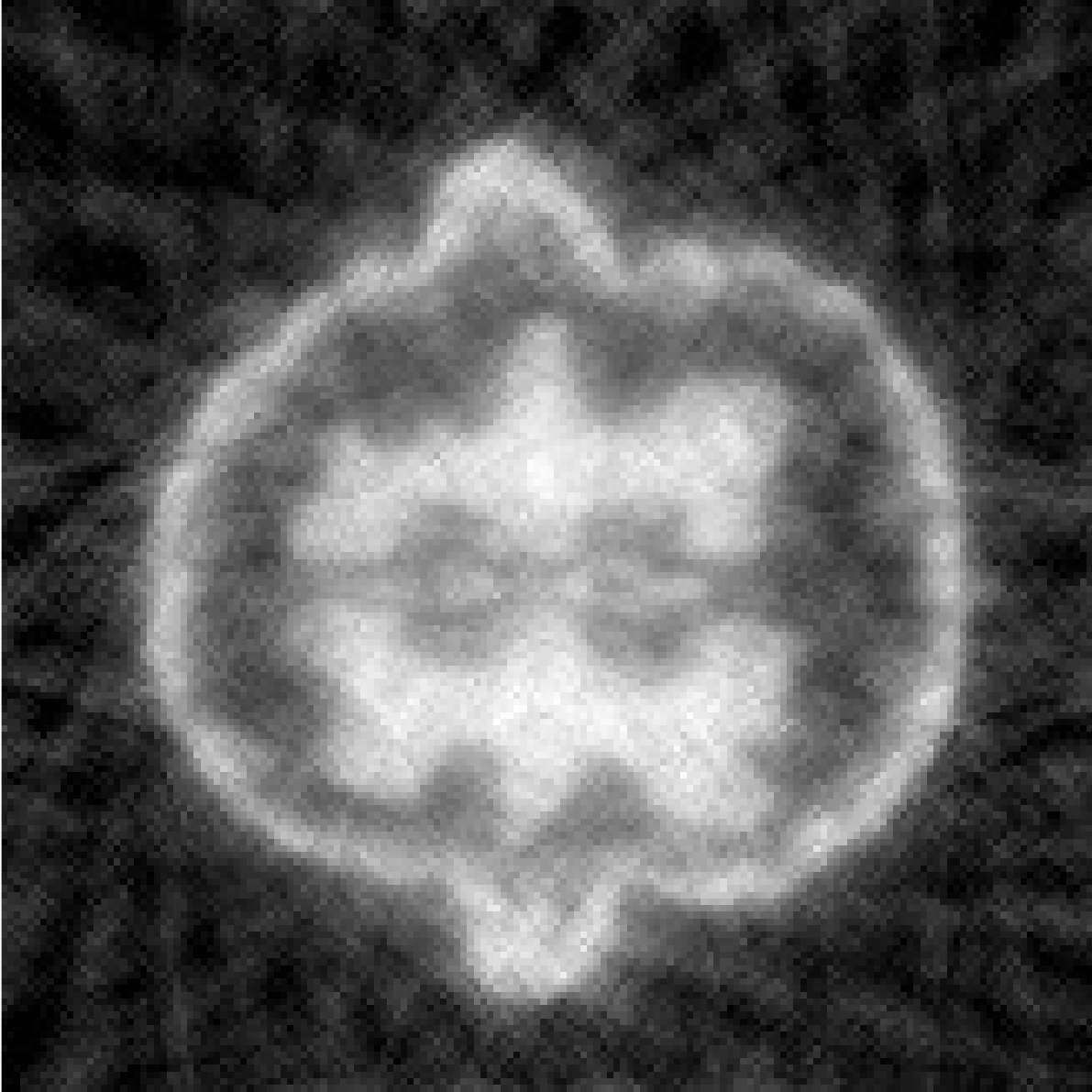}}\hspace{1.1cm}
\subfigure[UPRE: $15$]{\includegraphics[width=.24\textwidth]{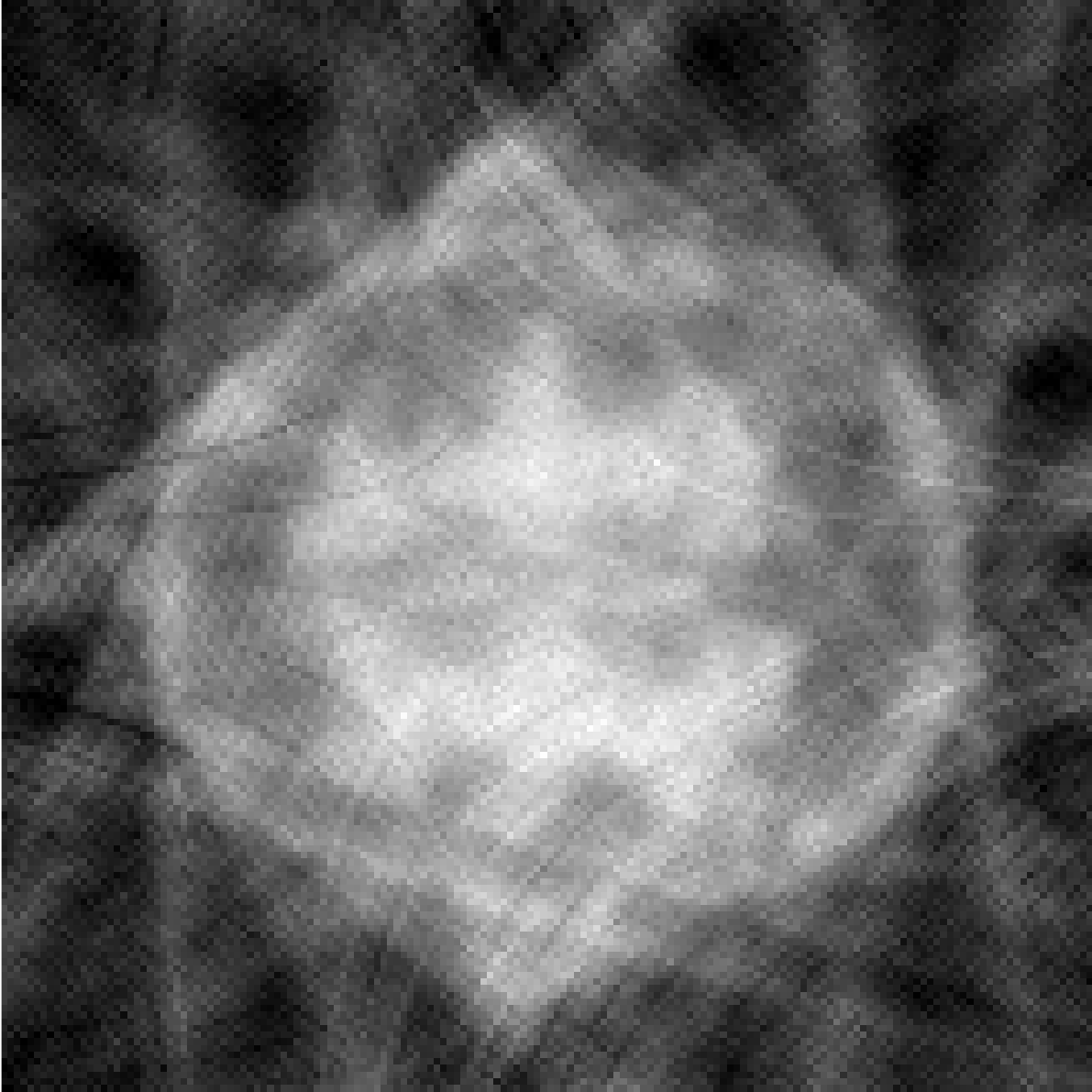}}\hspace{1.1cm}
\end{center}
\caption{Solutions for walnut with resolution $164\times 164$ without IRR, i.e. at $k=0$, with $\tmin=18$ yielding  $\toptp=25$, $29$, and $40$, automatically determined for sampling intervals  $6^\circ$,   $12^\circ$ and $24^\circ$.  \label{walnut164tmin18}}
\end{figure}

Figures~\ref{walnut16460}-\ref{walnut16430} show results for one set of data at increasing sparsity, compare with \cite[Figure 6.6 and Figure 6.7]{walnut2}, which give results with resolution for $N=128$ and $256$, respectively, and angle separation $2^\circ$, $4^\circ$, $6^\circ$ and $12^\circ$. The approach in \cite{walnut2}  uses selected choices for the regularization parameter based on a sparsity argument with prior information and seeks to support the use of the sparsity argument for reconstruction of sparse data sets. The images exhibit the rather standard total variation  blocky structures when applied for truly sparsely sampled data. Our results show robust reconstructions with the automatically determined solutions, after first examining the plot for $\rho(t)$. IRR generates marginal improvements in qualitative solutions. To show the impact of the correct choice of $\topt$ on the solution we show a set of results at iteration $k=0$ using $\tmin=18$ in Figure~\ref{walnut164tmin18}, with the positive constraint. The UPRE  yields solutions qualitatively similar to the case with $\tmin=5$.  In the examples here we do impose an additional positivity constraint on the solutions at each step, before calculating the iterative weighting matrix. 

The results demonstrate that the projected problem with automatic determination of $\zetaopt$ can be used to reconstruct sparsely sampled tomographic data, provided that an initial estimate for $\tmin$ is manually determined by consideration of the plot of $\rho(t)$.  Further, IRR stabilizes the solutions. For the sparse data sets the solutions do not exhibit the characteristic blocky reconstructions of total variation image reconstructions, as seen in \cite{walnut2}. 

\section{Conclusions}\label{conclusions}
We have demonstrated that regularization parameter estimation by the method of UPRE can be effectively applied for regularizing the projected problem. Our results also motivate a choice for the weighting parameter in the WGCV. These results apply for the concept of \textit{full regularization} which was recently introduced in \cite{HJ:16}. It was argued there, however, that in the case of full regularization, no additional regularization of the iterate is required, because the LSQR iterate approximates the TSVD solution of the full problem. It is known, however, that even with truncation additional filtering through the use of a FTSVD solution is required. Moreover, while the projected solution may not actually need regularization, without prior knowledge of the optimal subspace size $t$, it is difficult to determine the point at which semi-convergence contaminates the projected solution. Hence effectively regularizing via the hybrid LSQR is still necessary. Our results demonstrate that the regularization estimators will find effective regularization of the FTSVD solution of the full problem. In the case of the \textit{partial regularization}, i.e. in the cases in which a small Ritz value appears before the LSQR iterate has captured the dominant SVD components of $A$, as discussed in \cite{HJ:16},  regularization on the projected problem will not adequately regularize the full problem. Here we handled the situation in which LSQR does not sufficiently capture the dominant singular space by restricting the range for the regularization parameter dependent on the singular value for $\gamma_{\topt}$, so as to effectively use a FTSVD solution of the projected problem.  For future investigation we suggest that  it is important to identify the extent to which the LSQR iterate captures a right singular subspace for $A$, and to then potentially use more than one regularization parameter, one which is chosen to regularize the dominant terms of the spectrum, and one which handles the small singular values of $B_t$,  extending the windowed regularization parameter techniques \cite{ChEaOl:11}. 

This work  also demonstrates that  edge preserving regularization, via  iterative reweighting, can  be applied to stabilize regularized solutions of the projected problem. Our results suggest manual estimation of a minimal subspace size can then lead to useful estimates for an optimal projected space, with the use of the IRR leading to improvements in the solutions when $\topt$ is found by different methods, including the use of $\toptp$, $\toptmin$ and $\toptG$, hence making the determination of this $\topt$ less crucial in providing an acceptable solution. Future work on this topic should include extending 
use of more general iteratively  reweighted regularizers accounting for edges in more than one direction in conjunction with the projected solutions. 

\appendix
\section{Expansion Solutions}\label{SVD}
Suppose the SVD of matrix $A$, 
$A \in\mathcal{R}^{m \times n} $, is given by   $A =U\Sigma V^T$,  where  the singular values are ordered $\sigma_1\ge \sigma_2 \ge \cdots \ge \sigma_{m^*}>0$ and occur on the diagonal of $\Sigma \in \mathcal{R}^{m\times n}$ with $n-m$ zero columns (when $m < n$) or $m-n$ zero rows (when $m > n$),  and $ U \in\mathcal{R}^{m \times m}$, and $ V \in\mathcal{R}^{n \times n}$ are orthogonal matrices, \cite{GoLo:96}.     Then 
\begin{align}\label{svdsoln}
\bfx(\alpha) &= \sum_{i=1}^{m^*} \frac{\sigma^2_i}{\sigma_i^2+\alpha^2} \frac{\bfu^T_{i}\bfb }{\sigma_i} \bfv_{i} = \sum_{i=1}^{m^*} \phi_i(\alpha) \frac{\hat{b}_i}{\sigma_i} \bfv_{i}, \quad \hat{b}_i=\bfu^T_{i}\bfb. \end{align}
 For the projected case $ B_t \in\mathcal{R}^{(t+1) \times t} $,  i.e. $m > n$,  and the expression still applies with  $\|\bfb\|_2\bfe^{(t+1)}_1$ replacing $\bfb$,  $\zeta$ replacing $\alpha$, $\gamma_i$ replacing $\sigma_i$ and  $m^*=t$ in \eqref{svdsoln}.
\section{Regularization Parameter Estimation}\label{appB}
All formulae apply using the SVD for  $B_t$ replacing that for matrix $A$. 
\subsection{Unbiased Predictive Risk Estimator}
The UPRE function   is given by
\begin{align*}
U(\alpha)=\sum_{i=1}^{{m^*}} \left( \frac{1}{\sigma_{i}^{2}\alpha^{-2}+1}\right)^2 \hat{b}_{i}^{2}+2\left(\sum_{i=1}^{m^*}\phi_{i}\right)-m. 
\end{align*}
\subsection{Morozov Discrepancy Principle}
The MDP function is given by
\begin{align*}
\sum_{i= 1}^{m^*} \left( \frac{1}{\sigma_{i}^{2}\alpha^{-2}+1}\right)^2 \hat{b}_{i}^{2}+\sum_{i=n+1}^{m}\hat{b}_{i}^{2} = \delta. 
\end{align*}
 For the projected case   $ \delta_{\mathrm{proj}}$ replaces $\delta$. 
\subsection{Generalized Cross validation}
Using the SVD for $B_t$ the  WGCV function  is given by 
\begin{align*}
G(\zeta,\omega)= \frac{\sum_{i=1}^{t} \left( \frac{1}{\gamma_{i}^{2}\zeta^{-2}+1}\right)^2 \hat{b}_{i}^{2}+\sum_{i=t+1}^{t+1}\hat{b}_{i}^{2} }{((1 +t   -\omega t ) +\omega \zeta^2\sum_{i=1}^t\frac{1}{\gamma^2_i+\zeta^2})^2}.
\end{align*}
With $\omega = 1$ this reduces to the expression for the projected GCV, \eqref{projgcv}.

\end{document}